\documentclass[12pt,a4paper]{article}

\usepackage{bbm}
\usepackage{amsmath}
\usepackage[english]{babel}
\usepackage{amstext}
\usepackage{amsfonts}
\usepackage{subfigure}
\usepackage{multirow}
\usepackage{microtype}% makes pdf look better
\usepackage{graphicx}
\usepackage{float}
\usepackage{color}
\usepackage{epstopdf}
\usepackage[english]{babel}
\usepackage{amsthm,accents}
\usepackage{amsmath}
\usepackage{mathtext}
\usepackage{amstext}
\usepackage{amsfonts}
\usepackage{graphicx}
\usepackage{caption}
\usepackage{wrapfig}
\usepackage{amssymb}
\usepackage{thmtools}
\usepackage{mathtools}
\usepackage{breqn}
\usepackage{indentfirst}
\usepackage{longtable}
\usepackage{changepage}
\usepackage{multirow}
\usepackage{url}
\usepackage{subfig}
\usepackage{sidecap}
\usepackage{float}
\usepackage{color}
\usepackage{booktabs}
\usepackage{multirow}

\usepackage{epsfig}
\usepackage{calrsfs}
\usepackage{relsize}
\usepackage{amssymb}

\renewcommand{\restriction}{\mathord{\upharpoonright}}

\newcommand{\vertiii}[1]{{\left\vert\kern-0.25ex\left\vert\kern-0.25ex\left\vert #1\right\vert\kern-0.25ex\right\vert\kern-0.25ex\right\vert}}

\newtheorem{proposition}{Proposition}[section]

 \newtheorem{remark}{Remark}[section]

\newtheorem {definition}{Definition}[section]
\newcommand\abs[1]{\left|#1\right|}
\let\div\undefined
\DeclareMathOperator{\div}{div}

\DeclareMathAlphabet{\pazocal}{OMS}{zplm}{m}{n}

\numberwithin{equation}{section}

\newcommand{\hs}{\hspace{6mm}}

\DeclareMathOperator{\arsinh}{arsinh}

\begin{document}
\title{Reliable computer simulation methods for electrostatic biomolecular models
based on the Poisson-Boltzmann equation}

\author{J. Kraus\footnote{University of Duisburg-Essen, Germany},
S. Nakov\footnote{RICAM, Austrian Academy of Sciences},
S. Repin\footnote{University of Jyv\"askyl\"a, Finland,
V.A. Steklov Institute of Mathematics at St. Petersburg}}
\date{}
\maketitle
\thispagestyle{empty}
\begin{abstract}
\noindent 
%The Poisson-Boltzmann equation models the electrostatic interaction between bio molecules in a solution consisting of moving ions dissloved in generally high dielectric, such as water. It is widely used in biology and biophysics, for example for calculation of the excitonic couplings between two chromophores attached to both ends of a protein. Despite its wide acceptance, it was not known until 2007 if this equation always has a unique solution and to what functional space it belongs. The main purpose of this article is to examine/ prove the existence and uniqueness of a solution to the general PBE. 
In this paper we have derived explicitly computable bounds on the error in energy norm for the nonlinear Poisson-Boltzmann equation.
Together with the computable bounds, we have also obtained efficient error indicators which can serve as a basis for a reliable adaptive finite
element algorithm. 
 
\end{abstract}

\setcounter{page}{1}
\pagenumbering{roman}

%\newpage
%\tableofcontents
%\newpage
%\listoftables
%\addcontentsline{toc}{subsection}{List of Tables}
%\newpage
%\listoffigures
%\addcontentsline{toc}{subsection}{List of Figures}
%\newpage

%\listofsymbols
%\newpage
\setcounter{page}{1}
\pagenumbering{arabic}
\section{Introduction}

Biomolecular electrostatics plays an important role in the analysis of the molecular structure
of biological macromolecules such as proteins, RNA or
DNA~\cite{Li2015DNA,Tan2009Predicting,Lipfert2014Understanding}.
When modeling various electrostatic effects, a commonly accepted and widely used approach
is based on solving the nonlinear Poisson Boltzmann equation (PBE). Applications include
computations of the electrostatic potential of biomolecules in solution, the encounter rate coefficient,
free energy of association in conjunction with its salt dependence, or pKa values of such molecules.
Biomolecular association, e.g., the association of ligand and proteins, depends in a complex
manner on the shape of the molecules and their electrostatic fields. Therefore, predictions by
mathematical models have to take into account both shape and charge distribution effects,
cf.~\cite{Fogolari_Brigo_Molinari_2002}.

The Poisson-Boltzmann equation introduced by Gouy~\cite{Gouy1910}
and Chapman~\cite{Chapman1913} describes the electrochemical potential of ions
in the diffuse layer caused by a charged solid that comes into contact with an ionic solution,
creating a layer of surface charges and counter-ions in the form of a double layer. The model
accounts for the thermal motion of ions that behave as point charges. It has been generalized
by Debye and Huckel to provide a theory for the electrostatic interaction of ions in electrolyte
solutions~\cite{DebyeHuckel1923}.

Simple-shape molecular models, e.g., electrostatic models for globular proteins as used
in~\cite{Kirkwood1934}, had been replaced in the early 1980s by models based on more
complex geometries. This  development was driven by the progress of finite element (FE),
boundary element~(BE), and finite difference (FD) methods for solving nonlinear partial
differential equations (PDE), see e.g. \cite{RecentProgress}.
Numerous software packages for the simulation of biomolecular electrostatic effects that
are presently available, such as APBS, CHARMM, DelPhi and UHBD, reflect the popularity
and success of the PBE model. 

Major advances in the quality of the numerical solution of the PBE regarding accuracy and
efficiency are due to proper regularization and mesh adaptation techniques, see,
e.g.,~\cite{Gilson1993,Chen2006b,Holst2012}. 
Adaptive FE methods exploit error indicators, which must be reliable and efficient in
that upon multiplication by constants of the same order they provide bounds for the
actual error from above and below. 
Efficient error indicators can be constructed by different methods closely related to different
approaches to the a posteriori error estimation problem. In this context, we mention residual
based methods, goal-oriented methods, methods based on post-processing of numerical
solutions (e.g., averaging or equilibration), and functional type methods.
The latter have been developed in the framework of duality theory for convex variational
problems~\cite{Ekeland_Temam,Repin_2000,Repin_noninear_var_problems_2012}.
They provide estimates that generate guaranteed tight bounds on the distance to the exact
solution valid for the whole class of energy admissible functions (see, e.g.,~\cite{Repin}).
These estimates contain neither mesh dependent constants nor do they rely on any special
conditions or assumptions on the exact solution (e.g., higher regularity) or approximation (e.g.,
Galerkin orthogonality), which means that they are fully
computable. For these
reasons, they are very convenient to use and the error analysis presented in this paper is based
upon this approach. 

\begin{figure}[!htb]
    \centering
    \begin{minipage}{1\textwidth}
        \centering
        \captionsetup{width=0.8\linewidth}
      \includegraphics[width=0.5\linewidth]{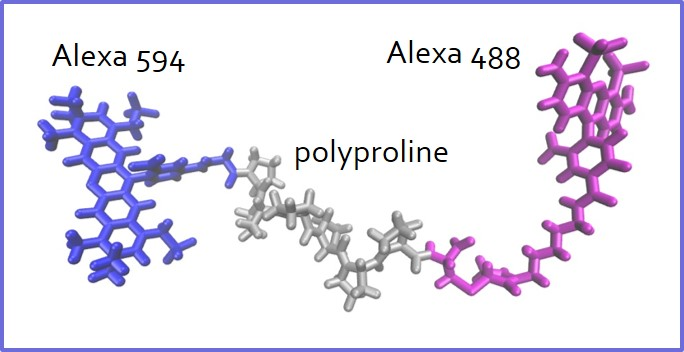}
        \caption{{\small Polyproline 6 labeled with Alexa 488 and Alexa 594~\cite{Sobakinskaya_Busch_Renger_2018}.}}
        \label{polyproline03}
    \end{minipage}%
    \end{figure}

The present paper is a continuation of a recent work by the authors (\cite{Kraus_Nakov_Repin_PBE1_2018})
and is devoted to adaptive modeling of electrostatic interactions of biomolecules. We use two test systems
on which the theoretical findings are demonstrated. The first system consists of two chromophores Alexa 488
and Alexa 594 (Figure \ref{polyproline03}). These chromophores are frequently used for protein labeling in
biophysical experiments. Here we are interested in calculating the electrostatic interaction between them.
The interaction of the dyes, especially the charged ones, such as in our case, influences their conformational states
and orientations, that impacts results of the FRET (F{\"o}rster Resonance Energy Transfer) experiment.
Thus, interpretation and prediction of the experimental results depends on detailed understanding of the
chromophores dynamics. The second test is performed on an insulin protein with a PDB ID 1RWE. This is a
small protein that functions in the hormonal control of metabolism \cite{Wan_et_al_2004}. Because of its 
mportance in the treatment of diabetes mellitus, this protein has attracted attention as a target of protein
engineering. In recent years analogues have gained widespread clinical acceptance \cite{Buse_2001}.
Despite such empirical success, how insulin binds to the insulin receptor is not well understood.
The important contribution to the binding may come from electrostatic interactions between two molecules.
The Poisson-Boltzmann equation can be used to calculate the electrostatic surface of insulin, which would
help to determine the binding sites. Then, the distribution of the electrostatic potential around the protein can be
used in simulations of binding dynamics. Typically, the exact solution of such problems behaves very differently
in different parts of the domain and it is often impossible to a priori locate the zones with complicated behavior
(e.g., high gradients, oscillations, singularities, etc.). Therefore, a crucial component in this approach is mesh
adaptation, which requires robust and efficient error indicators. 

The main contribution of this work is to develop error control methods that allow for a fully reliable mathematical
modeling of the class of problems in question.
The paper is organized as follows. In Section~\ref{Section_Problem_Formulation} a class of nonlinear interface
problems describing the electrostatic potential of biomolecules is presented. The general problem, which is
governed by the nonlinear PBE, is first posed in a classical form, discussing also different regularization
techniques based on two- and three-term splittings, and then in a variational form. 
Section~\ref{Section_A_ Posteriori_Error_Estimates}
focuses on the derivation of error majorants and minorants for the individual components of the solution that
appear in the different splittings leading to reliable and fully computable a posteriori estimates as well as efficient
and robust indicators for the overall error.
Near best approximation results are also proven in Section~\ref{Section_A_ Posteriori_Error_Estimates} for
different regularization techniques. Finally, in Section~\ref{Section_Numerical_Results}, theoretical discussions
are complemented by extensive numerical tests that demonstrate the reliability of the presented methods.

\section{Problem formulation}\label{Section_Problem_Formulation}

\subsection {Classical form of the problem}
In this paper we consider an interface problem describing the electrostatic potential in a system consisting of a (macro)molecule embedded
in a solution, e.g., of a solvent like water and a solute like $\rm {NaCl}$. The computational domain  $\Omega\subset \mathbb R^d$ is assumed
to be bounded with Lipschitz boundary $\partial \Omega$. The domain containing the molecule is denoted by $\Omega_m\subset\mathbb R^d$
and assumed to be strictly inside $\Omega$, i.e. $\overline \Omega_m \subset \Omega$ and also with Lipschitz boundary.
The domain containing the solution with the moving ions of the solute is denoted by $\Omega_s$ and is defined by
$\Omega_s = \Omega \setminus \overline\Omega_m$. The interface of $\Omega_m$ and $\Omega_s$ is denoted by $\Gamma=\overline \Omega_m\cap\overline\Omega_s=\partial\Omega_m$, and the outward unit normal vector on $\partial \Omega_m = \Gamma$ by $n$.
\begin{figure}[h]
\centering
    \includegraphics[width=0.35\linewidth]{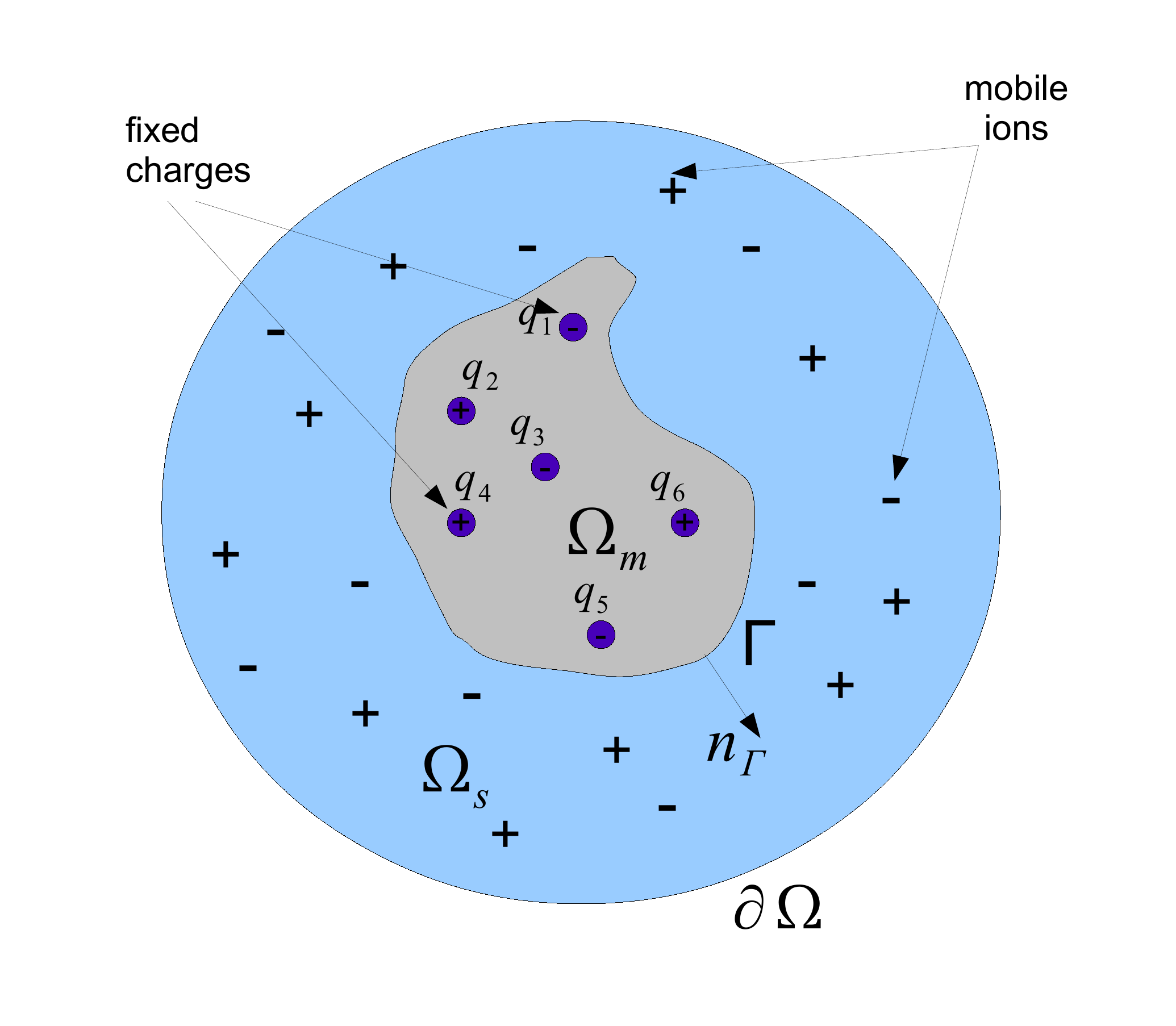} 
%\caption{The simple case of one molecular domain and one solution domain as considered by Holst, %McCammon, Yu, Zhou, Zhu- 2012 and by Holst, Chen, Xu-2007.}
\caption{Computational domain $\Omega$ with molecular domain $\Omega_m$ and solution domain $\Omega_s$.}
  \label{Domain1}
\end{figure}
Assuming that we have only two ion species in the solution with the same concentration $c(x)$ (which are univalent but with opposite charge, i.e $q_j=(-1)^j e_c,\,j=1,2,\,e_c>0$), the interface problem reads as follows
\begin{subequations}\label{PBE_special_form}
\begin{eqnarray}
-\nabla\cdot\left(\epsilon\nabla \phi\right)+8\pi ce_c\,\sinh\left(\frac{e_c \phi}{k_B T}\right)&=&4 \pi\sum\limits_{i=1}^{N}{q_i\delta_{x_i}(x)}\quad \text{in }\Omega_m\cup\Omega_s, \label{PBE_special_form_1}\\
	\left[\phi\right]_\Gamma&=&0,\label{PBE_special_form_2}\\
	\left[\epsilon\nabla \phi\cdot n\right]_\Gamma&=&0,\label{PBE_special_form_3}\\
	\phi&=&\psi\quad\text{on }\partial\Omega,\label{PBE_special_form_4}
	\end{eqnarray}
	\end{subequations}
where $\psi \in C^{0,1}(\partial\Omega)$ and $e_c$ is the electron charge. The function $c(x)\ge 0$ is given by
\begin{align*}
&c(x)=\left\{
\begin{array}{lll}
&0,\, &x\in\Omega_m ,\\
&c^\infty,\,&x \in\Omega_s ,
\end{array}
\right.
\end{align*}
where $c^\infty$ is a positive constant and represents the bulk concentration of the two ion species in solution, $k_B$ is the Boltzmann constant, $T$ is the absolute temperature (constant), $x_i$ is the position of the $i$-th fixed point charge in the molecular domain $\Omega_m$, $\delta_{x_i}$ is the delta distribution centered at $x_i$, and $\phi(x)$ is the unknown electrostatic potential.
The coefficient $\epsilon(x)$ (dielectric constant) is piecewise constant, i.e.
\begin{align}\label{definition_of_eps}
&\epsilon(x)=\left\{
\begin{aligned}
&\epsilon_m, &x&\in \Omega_m,\\
&\epsilon_s, &x&\in\Omega_s.
\end{aligned}
\right.
\end{align}
Introduce the new variable $\tilde \phi=\frac{e_c\phi}{k_BT}$. Then writing $q_i=z_ie_c$, where $z_i\in \mathbb Z$ we obtain the dimensionless form of \eqref{PBE_special_form}:
\begin{subequations}\label{PBE_dimensionless}
\begin{eqnarray}
-\nabla\cdot\left(\epsilon(x)\nabla \tilde \phi\right)+k^2(x) \sinh\left(\tilde \phi\right)&=&\frac{4\pi e_c^2}{k_BT}\sum\limits_{i=1}^{N}{z_i\delta_{x_i}(x)} \quad\text{in } \Omega_m\cup \Omega_s, \label{PBE_dimensionless_1}\\
	\left[\tilde \phi\right]_\Gamma&=&0,\label{PBE_dimensionless_2}\\
	\left[\epsilon\nabla \tilde \phi\cdot n\right]_\Gamma&=&0,\label{PBE_dimensionless_3}\\
	\tilde \phi&=&g \quad\text{on } \partial\Omega , \label{PBE_dimensionless_4}
	\end{eqnarray}
	\end{subequations}
where $g=\frac{e_c\psi}{k_BT}$ and the coefficient $k(x)$ is given by
\begin{align}\label{definition_of_k}
&k^2(x)=\frac{8\pi c(x) e_c^2}{k_BT}=\left\{
\begin{aligned}
&k_m^2=0, &x&\in \Omega_m,\\
&k_s^2=\frac{8\pi c^\infty e_c^2}{k_BT}, &x&\in\Omega_s .
\end{aligned}
\right.
\end{align}
Equation \eqref{PBE_special_form_1} is often referred to as the Poisson-Boltzmann equation (PBE) \cite{Sharp_Honig_1990}, \cite{Oberoi_Allewell_1993}, \cite{Holst2012}. Equations \eqref{PBE_special_form_2} and \eqref{PBE_special_form_3} are the interface conditions. Here $\left[f\right]_\Gamma$ denotes the jump of a function $f:\Omega_m\cup\Omega_s\to\mathbb R$
that is uniformly continuous in $U_\Gamma \cap \Omega_m$ and  $U_\Gamma \cap \Omega_s$, where $U_\Gamma$
is a neighborhood of $\Gamma$, that is,
\begin{align*}
\left[f\right]_\Gamma:=\gamma_\Gamma\left(f{\restriction_{\Omega_m}}\right)-\gamma_\Gamma\left(f{\restriction_{\Omega_s}}\right)
\end{align*}
with $\gamma_\Gamma(f{\restriction_X})$ denoting the unique extension
of $f{\restriction_X}$ by continuity to $\Gamma$ and $f{\restriction_X}$ denoting the restriction
of $f$ to $X\in\{\Omega_m,\Omega_s\}$.

%where $f{\restriction_X}$ denotes the restriction of $f$ to $X\in\{\Omega_m,\Omega_s\}$ and $\gamma_\Gamma(f{\restriction_X})$ the trace of this restriction on $\Gamma\subset \partial X$ with $\partial X$ denoting the boundary of $X$.\\
%
%\begin{remark}

We notice that in fact the physical problem prescribes a vanishing potential at infinite distance from
the boundary of $\Omega_m$, i.e., $\lim_{|x|\to \infty}{\phi(x)}=0$. In practice, one uses a bounded computational domain and imposes the boundary condition \eqref{PBE_special_form_4} instead, where the function $\psi$ can usually be calculated accurately enough by solving a simpler problem, possibly with a known analytical solution.

\subsection{2-term and 3-term splittings of the solution}
A commonly used technique (see, e.g. \cite{Chen2006b}) to solve and analyze problem \eqref{PBE_dimensionless} is to split the solution according to $\tilde \phi=G+u$, where $G$ is the analytically
known solution of the problem
\begin{align}\label{PDE_for_G}
-\nabla\cdot(\epsilon_m \nabla G)=\frac{4\pi e_c^2}{k_BT}\sum\limits_{i=1}^{N}{z_i\delta_{x_i}(x)},\quad\text{on }\mathbb R^d,\, d\in\{2,3\}.
\end{align} 
If $d=2$, then 
\begin{align}\label{expression_for_G_2d}
G=\sum\limits_{i=1}^{N}{G_i}=-\frac{2\pi e_c^2}{\epsilon_m k_BT}\sum\limits_{i=1}^{N}{z_i \ln{|x-x_i|}}
\end{align}
and if $d=3$
\begin{align}\label{expression_for_G}
G=\sum\limits_{i=1}^{N}{G_i}=\frac{e_c^2}{\epsilon_m k_BT}\sum\limits_{i=1}^{N}{\frac{z_i}{|x-x_i|}}.
\end{align}
The regular component $u$ is assumed to be a function in $H^1(\Omega)$.
Since $\epsilon_m$ and $\epsilon_s$ are constants, the problem for $u$ is as follows:
 \begin{subequations}\label{2term_regular_component}
 \begin{eqnarray}
-\nabla\cdot\left(\epsilon\nabla u\right)+k^2 \sinh\left(u+G\right)&=&0\quad\text{in } \Omega_m\cup \Omega_s,\label{2term_regular_component_1}\\
\left[u\right]_\Gamma&=&0, \label{2term_regular_component_2}\\
\left[\epsilon\nabla u\cdot n\right]_\Gamma&=&-\left[\epsilon\nabla G\cdot n\right]_\Gamma, \label{2term_regular_component_3}\\
u&=&g-G\quad \text{on}\quad \partial\Omega. \label{2term_regular_component_4}
\end{eqnarray}
\end{subequations}
For further analysis of \eqref{2term_regular_component}, it is convenient to split $u$ into $u^N+u^L$, where $u^L$ solves
 \begin{subequations}\label{PBE_2term_regular_linear_part}
	 \begin{eqnarray}
-\nabla\cdot\left(\epsilon \nabla u^L\right)&=&0 \quad\text{in } \Omega_m\cup \Omega_s, \label{PBE_2term_regular_linear_part_1}\\
	\left[u^L\right]_\Gamma&=&0,\label{PBE_2term_regular_linear_part_2}\\
	 \left[\epsilon\frac{\partial u^L}{\partial n}\right]_\Gamma&=&-\left[\epsilon \nabla G\cdot n\right], \label{PBE_2term_regular_linear_part_3}\\
	 u^L&=&g-G,\quad \text{on } \partial \Omega , \label{PBE_2term_regular_linear_part_4}
	\end{eqnarray}
\end{subequations}
 and $u^N$ solves the homogeneous nonlinear problem%~\eqref{PBE_special_form_regular_nonlinear_part}
\begin{subequations}\label{PBE_2term_regular_nonlinear_part}
\begin{eqnarray}
-\nabla\cdot\left(\epsilon \nabla u^N\right)+k^2 \sinh(u^N+G+u^L)&=&0 \quad\text{in } \Omega_m\cup \Omega_s, \label{PBE_2term_regular_nonlinear_part_1}\\
	\left[u^N\right]_\Gamma&=&0,\label{PBE_2term_regular_nonlinear_part_2}\\
	 \left[\epsilon\frac{\partial u^N}{\partial n}\right]_\Gamma&=&0, \label{PBE_2term_regular_nonlinear_part_3}\\
	 u^N&=&0,\quad \text{on } \partial \Omega. \label{PBE_2term_regular_nonlinear_part_4}
	  \end{eqnarray}
	 \end{subequations}
Notice that the problem \eqref{PBE_2term_regular_nonlinear_part} uses the exact solution $u^L$ of \eqref{PBE_2term_regular_linear_part}. 

However, the splitting $\tilde \phi=G+u$ causes numerical instability when $G$ and $u$ are with opposite signs but
nearly of the same absolute value in $\Omega_s$. This mostly happens when $\epsilon_m \ll \epsilon_s$
(see, e.g.~\cite{Holst2012}). In order to overcome this difficulty one may use a splitting of $\tilde \phi$ into $3$ components,
two of which add up to zero in $\Omega_s$. Such a splitting is given by
 \begin{equation}\label{3_term_reg}
 \tilde \phi=G+u^H+u, 
 \end{equation}
 such that $\tilde \phi=u$ in $\Omega_s$, i.e. $u^H=-G$ in $\Omega_s$, and has been used in~\cite{Holst2012}.
%\begin{remark}
%Note that \eqref{PBE_dimensionless} does not satisfy the requirements that are necessary to guarantee a classical solution in the usual sense due to the presence of the delta distributions. By introducing the splitting of $\tilde \phi$ into $G$ and a regular remainder, however, the resulting problem for the regular component $u$ in general has a classical solution.
%\end{remark}

The no-jump condition \eqref{PBE_dimensionless_2} on $\Gamma$ can be expressed as:
\begin{align*}
\left[\tilde \phi\right]_\Gamma:=\gamma_\Gamma\left((G+u^H+u){\restriction_{\Omega_m}}\right)-\gamma_\Gamma\left(u{\restriction_{\Omega_s}}\right)=0.
\end{align*}

If we further require $u^H$ to be continuous across $\Gamma$, then $\left[u\right]_\Gamma=0$. Consequently, we arrive at the following system of equations:
\begin{eqnarray}\label{equation_in_omega_s}
u^H&=&-G \quad \text{in } \Omega_s, \nonumber \\
-\nabla\cdot\left(\epsilon_s\nabla u\right)+k^2 \sinh(u)&=&0 \quad \text{in } \Omega_s,
\end{eqnarray}
and
\begin{equation}\label{equation_in_omega_m}
-\nabla\cdot\left(\epsilon_m\nabla\left(u+u^H+G \right)\right)= \frac{4\pi e_c^2}{k_BT}\sum\limits_{i=1}^{N}{z_i\delta_{x_i}(x)} \quad \text{in } \Omega_m.
\end{equation}
Hence, defining $u^H$ in $\Omega_m$ to be the solution of 
\begin{subequations}\label{PBE_special_form_harmonic}
\begin{eqnarray}
	-\Delta u^H&=&0 \quad\text{in}\quad\Omega_m, \label{PBE_special_form_harmonic_a}\\
	u^H&=&-G \quad\text{on} \quad \Gamma=\partial\Omega_m, \label{PBE_special_form_harmonic_b}
\end{eqnarray}
\end{subequations}
taking into account \eqref{PDE_for_G} and recalling that $\epsilon_m$ is a constant, we get 
\begin{align}\label{equation_in_omega_m_regular_part_3_term}
	-\nabla\cdot\left(\epsilon_m \nabla u\right)=0\quad\text{in }\Omega_m.
\end{align}
In view of  \eqref{definition_of_eps} and \eqref{definition_of_k} we can represent \eqref{equation_in_omega_s} and \eqref{equation_in_omega_m_regular_part_3_term} in one common
form, namely,
\begin{align*}
	-\nabla\cdot\left(\epsilon \nabla u\right)+k^2 \sinh(u)=0\quad\text{in }\Omega_m\cup\Omega_s.
\end{align*}
In order to find the interface condition for the flux $\epsilon\nabla u\cdot n$, we note that 
%\begin{align*}
$
\left[\epsilon\frac{\partial \tilde \phi}{\partial n}\right]=0.
$
%\end{align*}
From this condition, we deduce the relation
$$
\gamma_\Gamma\left({\epsilon_m\frac{\partial u}{\partial n}}{\restriction_{\Omega_m}}\right)+\gamma_\Gamma\left(\epsilon_m\frac{\partial\left(u^H+G\right)}{\partial n}{\restriction_{\Omega_m}}\right)=\gamma_\Gamma\left(\epsilon_s\frac{\partial u}{\partial n}{\restriction_{\Omega_s}}\right)
$$
or, equivalently,
$$
\left[\epsilon\frac{\partial u}{\partial n}\right]_\Gamma=-\gamma_\Gamma \left(\epsilon_m\frac{\partial\left(u^H+G\right)}{\partial n}{\restriction_{\Omega_m}}\right)=:g_\Gamma
$$
Thus, we arrive at the following interface problem for the regular component $u$:\\
\begin{subequations}\label{PBE_special_form_regular}
\begin{eqnarray}
-\nabla\cdot\left(\epsilon \nabla u\right)+k^2 \sinh(u)&=&0 \quad\text{in } \Omega_m\cup \Omega_s, \label{PBE_special_form_regular_1}\\[0.5em]
	\left[u\right]_\Gamma&=&0,\label{PBE_special_form_regular_2}\\
	 \left[\epsilon\frac{\partial u}{\partial n}\right]_\Gamma&=&g_\Gamma, \label{PBE_special_form_regular_3}\\
	 u&=&g,\quad \text{on } \partial \Omega. \label{PBE_special_form_regular_4}
	\end{eqnarray}
\end{subequations}
For further analysis of problem~\eqref{PBE_special_form_regular}, we split the regular component $u$ into $u^N+u^L$, where  $u^L$ solves the linear nonhomogeneous interface problem%~\eqref{PBE_special_form_regular_linear_part}
 \begin{subequations}\label{PBE_special_form_regular_linear_part}
	 \begin{eqnarray}
-\nabla\cdot\left(\epsilon \nabla u^L\right)&=&0 \quad\text{in } \Omega_m\cup \Omega_s, \label{PBE_special_form_regular_linear_part_1}\\
	\left[u^L\right]_\Gamma&=&0,\label{PBE_special_form_regular_linear_part_2}\\
	 \left[\epsilon\frac{\partial u^L}{\partial n}\right]_\Gamma&=&g_\Gamma, \label{PBE_special_form_regular_linear_part_3}\\
	 u^L&=&g,\quad \text{on } \partial \Omega \label{PBE_special_form_regular_linear_part_4}
	\end{eqnarray}
\end{subequations}
 and $u^N$ solves the homogeneous nonlinear problem%~\eqref{PBE_special_form_regular_nonlinear_part}
\begin{subequations}\label{PBE_special_form_regular_nonlinear_part}
\begin{eqnarray}
-\nabla\cdot\left(\epsilon \nabla u^N\right)+k^2 \sinh(u^N+u^L)&=&0 \quad\text{in } \Omega_m\cup \Omega_s, \label{PBE_special_form_regular_nonlinear_part_1}\\
	\left[u^N\right]_\Gamma&=&0,\label{PBE_special_form_regular_nonlinear_part_2}\\
	 \left[\epsilon\frac{\partial u^N}{\partial n}\right]_\Gamma&=&0, \label{PBE_special_form_regular_nonlinear_part_3}\\
	 u^N&=&0,\quad \text{on } \partial \Omega . \label{PBE_special_form_regular_nonlinear_part_4}
	  \end{eqnarray}
	 \end{subequations}
Again, \eqref{PBE_special_form_regular_nonlinear_part} includes the exact solution of \eqref{PBE_special_form_regular_linear_part} as a known function.  %In what follows, we will present most of our analysis with respect to the splitting $\tilde \phi=G+u^H+u$ and we note that all the results, with some obvious modifications, apply also to the splitting $\tilde \phi=G+u$.

\subsection{Variational form of the problem}\label{Section_Variational_form_of_the_problem}

It is easy to see that the generalized solution $$u^L\in H_g^1(\Omega)=\{v\in H^1(\Omega): \gamma_{\partial\Omega}(v)=g \text{ on } \partial \Omega \}$$ of \eqref{PBE_special_form_regular_linear_part} solves the variational problem 
\begin{equation}\label{ul_weak_formulation}
a(u^L,v)=\langle g_\Gamma, \gamma_\Gamma(v)\rangle_\Gamma, \forall v\in H_0^1(\Omega),
\end{equation}
where 
\begin{equation}\label{def_of_bilinear_form_a}
a(u,v):=\int_{\Omega}{\epsilon \nabla u\cdot\nabla v dx},
\end{equation}
$\gamma_{\partial\Omega}(v)$ is the trace of $v$ on $\partial\Omega$, $\gamma_\Gamma(v)$ is the trace of $v$ on $\Gamma$, and $\langle .,.\rangle_\Gamma$ denotes the duality pairing in $H^{-1/2}(\Gamma)\times H^{1/2}(\Gamma)$.

If the distribution $g_\Gamma\in H^{-1/2}(\Gamma)$ is  regular so that the action on any function $v\in H^{1/2}$ can be represented by the integral $\int_\Gamma g_\Gamma v ds$,
where $g_\Gamma\in L^2$, then we the right hand side in~\eqref{ul_weak_formulation} can be written in the form $\int_\Gamma g_\Gamma \gamma_\Gamma ds$ for all $v\in H_0^1(\Omega)$. \\
For the $3$-term regularization, if $\nabla u^H$ is uniformly continuous in a neighborhood of the interface $\Gamma$, since $\nabla G$ is smooth in a neighborhood of the interface $\Gamma$, then we have 
\begin{align*}
g_\Gamma=-\gamma_\Gamma\left(\epsilon_m\frac{\partial(u^H+G)}{\partial n}{\restriction_{\Omega_m}}\right)\in L^2(\Gamma)
\end{align*}
and in this case we can use $\int\limits_{\Gamma}{g_\Gamma v ds}$ on the RHS of \eqref{ul_weak_formulation}. We can also write
\begin{align}
g_\Gamma=-\gamma_\Gamma(\epsilon_m \nabla u^H\cdot n{\restriction_{\Omega_m}})-\gamma_\Gamma(\epsilon_m\nabla G\cdot n{\restriction_{\Omega_s}}).
\end{align}
Therefore, we see that if $\nabla u^H$ is only in $H(\div;\Omega_m)$, then the functional $g_\Gamma\in H^{-1/2}(\Gamma)$ is defined as follows
\begin{align}\label{definition_of_functional_gGamma_3term_regularization}
\langle g_\Gamma,v\rangle_\Gamma = -\langle \gamma_{n,\Omega_m}(\epsilon_m\nabla u^H),v\rangle_\Gamma + \langle \gamma_{n,\Omega_s}(\epsilon_m\nabla G),v\rangle_\Gamma,\,\forall v\in H^{1/2}(\Gamma),
\end{align}
where $\gamma_{n,X}$ denotes the normal trace in the space $H(\div;X)$. Now, using the divergence formula, the weak formulation \eqref{ul_weak_formulation} can be rewritten as:
\begin{align}\label{definition_ul_weak_formulation_with_functional_gGamma_3term_regularization}
\begin{aligned}
&\text{Find }u^L\in H_g^1(\Omega) \text{ such that}&\\
&\int\limits_{\Omega}{\epsilon\nabla u^L\cdot\nabla v dx}=-\int\limits_{\Omega_m}{\div{(\epsilon_m\nabla u^H)}v dx}-\int\limits_{\Omega_m}{\epsilon_m\nabla u^H\cdot\nabla v dx}\\
&+\int\limits_{\Omega_s}{\div{(\epsilon_m\nabla G)}v dx}+\int\limits_{\Omega_s}{\epsilon_m\nabla G\cdot\nabla v dx}\\
&=-\int\limits_{\Omega_m}{\epsilon_m\nabla u^H\cdot\nabla v dx}+\int\limits_{\Omega_s}{\epsilon_m\nabla G\cdot\nabla v dx},\,\forall v\in H_0^1(\Omega).
\end{aligned}
\end{align}
For the $2$-term regularization, $g_\Gamma$ is known exactly, and it is given by the relation 
\[
g_\Gamma=-[\epsilon\nabla G\cdot n]\in L^2(\Gamma)
\]
Here is used the fact that $\nabla G$ is smooth in a neighborhood of the interface $\Gamma$. 
Since
\begin{align}\label{definition_of_functional_gGamma_2term_regularization}
\langle g_\Gamma, \gamma_\Gamma(v) \rangle_\Gamma = \int\limits_{\Omega}{(\epsilon_m-\epsilon)\nabla G\cdot \nabla v dx}=\int\limits_{\Omega_s}{(\epsilon_m-\epsilon_s)\nabla G\cdot \nabla v dx},\,\forall v\in H_0^1(\Omega),
\end{align}
the integral relation that defines $u^L$ in the $2$-term regularization comes in the form:
\begin{align}\label{definition_ul_weak_formulation_with_functional_gGamma_2term_regularization}
\begin{aligned}
&\text{Find }u^L\in H_{g-G}^1(\Omega) \text{ such that}\\
&\int\limits_{\Omega}{\epsilon\nabla u^L\cdot\nabla v dx}=\int\limits_{\Omega_s}{(\epsilon_m-\epsilon_s)\nabla G\cdot \nabla v dx},\,\forall v\in H_0^1(\Omega),\,\forall v\in H_0^1(\Omega).
\end{aligned}
\end{align}
The well-posedness of~\eqref{ul_weak_formulation} (or, equivalently, \eqref{definition_ul_weak_formulation_with_functional_gGamma_3term_regularization}) follows from the Lax-Milgram Lemma. Moreover, if $\partial\Omega$ and $\Gamma$ are sufficiently smooth, for example, if $\partial\Omega$ is Lipschitz and $\Gamma$ is $C^1$, then from \cite{Optimal_regularity_for_elliptic_transmission_problems} it follows that $u^L\in W^{1,p}(\Omega)$ for some $p>d=3$. Using the embedding theorems, we conclude that $u^L\in L^\infty(\Omega)$: Denote by $V^*$ the topological dual of a Banach space $V$ and by $p'=\frac{p}{p-1}$ the H{\"o}lder conjugate of $p$. It has been shown in~\cite{Optimal_regularity_for_elliptic_transmission_problems} that for $\partial\Omega$ being Lipschitz and $\Gamma\in C^1$,
$\Gamma$ is not touching $\partial \Omega$, there exists $p>3$ such that $-\nabla\cdot(\epsilon(x)\nabla u): W_0^{1,q}(\Omega)\to W^{-1,q}(\Omega)=\left(W_0^{1,q'}\right)^*$ is a topological isomorphism for all $q\in (p',p)$. In addition, if $\partial \Omega$ is also $C^1$, then $p$ may be taken to be $\infty$. This result is useful, because for $q>d=3$ we know that the functions in $W^{1,q}(\Omega)$ are H{\"o}lder continuous and thus in $L^\infty(\Omega)$. 

We can then apply this result to the homogenized version of problem~\eqref{ul_weak_formulation}: find $u_0^L\in H_0^1(\Omega)$ such that 
$$a(u_0^L,v)=-a(u_g^L,v)+\langle g_\Gamma, \gamma_\Gamma(v)\rangle_\Gamma,\,\forall v\in H_0^1(\Omega),$$ 
where $u^L=u_g^L+u_0^L\in H_g^1(\Omega)$ with $u_g^L\in W^{1,q}(\Omega)\cap L^\infty(\Omega)$, $q$ is a fixed number such that $3<q<p$, and $\gamma_{\partial\Omega}(u_g^L)=g$. Note that in order to apply Theorem 1.1 in \cite{Optimal_regularity_for_elliptic_transmission_problems} the functionals $\langle g_\Gamma,\gamma_\Gamma(.)\rangle_\Gamma$ and $-a(u_g^L,.)$ need to be well defined, bounded and linear over $W_0^{1,q'}$, where $q'=\frac{q}{q-1}<\frac{3}{2}$. 

To show the boundedness of these functionals we assume additionally that $\Gamma\in C^{1,1}$. In this case, by applying Theorem 2.4.2.5 from \cite{Grisvard_Elliptic_Problems_in_Nonsmooth_Domains} we get that $u^H\in W^{2,2}(\Omega_m)$ and thus by the Sobolev embedding theorem for $d=3$, $\nabla u^H\in L^6(\Omega_m)$ and for $d=2$, $\nabla u^H\in L^r(\Omega_m),\,\forall r: 1\leq r<\infty$. Then, for the $3$-term splitting, using~\eqref{definition_ul_weak_formulation_with_functional_gGamma_3term_regularization} and applying H\"{o}lder inequality, we obtain
\begin{align*}
&\left|\langle g_\Gamma,\gamma_\Gamma(v)\rangle_\Gamma\right|\leq \epsilon_m\left(\|\nabla u^H\|_{L^q(\Omega_m)}+\|\nabla G\|_{L^q(\Omega_s)}\right)\|v\|_{W^{1,q'}(\Omega)},\, \forall v\in W_0^{1,q'},\\
&\left|-a(u_g^L,v)\right|\leq \epsilon_\max \|u_g^L\|_{W^{1,q}(\Omega)}\|v\|_{W^{1,q'}(\Omega)},\, \forall v\in W_0^{1,q'}(\Omega).
\end{align*}
For the $2$-term splitting we will have $u^L=u_{g-G}^L+u_0^L$ where $u_{g-G}^L\in W^{1,q}(\Omega)\cap L^\infty(\Omega)$ with $\gamma_{\partial\Omega}(u_{g-G}^L)=g-G$ and $u_0^L\in H_0^1(\Omega)$. Thus, using \eqref{definition_ul_weak_formulation_with_functional_gGamma_2term_regularization} we obtain, as before,
\begin{align*}
&\left|\langle g_\Gamma,\gamma_\Gamma(v)\rangle_\Gamma\right|\leq (\epsilon_m-\epsilon_s)\|\nabla G\|_{L^q(\Omega_s)}\|v\|_{W^{1,q'}(\Omega)},\, \forall v\in W_0^{1,q'},\\
&\left|-a(u_{g-G}^L,v)\right|\leq \epsilon_\max \|u_{g-G}^L\|_{W^{1,q}(\Omega)}\|v\|_{W^{1,q'}(\Omega)},\, \forall v\in W_0^{1,q'}(\Omega).
\end{align*}
Another way to see that $u^L\in L^\infty(\Omega)$ without assuming that $\Gamma$ is $C^1$ and only assuming that it is Lipschitz is to apply Theorem B.2 from Kinderlehrer and Stampacchia \cite{Kinderlehrer_Stampacchia}. For this, we need only to ensure that for some $s>d$ it holds $\nabla G\in \left[L^s(\Omega_s)\right]^d$, $\nabla u^H\in \left[L^s(\Omega_m)\right]^d$ in the case of the
3-term splitting and that $\nabla G\in \left[L^s(\Omega_s)\right]^d$ in the case of the 2-term splitting (just apply the result of Theorem B.2 to the homogenized versions of \eqref{definition_ul_weak_formulation_with_functional_gGamma_3term_regularization} and \eqref{definition_ul_weak_formulation_with_functional_gGamma_2term_regularization}). Indeed, $\nabla G\in \left[L^s(\Omega_s)\right]^d$ since $G$ is smooth in $\Omega_s$ and $u^H\in W^{1,p}(\Omega)$ for some $p>d$ according to \cite{Optimal_regularity_for_elliptic_transmission_problems}.

\begin{definition}
If $u^N\in H_0^1(\Omega)$ is such that $b(x,u^N+u^L)v\in L^1(\Omega)\,\forall v\in H_0^1(\Omega)\cap L^\infty(\Omega)$ and 
\begin{align}\label{un_weak_formulation}
a(u^N,v)+\int_{\Omega}{b(x,u^N+u^L)v dx}=0,\,\forall v\in H_0^1(\Omega)\cap L^\infty(\Omega),
\end{align}
where $a(.,.)$ is defined in \eqref{def_of_bilinear_form_a} and $b(x,s):=k^2(x)\sinh(s)$ then $u^N$ is called a weak solution of \eqref{PBE_special_form_regular_nonlinear_part}.
Similarly, we call $u^N$ a weak solution of \eqref{PBE_2term_regular_nonlinear_part} if $u^N\in H_0^1(\Omega)$ and $u^N$ is such that $b(x,u^N+G+u^L)v\in L^1(\Omega)\,\forall v\in H_0^1(\Omega)\cap L^\infty(\Omega)$ and 
\begin{align}\label{2_term_regularization_un_weak_formulation}
a(u^N,v)+\int_{\Omega}{b(x,u^N+G+u^L)v dx}=0,\,\forall v\in H_0^1(\Omega)\cap L^\infty(\Omega).
\end{align}
\end{definition}
%We note that the weak formulation for problem \eqref{PBE_2term_regular_nonlinear_part} is similar with the only difference that instead of $u^L$ in \eqref{un_weak_formulation} we have $G+u^L$.
According to \cite{Kraus_Nakov_Repin_PBE1_2018}, we have the following proposition.
\begin{proposition}\label{Proposition_1}
Problem \eqref{un_weak_formulation} has a unique weak solution $u^N\in H_0^1(\Omega)$, which belongs to $L^\infty(\Omega)$ and satisfies $\|u^N\|_{L^\infty(\Omega)}\leq \|u^L\|_{L^\infty(\Omega_s)}$. Similarly, problem \eqref{2_term_regularization_un_weak_formulation} has an unique weak solution $u^N$ which belongs to $L^\infty(\Omega)$ and satisfies $\|u^N\|_{L^\infty(\Omega)}\leq \|G+u^L\|_{L^\infty(\Omega_s)}$. As a consequence, the test functions in \eqref{un_weak_formulation} and \eqref{2_term_regularization_un_weak_formulation} can be taken in $H_0^1(\Omega)$.
\end{proposition}
In certain situations, when using the 3-term splitting, it is better not to split the regular component $u$ additionally into $u^L+u^N$. Such a situation may arise if the norm of $u^L+u^N$ is much smaller than the norm of at least one of $u^L$ and $u^N$ since then small relative errors in $u^L$ and $u^N$ become substantial relative errors in the sum $u^L+u^N$.
Therefore, we also consider solving equation~\eqref{PBE_special_form_regular}.
\begin{definition}
$u$ is called a weak solution of \eqref{PBE_special_form_regular} if $u\in H_g^1(\Omega)$ and $u$ is such that $b(x,u)v\in L^1(\Omega)\,\forall v\in H_0^1(\Omega)\cap L^\infty(\Omega)$ and 
\begin{align}\label{u_weak_formulation}
\begin{aligned}
&a(u,v)+\int_{\Omega}{b(x,u)v dx}=\langle g_\Gamma, \gamma_\Gamma(v)\rangle_\Gamma\\
&=-\int\limits_{\Omega_m}{\epsilon_m\nabla u^H\cdot\nabla v dx}+\int\limits_{\Omega_s}{\epsilon_m\nabla G\cdot\nabla v dx}, \text{ for all } v\in H_0^1(\Omega)\cap L^\infty(\Omega).
\end{aligned}
\end{align}
\end{definition}
To see that \eqref{u_weak_formulation} has a solution, we can define the energy functional $J$ over $H_g^1(\Omega)$, like in \cite{Kraus_Nakov_Repin_PBE1_2018}
\begin{align}\label{definition_of_J_Regular_Part_3_term}
J(v):=\left\{
\begin{aligned}
&\int\limits_{\Omega}{\left[\frac{\epsilon(x)}{2}\abs{\nabla v}^2+k^2\cosh(v)+\epsilon_m\nabla u^H\cdot\nabla v \mathbbm{1}_{\Omega_m}- \epsilon_m\nabla G\cdot\nabla v\mathbbm{1}_{\Omega_s} \right]dx},\\
&\text{ if }  k^2\cosh(v)\in L^1(\Omega),\\\\
&+\infty, \text{ if } k^2\cosh(v) \notin L^1(\Omega),
\end{aligned}
\right.
\end{align}
 and show the existence of a unique minimizer $u\in H_g^1(\Omega)$, where $\mathbbm{1}_{X}$ is the indicator function of the set $X$.
 Then, for the minimizer $u$, using Lebesgue Dominated Convergence Theorem, we can prove that it is indeed a solution to \eqref{u_weak_formulation}. The uniqueness of the solution $u$ 
 of~\eqref{u_weak_formulation} is proven in a similar way to the approach in \cite{Kraus_Nakov_Repin_PBE1_2018}. 
However, an easier approach is to take advantage of the fact that we have already shown existence and uniqueness of a solution to problems \eqref{definition_ul_weak_formulation_with_functional_gGamma_3term_regularization} and \eqref{un_weak_formulation}. Moreover, since $u=u^L+u^N$, where $u^L\in L^\infty(\Omega)$ is a solution to \eqref{definition_ul_weak_formulation_with_functional_gGamma_3term_regularization} and $u^N\in L^\infty(\Omega)$ is a solution to \eqref{un_weak_formulation}, it follows that $u\in H_g^1(\Omega) \cap L^\infty(\Omega)$.

Hence we have the following proposition.
\begin{proposition}\label{Proposition_Existence_uniqueness_3_term_uRegular}
Problem \eqref{u_weak_formulation} has a unique weak solution $u\in H_0^1(\Omega)$ which belongs to $L^\infty(\Omega)$. As a consequence, the test functions in \eqref{u_weak_formulation} can be taken to be only in $H_0^1(\Omega)$.
\end{proposition}

\section{A posteriori error estimates}\label{Section_A_ Posteriori_Error_Estimates}
%{\color{red}\textbf{This part is more like a note:}
%From now on we will assume for simplicity that the interface $\Gamma$ is a triangulated surface which is also Lipschitz (and is conforming with the meshes of the underlying finite element spaces that will be further involved - this is not needed since everywhere in our a posteriori error estimates we do not assume a particular Galerkin FEM approximation - just conforming approximation !). In this case, $\Gamma$ is not $C^{1,1}$ anymore and thus $u^L$ is not guaranteed to be in $L^\infty(\Omega)$. However, this is not a problem for the analysis that follows. Indeed, the solution $\tilde u^L$ is again unique by the same reasoning as for problem \eqref{un_weak_formulation} or equivalently \eqref{variational_problem_for_un} and $\tilde u^L_{h_L}\in L^\infty(\Omega)$ since $\tilde u^L_{h_L}$ is a finite element solution. Therefore, $\tilde u^N$ will also be in $L^\infty(\Omega)$ by the same reasoning as for problem \eqref{un_weak_formulation}. Thus, for the exact flux $\tilde p_N^*=\epsilon\nabla \tilde u^N$ we will still have $\div \tilde p_N^*=\div(\epsilon\nabla \tilde u^N)\in L^\infty(\Omega)$ and there will be no problem with the functional a posteriori error estimate for the nonlinear component $\tilde u^N$ (we can also use the uniformly convex functional in a ball in $L^\infty$ for $F^*$).}
\subsection{Harmonic component $u^H$}\label{Section_estimate_for_harmonic_part}
To get fully reliable error bounds for approximate solutions of problem~\eqref{PBE_dimensionless}, we need first to derive a posteriori estimates for the quantities 
$$\|\nabla (u^H-\tilde u^H)\|_{L^2(\Omega_m)} \quad \mbox{and} \quad \|\nabla u^H-T(\nabla\tilde u^H)\|_{L^2(\Omega_m)},$$
where $\tilde u^H$ is a conforming approximation of $u^H$ and $T(\nabla \tilde u^H)\in H(\div;\Omega_m)$ with $T$ being a regularization operator that maps the
numerical flux $\nabla \tilde u^H$ into $H(\div;\Omega_m)$. 

For the first quantity, we have (see \cite{Repin_RADON_series})
\begin{align}\label{A_posteriori_estimate_primal_energy_norm_harmonic_component}
\|\nabla (u^H-\tilde u^H)\|_{L^2(\Omega_m)}\leq C_{F\Omega_m}\|\div \left(T(\nabla\tilde u^H)\right)\|_{L^2(\Omega_m)}+\|\nabla \tilde u^H-T(\nabla \tilde u^H)\|_{L^2(\Omega_m)}.
\end{align}
For the second quantity, we proceed as follows:
\begin{eqnarray}\label{extended_Prager_Synge_equality}
\|\nabla \tilde u^H-T(\nabla \tilde u^H)\|_{L^2(\Omega_m)}^2&=&\|\nabla(u^H-\tilde u^H)\|_{L^2(\Omega_m)}^2 \nonumber \\ &+&\|\nabla u^H-T(\nabla\tilde u^H)\|_{L^2(\Omega_m)}^2\\
&-&2\int\limits_{\Omega_m}{\left(T(\nabla\tilde u^H)-\nabla u^H\right)\cdot\nabla (\tilde u^H-u^H)dx}.\nonumber
\end{eqnarray}
Thus, using the Cauchy-Schwartz inequality, we obtain
\begin{align*}
&\|\nabla(u^H-\tilde u^H)\|_{L^2(\Omega_m)}^2+\|\nabla u^H-T(\nabla\tilde u^H)\|_{L^2(\Omega_m)}^2\\
&\leq \|\nabla \tilde u^H-T(\nabla \tilde u^H)\|_{L^2(\Omega_m)}^2+2\|\div\left(T(\nabla\tilde u^H)\right)\|_{L^2(\Omega_m)}C_{F\Omega_m}\|\nabla(\tilde u^H-u^H)\|_{L^2(\Omega_m)}\\
&\leq \|\nabla \tilde u^H-T(\nabla \tilde u^H)\|_{L^2(\Omega_m)}^2+\|\nabla (u^H-\tilde u^H)\|_{L^2(\Omega_m)}^2+C_{F\Omega_m}^2\|\div \left(T(\nabla \tilde u^H)\right)\|_{L^2(\Omega_m)}^2
\end{align*}
Finally,
\begin{align}\label{A_posteriori_estimate_dual_energy_norm_harmonic_component}
%\begin{aligned}
&\mathrel{\phantom{=}}\|\nabla u^H-T(\nabla\tilde u^H)\|_{L^2(\Omega_m)} \nonumber\\
&\leq \left(\|\nabla \tilde u^H-T(\nabla \tilde u^H)\|_{L^2(\Omega_m)}^2+C_{F\Omega_m}^2\|\div \left(T(\nabla \tilde u^H)\right)\|_{L^2(\Omega_m)}^2\right)^{1/2}\\ &=:M_{\oplus,H}\left(\tilde u^H,T(\nabla \tilde u^H)\right)\nonumber
%\end{aligned}
\end{align}
If $T(\nabla\tilde u^H)$ is additionally equilibrated (for example using the patchwise equilibration technique in \cite{Braess_Schoberl_2006}) then both estimates \eqref{A_posteriori_estimate_primal_energy_norm_harmonic_component} and \eqref{A_posteriori_estimate_dual_energy_norm_harmonic_component} follow from \eqref{extended_Prager_Synge_equality} (Prager-Synge estimate).

\subsection{Linear nonhomogeneous problem}\label{Section_estimate_for_linear_part}

In this section, we show how to obtain a guaranteed bound on the energy norm of the error $\nabla(\tilde u^L-v)$. Here $v$ is some conforming approximation of $\tilde u^L$, the weak solution of the interface problem \eqref{definition_ul_weak_formulation_with_functional_gGamma_3term_regularization} with $\nabla u^H$ replaced by $T(\nabla \tilde u^H)$. The function $\tilde u^H$ is some conforming approximation of $u^H$ and $T$ is some operator that maps the numerical flux $\epsilon \nabla \tilde u^H$ into $H(\div;\Omega_m)$. 
If $\nabla \tilde u^H$ is already in $H(\div;\Omega_m)$, then we can take $T$ to be the identity. The error estimate that we derive here is similar to the one derived in Chapter 4 from \cite{Repin_RADON_series} with the exception that here we avoid involving the trace constant in $\Omega_m$ by exactly prescribing the jump condition on the interface $\Gamma$. 

The function $\tilde u^L$ satisfies the weak formulation:
\begin{align}\label{tildeul_weak_formulation}
&\text{Find }\tilde u^L\in H_g^1(\Omega) \text{ such that}\nonumber\\
&\mathrel{\phantom{=}}\underline{\int\limits_{\Omega}{\epsilon\nabla \tilde u^L\cdot\nabla \phi dx}=-\int\limits_{\Omega_m}{\epsilon_m T(\nabla\tilde u^H)\cdot\nabla \phi dx}+\int\limits_{\Omega_s}{\epsilon_m\nabla G\cdot\nabla \phi dx}}\nonumber\\
&=-\int\limits_{\Omega_m}{\div(\epsilon_m T(\nabla \tilde u^H))\phi dx}-\int\limits_{\Omega_m}{\epsilon_m T(\nabla\tilde u^H)\cdot \nabla \phi dx}\nonumber\\
&+\int\limits_{\Omega_s}{\div(\epsilon_m\nabla G)\phi dx}+\int\limits_{\Omega_s}{\epsilon_m\nabla G\cdot\nabla \phi dx}+\int\limits_{\Omega_m}{\div(\epsilon_m T(\nabla \tilde u^H))\phi dx}\nonumber\\
&=-\langle \gamma_{n,\Omega_m}(\epsilon_m T(\nabla \tilde u^H)),\gamma_\Gamma(\phi)\rangle_\Gamma + \langle \gamma_{n,\Omega_s}(\epsilon_m\nabla G),\gamma_\Gamma(\phi)\rangle_\Gamma\nonumber\\
&+\int\limits_{\Omega_m}{\div(\epsilon_m T(\nabla \tilde u^H))\phi dx},\,\forall \phi\in H_0^1(\Omega).
\end{align}
Now, let $\tilde y_L^*=y^*_{\tilde g_\Gamma}+\tilde y_0^*$, where 
\begin{align}\label{definition_of_yStar_tildeuL_3term_regularization}
\begin{aligned}
&y^*_{\tilde g_\Gamma}=-\epsilon_m T(\nabla \tilde u^H)\mathbbm{1}_{ \Omega_m }+\epsilon_m\nabla G  \mathbbm{1}_{\Omega_s},\\
&\tilde y_0^*\in H(\div;\Omega) ,
\end{aligned}
\end{align}
and $\mathbbm{1}_{X}$ is the indicator function of the set $X$. From \eqref{definition_of_yStar_tildeuL_3term_regularization} it follows that ${\tilde  {y}_L^*}_{\restriction_{\Omega_m}}\in H(\div;\Omega_m)$ and ${\tilde  {y}_L^*}_{\restriction_{\Omega_s}}\in H(\div;\Omega_s)$. From \eqref{tildeul_weak_formulation}, testing with functions $\phi\in H_0^1(\Omega_m)$ and $\phi\in H_0^1(\Omega_s)$,
we see that 
\begin{align}\label{properties_of_tildepL_3term_regularization}
\begin{aligned}
(\epsilon\nabla\tilde u^L)_{\restriction_{\Omega_m}}&\in H(\div;\Omega_m),\,\div\left((\epsilon\nabla\tilde u^L)_{\restriction_{\Omega_m}}\right)=-\div\left(\epsilon_m T(\nabla \tilde u^H)\right),\\
(\epsilon\nabla\tilde u^L)_{\restriction_{\Omega_s}}&\in H(\div;\Omega_s),\,\div\left((\epsilon\nabla\tilde u^L)_{\restriction_{\Omega_m}}\right)=0.
\end{aligned}
\end{align}
Let $v\in H_g^1(\Omega)$ be a conforming approximation of $\tilde u^L$. We proceed with the derivation of an a posteriori estimate for $\vertiii{\nabla(\tilde u^L-v)}$, where 
\[
\vertiii{q}^2=\int\limits_{\Omega}{\epsilon|q|^2dx},\quad \vertiii{q}_*^2=\int\limits_{\Omega}{\frac{1}{\epsilon}|q|^2 dx} , \quad \forall q\in [L^2(\Omega)]^3 .
\]
Furthermore,
\begin{align*}
&(\epsilon\nabla(\tilde u^L-v),\nabla\phi)=-\langle\gamma_{n,\Omega_m}\left(\epsilon_m T(\nabla\tilde u^H) \right),\gamma_\Gamma(\phi) \rangle_\Gamma+\langle\gamma_{n,\Omega_s}\left(\epsilon_m \nabla G \right),\gamma_\Gamma(\phi) \rangle_\Gamma\\
&+\int\limits_{\Omega_m}{\div\left(\epsilon_m T(\nabla\tilde u^H)\right)\phi dx}+\int\limits_{\Omega_m}{\div(y^*_{\tilde g_\Gamma}+\tilde  y_0^*)\phi dx}+\int\limits_{\Omega_s}{\div(y^*_{\tilde g_\Gamma}+\tilde y_0^*)\phi dx}\\
&-\langle\gamma_{n,\Omega_m}\left(y^*_{\tilde g_\Gamma}+\tilde y_0^*\right),\gamma_\Gamma(\phi) \rangle_\Gamma-\langle\gamma_{n,\Omega_s}\left(y^*_{\tilde g_\Gamma}+\tilde  y_0^* \right),\gamma_\Gamma(\phi) \rangle_\Gamma -(\epsilon\nabla v-\tilde y_L^*,\nabla \phi)\\
&=\int\limits_{\Omega_m}{\div\left(\epsilon_m T(\nabla\tilde u^H) \right)\phi dx}+\int\limits_{\Omega_m}{\div(y^*_{\tilde g_\Gamma}+\tilde y_0^*)\phi dx}+\int\limits_{\Omega_s}{\div(y^*_{\tilde g_\Gamma}+\tilde y_0^*)\phi dx}\\
&-(\epsilon\nabla v-\tilde y_L^*,\nabla \phi)=\int\limits_{\Omega_m}{\div \tilde y_0^*\phi dx}+\int\limits_{\Omega_s}{\div \tilde y_0^*\phi dx}-(\epsilon\nabla v-\tilde y_L^*,\nabla \phi)\\
&=(\div \tilde y_0^*,\phi)-(\epsilon\nabla v-\tilde y_L^*,\nabla \phi),
\end{align*}
where we have used \eqref{definition_of_yStar_tildeuL_3term_regularization}. Now, taking $\phi=\tilde u^L-v\in H_0^1(\Omega)$ we have that
\begin{align*}
\vertiii{\nabla(\tilde u^L-v)}^2\leq \frac{C_{F\Omega}}{\sqrt{\epsilon_\min}}\|\div \tilde y_0^*\|_{L^2(\Omega)}\vertiii{\nabla(\tilde u^L-v)}+\vertiii{\epsilon\nabla v-\tilde y_L^*}_*\vertiii{\nabla(\tilde u^L-v)}
\end{align*}
and thus by dividing by $\vertiii{\nabla(\tilde u^L-v)}$ we obtain
\begin{align}\label{A_posteriori_error_estimate_for_tildeuL_3term_regularization}
\vertiii{\nabla(\tilde u^L-v)}\leq \frac{C_{F\Omega}}{\sqrt{\epsilon_\min}}\|\div \tilde y_0^*\|_{L^2(\Omega)}+\vertiii{\epsilon\nabla v-\tilde y^*_{\tilde g_\Gamma}-\tilde y_0^*}_*=:M_{\oplus,L}(v,\tilde y_0^*)
\end{align}
Now, we show that the estimate \eqref{A_posteriori_error_estimate_for_tildeuL_3term_regularization} is sharp. For this, take $\tilde y_0^*=\tilde p_L^*-y^*_{\tilde g_\Gamma}$, where $\tilde p_L^*=\epsilon\nabla\tilde u^L$. We claim that $\tilde y_0^* \in H(\div;\Omega)$. Indeed, let $\phi\in H_0^1(\Omega)$. Then using \eqref{definition_of_yStar_tildeuL_3term_regularization}, \eqref{properties_of_tildepL_3term_regularization}, and the fact that $G$ is harmonic in $\Omega_s$, we obtain
\begin{align*}
&\int\limits_{\Omega}{\tilde y_0^*\cdot\nabla\phi dx}=\int\limits_{\Omega}{\tilde p_L^*\cdot\nabla\phi dx}-\int\limits_{\Omega}{y^*_{\tilde g_\Gamma}\cdot\nabla\phi dx}\\
&=-\langle \gamma_{n,\Omega_m}\left(\epsilon_m T(\nabla\tilde u^H)\right),\gamma_\Gamma(\phi) \rangle_\Gamma+\langle \gamma_{n,\Omega_s}\left(\epsilon_m \nabla G\right),\gamma_\Gamma(\phi) \rangle_\Gamma\\
&+\int\limits_{\Omega_m}{\div\left(\epsilon_m T(\nabla\tilde u^H) \right)\phi dx}-\int\limits_{\Omega_m}{\div\left(\epsilon_m T(\nabla\tilde u^H) \right)\phi dx}+\langle \gamma_{n,\Omega_m}\left(\epsilon_m T(\nabla\tilde u^H)\right),\gamma_\Gamma(\phi) \rangle_\Gamma\\
&+\int\limits_{\Omega_s}{\div(\epsilon_m\nabla G)\phi dx}-\langle \gamma_{n,\Omega_s}\left(\epsilon_m \nabla G\right),\gamma_\Gamma(\phi) \rangle_\Gamma=0.
\end{align*}
Thus $\tilde y_0^*\in H(\div;\Omega)$ and $\div \tilde y_0^*=0$. Now, substituting $\tilde y_0^*$ into the estimate \eqref{A_posteriori_error_estimate_for_tildeuL_3term_regularization} and using again \eqref{definition_of_yStar_tildeuL_3term_regularization} and \eqref{properties_of_tildepL_3term_regularization}, we obtain that the RHS of \eqref{A_posteriori_error_estimate_for_tildeuL_3term_regularization} is equal to $\vertiii{\nabla(\tilde u^L-v)}$. In practice, to find a sharp bound on the error we can do a minimization of the majorant $M_{\oplus,L}(v,\tilde y_0^*)$ in $\tilde y_0^*$ over a finite dimensional subspace of $H(\div;\Omega)$. However, it is more convenient to minimize the squared majorant $M_{\oplus,L}^2(v,\tilde y_0^*;\alpha)$ simultaneously over $\alpha\in \mathbb R_{>0}:=\{x\in \mathbb R: x>0\}$ and $\tilde y_0^*$ in a finite dimensional subspace of $H(\div;\Omega)$, cf.~\cite{Repin_RADON_series}:
\begin{align}\label{A_posteriori_error_estimate_for_tildeuL_3term_regularization_squared_Majorant}
&\vertiii{\nabla(\tilde u^L-v)}^2\leq M_{\oplus,L}^2(v,\tilde y_0^*)\nonumber\\
&\leq (1+\alpha)\frac{C_{F\Omega}}{\epsilon_\min}\|\div \tilde y_0^*\|_{L^2(\Omega)}^2+\left(1+\frac{1}{\alpha}\right)\vertiii{\epsilon\nabla v-y^*_{\tilde g_\Gamma}-\tilde y_0^*}_*^2:=M_{\oplus,L}^2(v,\tilde y_0^*;\alpha)
\end{align}
Another approach to obtain a sharp bound on the error is to apply an appropriate flux reconstruction, similar to the one we use in Section~\ref{subsec:finding_regular_comp}.

For the $2$-term regularization, we can obtain in a similar way the estimate
\begin{align}\label{A_posteriori_error_estimate_for_tildeuL_2term_regularization}
\vertiii{\nabla(u^L-v)}\leq \frac{C_{F\Omega}}{\sqrt{\epsilon_\min}}\|\div y_0^*\|_{L^2(\Omega)}+\vertiii{\epsilon\nabla v-(\epsilon_m-\epsilon)\nabla G-y_0^*}_* ,
\end{align}
where we write $u^L$ instead of $\tilde u^L$ since there is no approximation error in the interface condition \eqref{PBE_2term_regular_linear_part_3} or equivalently in the RHS of the weak formulation \eqref{definition_ul_weak_formulation_with_functional_gGamma_2term_regularization} and $y_0^*\in H(\div;\Omega)$ is arbitrary. 

\subsection{Nonlinear homogeneous problem}\label{Section_estimate_for_nonlinear_part}
Now, we turn our attention to the problem \eqref{un_weak_formulation} which falls in the class of problems that we have considered in \cite{Kraus_Nakov_Repin_PBE1_2018}. Since in practice, we only
have an approximation $\tilde u^L_{h_L}$ to $u^L$, we consider problem \eqref{un_weak_formulation} with $\tilde u^L_{h_L}$ instead of $u^L$ and we assume that $\tilde u^L_{h_L}\in L^\infty(\Omega)$ which is the case if $\tilde u^L_{h_L}$ is for example a finite element approximation. We denote the exact solution of problem \eqref{un_weak_formulation} by $\tilde u^N$. Applying Proposition \ref{Proposition_1} to problem \eqref{un_weak_formulation} with $\tilde u^L_{h_L}\in L^\infty(\Omega)$ we see that $\|\tilde u^N\|_{L^\infty(\Omega)}\leq \|\tilde u^L_{h_L}\|_{L^\infty(\Omega_s)}$.

From \cite{Kraus_Nakov_Repin_PBE1_2018} we have the following error equality
\begin{align}\label{Explicit_form_of_Functional_error_equality}
\vertiii{\nabla(\tilde u^N-v)}^2&+\vertiii{\tilde p_N^*-\tilde y_N^*}_*^2 \nonumber \\
&+2 D_F(v,-\Lambda^*\tilde p_N^*)+ 2 D_F(\tilde u^N,-\Lambda^*\tilde y_N^*)=2M_{\oplus,N}^2(v,\tilde y_N^*)
\end{align}
where $v$ is an arbitrary conforming approximation of $\tilde u^N$, $\tilde p_N^*=\epsilon\nabla \tilde u^N$, $\tilde y_N^*\in H(\div;\Omega)$ with $\div \tilde y_N^*=0$ in $\Omega_m$ is an approximation of the exact flux $\tilde p_N^*$, $\Lambda^*\equiv -\div$,  and $\vertiii{.}$, $\vertiii{.}_*$ denote the primal and the dual energy norm, respectively, given by $\vertiii{q}^2:=\int\limits_{\Omega}{\epsilon|q|^2dx}$, $\vertiii{q}_*^2=\int\limits_{\Omega}{\frac{1}{\epsilon}|q|^2dx}$ for all $q\in \left[L^2(\Omega)\right]^d$. The majorant $2M_{\oplus,N}^2(v,\tilde y_N^*)$ is given by
\begin{align}\label{Explicit_form_of_Nonlinear_Majorant}
2M_{\oplus,N}^2(v,\tilde y_N^*)=\int\limits_{\Omega}{\eta^2(x)dx}=\vertiii{\epsilon\nabla v-\tilde y_N^*}_*^2+2D_F(v,-\Lambda^*\tilde y_N^*),
\end{align}
where 
\begin{align}
\label{Explicit_form_DF_v_minus_lambda_yStar}
&\mathrel{\phantom{=}}D_F(v,-\Lambda^*\tilde y_N^*)=-\int\limits_{\Omega_s}{\div \tilde y_N^*vdx}\\ \nonumber
&+\int\limits_{\Omega_s}{k^2\left(\cosh(v+\tilde u_{h_L}^L)+\frac{\div \tilde y_N^*}{k^2}\left(\arsinh\left(\frac{\div \tilde y_N^*}{k^2}\right)-\tilde u_{h_L}^L\right)-\sqrt{\left(\frac{\div \tilde y_N^*}{k^2}\right)^2+1}\right)dx}
\end{align}
The quantities $D_F(v,-\Lambda^*\tilde p_N^*)$ and $D_F(\tilde u^N,-\Lambda^*\tilde y_N^*)$ are non-negative and measure the error in $v-u$ and in $\div \tilde y_N^*-\div\tilde p_N^*$, respectively, as it is shown in \cite{Kraus_Nakov_Repin_PBE1_2018}. Since $D_F(v,-\Lambda^*\tilde p_N^*)\ge 0$ and $D_F(\tilde u^N,-\Lambda^*\tilde y_N^*)\ge 0$ we have an upper bound for the error in the combined energy norm. It is also easy to obtain a lower bound for the  same error (see \cite{Kraus_Nakov_Repin_PBE1_2018}). These two bounds can be written as
\begin{align}\label{strengthened_error_inequality_Th7_2_1}
&\frac{1}{2}\vertiii{\epsilon \nabla v-\tilde y_N^*}_*^2\leq\vertiii{\nabla(\tilde u^N-v)}^2+\vertiii{\tilde p_N^*-\tilde y_N^*}_*^2\leq 2M_{\oplus,N}^2(v,\tilde y_N^*)
\end{align}
In \cite{Kraus_Nakov_Repin_PBE1_2018}, the following practical estimation for the error in the combined energy norm is suggested
\begin{align}\label{practical_estimation_of_combined_energy_norm_error}
\vertiii{\nabla (v-\tilde u^N)}^2+\vertiii{\tilde y_N^*-\tilde p_N^*}_*^2\sim \vertiii{\epsilon \nabla v-\tilde y_N^*}_*^2.
\end{align}

Note that in practice, the term $2\int\limits_{\Omega}{(\tilde y_N^*-\tilde p_N^*)\cdot \nabla (v-\tilde u^N) dx}$ is much smaller than $\vertiii{\nabla (v-\tilde u^N)}^2+\vertiii{\tilde y_N^*-\tilde p_N^*}_*^2$.
For the 2-term regularization, instead of $\tilde u_{h_L}^L$ in \eqref{Explicit_form_DF_v_minus_lambda_yStar}, we have $G+u_{h_L}^L$, where $u_{h_L}^L$ is a conforming (not neccessarily finite element) approximation of $u^L$.

We end this section by recalling a near best approximation result (\cite{Kraus_Nakov_Repin_PBE1_2018} ). Contrary to the result in \cite[Theorem 6.2]{Chen2006b},
we do not make any restrictive assumptions on the meshes to ensure that the finite element approximations $\tilde u_{h_N}^N$ are uniformly bounded in the $L^\infty$ norm.  
Let $V_h\subset L^\infty(\Omega)$ be a closed subspace of $H_0^1(\Omega)$ and $\tilde u_{h_N}^N\in V_h$ be the Galerkin approximation of $\tilde u^N$ defined by:
\begin{align}\label{un_Galerkin_formulation}
\begin{aligned}
&\text{ Find } \tilde u_{h_N}^N\in V_h, \text{ such that }\\
&a(\tilde u_{h_N}^N,v)+\int_{\Omega}{b(x,\tilde u_{h_N}^N+\tilde u_{h_L}^L)v dx}=0,\,\forall v\in V_h
\end{aligned}
\end{align}
\begin{proposition}\label{Proposition_quasi_optimal_a_priori_error_estimate_1_3_term_uN}
Let $V_h\subset L^\infty(\Omega)$ be a closed subspace of $H_0^1(\Omega)$ and $\tilde u_{h_N}^N\in V_h$ be the Galerkin approximation of $\tilde u^N$ defined by \eqref{un_Galerkin_formulation}. Then
\begin{align*}
\vertiii{\nabla (\tilde u_{h_N}^N-\tilde u^N)}^2\leq \inf\limits_{v\in V_h}{\bigg\{\vertiii{\nabla (v-\tilde u^N)}^2+\int\limits_{\Omega_s}{k^2\left(\sinh(v+\tilde u_{h_L}^L)-\sinh(\tilde u^N+\tilde u_{h_L}^L)\right)^2 dx}\bigg\}}
\end{align*}
\end{proposition}
If $V_h$ is a finite element space, then using Proposition \ref{Proposition_quasi_optimal_a_priori_error_estimate_1_3_term_uN}, qualified and unqualified convergence of the finite element approximations $\tilde u_{h_N}^N$ can be proven since $\tilde u^N, \tilde u_{h_L}^L\in L^\infty(\Omega)$ and $\sinh$ is a locally Lipschitz function (see \cite{Kraus_Nakov_Repin_PBE1_2018}). Of course, an analogous result holds also for the component $\tilde u^N$ in the 2-term splitting.

\subsection{Nonlinear nonhomogeneous problem for $u$ in the 3-term splitting without performing additional splitting into $u^L+u^N$}\label{Section_estimate_for_uRegular_3_term}

As we mentioned in Section~\ref{Section_Variational_form_of_the_problem}, it would be better if we could estimate directly the error when solving \eqref{u_weak_formulation}. Here $v$ denotes some conforming approximation of $\tilde u$ - the weak solution of problem \eqref{u_weak_formulation} where instead of $\nabla u^H$ we have $T(\nabla \tilde u^H)$. The function $\tilde u^H$ is some conforming approximation of $u^H$ and $T$ is some operator that maps the numerical flux $\epsilon \nabla \tilde u^H$ into $H(\div;\Omega_m)$ (if $\nabla \tilde u^H$ is already in $H(\div;\Omega_m)$, then we can take $T$ to be the identity). By $\tilde y^*\in \left[L^2(\Omega)\right]^3$ we denote an arbitrary approximation of the exact flux $\tilde p^*=\epsilon\nabla\tilde u$. 
We briefly discuss the derivation of a functional error estimate similar to the one derived in \cite{Kraus_Nakov_Repin_PBE1_2018}. For this, we consider only the case of homogeneous boundary conditions $g=0$, which correspond to $\psi=0$ in \eqref{PBE_special_form_4}, and note that the case of nonhomogeneous boundary conditions can be easily treated (see \cite{Repin}). By $\tilde J$ we denote the functional $J$ in \eqref{definition_of_J_Regular_Part_3_term} but with $T(\nabla\tilde u^H)$ instead of $\nabla u^H$ in its definition. We rewrite the functional $\tilde J$ in the general form $\tilde J=G(\Lambda v)+F(v)$, where 
\begin{align*}
&G(\Lambda v):=\int\limits_{\Omega}{\frac{\epsilon}{2}\abs{\nabla v}^2dx},\\
& F(v):=\int\limits_{\Omega}{k^2 \cosh(v)dx}+\int\limits_{\Omega_m}{\epsilon_m T(\nabla \tilde u^H)\cdot\nabla v dx}-\int\limits_{\Omega_s}{\epsilon_m\nabla G\cdot\nabla v dx}
\end{align*}
and $\Lambda:=\nabla: H_0^1(\Omega)\to \left[L^2(\Omega)\right]^3$. We further denote by $\Lambda^*:=-\div: \left[L^2(\Omega)\right]^3 \to H^{-1}(\Omega)$ the adjoint operator to $\Lambda$, where $H^{-1}(\Omega)$ denotes the dual space of $H_0^1(\Omega)$ and with $\langle \cdot , \cdot \rangle$ the duality pairing in $H^{-1}(\Omega)\times H_0^1(\Omega)$. In order to apply the abstract framework from \cite{Kraus_Nakov_Repin_PBE1_2018}, which is based on the theory in \cite{Repin}, we compute the Fenchel conjugate of $F$ evaluated at $-\Lambda^*\tilde y^*$ for $\tilde y^*\in \left[L^2(\Omega)\right]^3$. Since the exact flux $\tilde p^*=\epsilon\nabla \tilde u$ satisfies the prescribed interface jump condition $\tilde g_\Gamma$ by the functions $\epsilon_m T(\nabla \tilde u^H)$ and $\epsilon_m\nabla G$, i.e
\begin{align}\label{definition_of_functional_tilde_gGamma_3term_regularization}
\langle \tilde g_\Gamma,v\rangle_\Gamma = -\langle \gamma_{n,\Omega_m} \left(\epsilon_m T\left(\nabla\tilde u^H\right)\right),v\rangle_\Gamma + \langle \gamma_{n,\Omega_s}(\epsilon_m\nabla G),v\rangle_\Gamma,\,\forall v\in H^{1/2}(\Gamma),
\end{align}
it is enough to compute $F^*(-\Lambda^*\tilde y^*)$ only for such $\tilde y^*\in \left[L^2(\Omega)\right]^3$ that can be represented in the form 
\begin{align}\label{special_form_yStar_uRegular}
\tilde y^*=-\epsilon_m T(\nabla\tilde u^H) \mathbbm{1}_{\Omega_m}+\epsilon_m\nabla G\mathbbm{1}_{\Omega_s}+\tilde y_0^*, \text{ for } \tilde y_0^*\in H(\div;\Omega).
\end{align}
 
\begin{align}\label{Computing_F_star}
\begin{aligned}
&F^*(-\Lambda^*\tilde y^*)=\sup\limits_{z\in H_0^1(\Omega)}{\left[\langle -\Lambda^* \tilde y^*,z\rangle - F(z)\right]}=
\sup\limits_{z\in H_0^1(\Omega)}{\left[( -\tilde y^*,\Lambda z) - F(z)\right]}\\
&=\sup\limits_{z\in H_0^1(\Omega)}{\int\limits_{\Omega}{\left[-\tilde y^*\cdot\nabla z-k^2\cosh(z)-\epsilon_m T(\nabla\tilde u^H) \mathbbm{1}_{\Omega_m}+\epsilon_m\nabla G\mathbbm{1}_{\Omega_s}\right]dx}}\\
&=\sup\limits_{z\in H_0^1(\Omega)}{\int\limits_{\Omega}{\left[-\tilde y_0^*\cdot\nabla z-k^2\cosh(z)\right]dx }}\\
&=\sup\limits_{z\in H_0^1(\Omega)}{\int\limits_{\Omega}{\left[\div\tilde y_0^* z-k^2\cosh(z)\right]dx }}\quad\left(\text{ finite if }\div \tilde y_0^*=0 \text{ in } \Omega_m\right)\\
&\leq \int\limits_{\Omega_s}{\sup\limits_{\xi\in\mathbb R}{\left(\div \tilde y_0^*(x)\xi-k^2\cosh(\xi)\right) dx}}\\
&=\int\limits_{\Omega_s}{\left(\div \tilde y_0^*\arsinh\left(\frac{\div \tilde y_0^*}{k^2}\right) - k^2\sqrt{\left(\frac{\div\tilde y_0^*}{k^2}\right)^2+1}\right)dx}.
\end{aligned}
\end{align}
We note that we actually have equalities everywhere in \eqref{Computing_F_star}. The proof for this is similar to the one in \cite{Kraus_Nakov_Repin_PBE1_2018} and we omit it here. We define the majorant $M_{\oplus}^2(v,\tilde y^*)$ for any $\tilde y^*$ of the form \eqref{special_form_yStar_uRegular} with $\div \tilde y_0^*=0$ in $\Omega_m$ by
\begin{align}\label{Explicit_form_of_uRegular_Majorant}
2M_{\oplus}^2(v,\tilde y^*)=\int\limits_{\Omega}{\tilde \eta^2(x)dx}=\vertiii{\epsilon\nabla v-\tilde y^*}_*^2+2D_F(v,-\Lambda^*\tilde y^*),
\end{align}
where 
\begin{eqnarray}
D_F(v,-\Lambda^*\tilde y^*)&=&F(v)+F^*(-\Lambda^*\tilde y^*) +\langle \Lambda^*\tilde y^*, v\rangle \nonumber \\
&=&\int\limits_{\Omega_s}{k^2\left(\cosh(v)+\frac{\div \tilde y_0^*}{k^2}\arsinh\left(\frac{\div \tilde y_0^*}{k^2}\right)-\sqrt{\left(\frac{\div \tilde y_0^*}{k^2}\right)^2+1}\right)dx} \nonumber \\
&-&\int\limits_{\Omega_s}{\div \tilde y_0^*v dx}. \label{Explicit_form_DF_v_minus_lambda_yStar_3_term_uRegular}
\end{eqnarray}
 The error estimate for the combined energy norm can be expressed as
\begin{align}\label{error_estimate_CEN_uRegular}
&\frac{1}{2}\vertiii{\epsilon \nabla v-\tilde y^*}_*^2\leq\vertiii{\nabla(\tilde u-v)}^2+\vertiii{\tilde p^*-\tilde y^*}_*^2\leq 2M_{\oplus}^2(v,\tilde y^*) .
\end{align}
To see that this estimate is sharp, we take $\tilde y_0^*:=\tilde p^*+\epsilon_m T(\nabla\tilde u^H) \mathbbm{1}_{\Omega_m}-\epsilon_m\nabla G\mathbbm{1}_{\Omega_s}$. It is easy to verify that this $\tilde y_0^*$ is in $H(\div;\Omega)$ with $\div y_0^*=0$ in $\Omega_m$. Then the corresponding $\tilde y^*$ is equal to $\tilde p^*$ and clearly $M_\oplus(v,\tilde y^*)^2 =J(\tilde u)- I^*(\tilde p^*)=0$, where $I^*(\tilde y^*)=-G^*(\tilde y^*)-F^*(-\Lambda^*\tilde y^*)$ and $G^*$ is the Fenchel conjugate of $G$ (see \cite{Repin, Kraus_Nakov_Repin_PBE1_2018}). 

Similarly to the near best approximation result for $\tilde u^N$ (Proposition \ref{Proposition_quasi_optimal_a_priori_error_estimate_1_3_term_uN}), we present such a result also for the regular
component $\tilde u$. Let $V_h\subset L^\infty(\Omega)$ be again a closed subspace of $H_0^1(\Omega)$ and let $\tilde u_h$ be the unique minimizer of $\tilde J$ over $V_h$, which is also
the unique solution of the Galerkin problem:
\begin{align}\label{Discrete_Weak_Formulation_uRegular_3_term}
&\text{Find } \tilde u_h\in V_h \text{ such that }\nonumber\\
&a(\tilde u_h,v)+\int\limits_{\Omega}{b(x,\tilde u_h)v dx}=-\int\limits_{\Omega_m}{\epsilon_m T(\nabla \tilde u^H)\cdot\nabla v dx}+\int\limits_{\Omega_s}{\epsilon_m\nabla G\cdot\nabla v dx},\\
&\text{ for all } v\in V_h.\nonumber
\end{align}
Now, if we denote $D_G(\Lambda v, \tilde y^*):=G(\Lambda v)+G^*(\tilde y^*)-(\tilde y^*,\Lambda v)\ge 0$ (see \cite{Repin}), using the abstract framework presented in \cite{Repin, Kraus_Nakov_Repin_PBE1_2018} we can write for any $v\in H_0^1(\Omega)$
\begin{align}\label{primal_part_of_error} 
&\tilde J(v)-\tilde J(\tilde u)=M_\oplus^2(v,\tilde p^*)=D_G(\Lambda v,\tilde p^*)+D_F(v,-\Lambda^*\tilde p^*).
\end{align}
Then, using \eqref{primal_part_of_error} and that $D_G(\Lambda v,\tilde p^*)=\frac{1}{2}\vertiii{\nabla (v-\tilde u)}^2$, for any $v \in V_h$ we can write
\begin{align*}
&\mathrel{\phantom{=}}\vertiii{\nabla (\tilde u_h-\tilde u)}^2+2D_F(\tilde u_h,-\Lambda^*\tilde p^*)=2\left(\tilde J(\tilde u_h)-\tilde J(\tilde u)\right)\\
&\leq 2\left(\tilde J(v)-\tilde J(\tilde u)\right)=\vertiii{\nabla (v-\tilde u)}^2+2D_F(v,-\Lambda^*\tilde p^*).
\end{align*}
Next, using $\div \tilde p_0^*=k^2\sinh(\tilde u)$, where $\tilde p_0^*=\tilde p^*+\epsilon_m T(\nabla\tilde u^H) \mathbbm{1}_{\Omega_m}-\epsilon_m\nabla G\mathbbm{1}_{\Omega_s}$, we calculate
\begin{align}
D_F(v,-\Lambda^*\tilde p^*)=\int\limits_{\Omega_s}{k^2\left(\cosh(v)-\cosh(\tilde u)+\tilde u\sinh(\tilde u)-v\sinh(\tilde u)\right)dx}.
\end{align}
Using Proposition 3.2 in \cite{Kraus_Nakov_Repin_PBE1_2018} and the fact that $2D_F(\tilde u_h,-\Lambda^*\tilde p^*)\ge 0$, we obtain the near best approximation
result for the regular component $\tilde u$.
\begin{proposition}\label{Proposition_quasi_optimal_a_priori_error_estimate_1_3_term_uRegular}
Let $V_h\subset L^\infty(\Omega)$ be a closed subspace of $H_0^1(\Omega)$ and $\tilde u_h\in V_h$ be the Galerkin approximation of $\tilde u$ defined by \eqref{Discrete_Weak_Formulation_uRegular_3_term}. Then
\begin{align}\label{quasi_optimal_a_priori_error_estimate_1}
\vertiii{\nabla (\tilde u_h-\tilde u)}^2\leq \inf\limits_{v\in V_h}{\bigg\{\vertiii{\nabla (v-\tilde u)}^2+\int\limits_{\Omega_s}{k^2(\sinh(v)-\sinh(\tilde u))^2 dx}\bigg\}}
\end{align}
\end{proposition}
If $V_h$ is a finite element space, then using Proposition \ref{Proposition_quasi_optimal_a_priori_error_estimate_1_3_term_uRegular}, qualified and unqualified convergence of the finite element approximations $\tilde u_h$ can be proven since $\tilde u\in L^\infty(\Omega)$ and $\sinh$ is a locally Lipschitz function (see \cite{Kraus_Nakov_Repin_PBE1_2018}).

\subsection{Overall error in the regular component $u$}

\subsubsection{Overall error in the regular component $u$ with the additional splitting in $u^L+u^N$}

Finally, we give a justification for how the a posteriori error estimates developed in Sections \ref{Section_estimate_for_harmonic_part}, \ref{Section_estimate_for_linear_part}, \ref{Section_estimate_for_nonlinear_part} can be put together and applied in practice. More precisely, we estimate the effect of using an approximation of $u^H$ in \eqref{definition_ul_weak_formulation_with_functional_gGamma_3term_regularization} to compute an approximation of $u^L$ which is then used in equation \eqref{un_weak_formulation} to compute an approximation of $u^N$ on the quality of the total approximate regular part of the potential $\tilde u_h\sim u$. Again, by $\tilde u^L$ we denote the exact solution to \eqref{tildeul_weak_formulation}, where $\tilde u^H$ is a conforming approximation of $u^H$ and $T(\nabla\tilde u^H)\in H(\div;\Omega)$. The finite element approximation of $\tilde u^L$ is denoted by $\tilde u_{h_L}^L$ and is assumed to be computed by some conforming FEM based on the weak formulation \eqref{tildeul_weak_formulation} on some mesh for which we use a subindex ${h_L}$ to distinguish the finite element functions corresponding to this mesh. This means that $\tilde u_{h_L}^L\in H_g^1(\Omega)$ and it can be regarded as a conforming approximation of both $u^L$ and $\tilde u^L$. With $u^N$ as before, we denote the exact solution to equation \eqref{un_weak_formulation} with the exact $u^L$ in it. Again, by $\tilde u^N$ we denote the exact solution to equation \eqref{un_weak_formulation} but with $\tilde u_{h_L}^L$ in it, and by $\tilde u_{h_N}^N$ we denote the conforming finite element approximation of $\tilde u^N$. The subindex $h_N$ in $\tilde u_{h_N}^N$ again means that this approximation is computed on a possibly different mesh than the one used for $\tilde u_{h_L}^L$. In this notation, the approximation $\tilde u_h$ of $u$ that we compute satisfies $\tilde u_h=\tilde u_{h_L}^L+\tilde u_{h_N}^N$. 
We want to estimate $|||\nabla(u-\tilde u_h)|||$.
\begin{eqnarray}
\vertiii{\nabla(u-\tilde u_h)}&=&\vertiii{\nabla(u^N+u^L-\tilde u_{h_N}^N-\tilde u_{h_L}^L)} \nonumber \\
&\leq& \vertiii{\nabla(u^N-\tilde u_{h_N}^N)}+\vertiii{\nabla(u^L-\tilde u_{h_L}^L)} \label{triangle_inequality_on_u_minus_tildeu_EnergyNorm}
\end{eqnarray}
For the first term on the RHS we have that
\begin{align}\label{estimate_on_the_error_in_computing_un_EnergyNorm}
\vertiii{\nabla(u^N-\tilde u_{h_N}^N)}\leq \vertiii{\nabla(u^N-\tilde u^N)}+\vertiii{\nabla(\tilde u^N-\tilde u_{h_N}^N)}.
\end{align}
The second term on the RHS in \eqref{estimate_on_the_error_in_computing_un_EnergyNorm} we estimate by the functional a posteriori error estimate \eqref{Explicit_form_of_Functional_error_equality} and the first term we estimate as follows:
\begin{align}\label{the_equations_that_un_and_tilde_un_satisfy}
\begin{array}{|ll}
&a(u^N,v)+(b(x,u^N+u^L),v)=0,\,\forall v\in H_0^1(\Omega)\\
&a(\tilde u^N,v)+(b(x,\tilde u^N+\tilde u_{h_L}^L),v)=0,\,\forall v\in H_0^1(\Omega)
\end{array}
\end{align}
Subtracting the second from the first equation above, we get
\begin{align}\label{weak_equation_estimating_uN_minus_tildeuN}
&a(u^N-\tilde u^N,v)=(b(x,\tilde u^N+\tilde u_{h_L}^L)-b(x,u^N+u^L),v)=0,\,\forall v\in H_0^1(\Omega).
\end{align}
Now take $v:=u^N+u^L-\tilde u^N-\tilde u_{h_L}^L\in H_0^1(\Omega)$ and obtain 
%If $\Gamma$ is $C^1$ then $u^L\in L^\infty(\Omega)$ and from Proposition \ref{Proposition_1} we know that also   $u^N\in L^\infty(\Omega)$. Thus $v\in L^\infty(\Omega)$ and we can test with it in \eqref{weak_equation_estimating_uN_minus_tildeuN}. If $\Gamma$ is just Lipschitz then, $v$ is not necessarily in $L^\infty(\Omega)$. In this case, we can apply the result of Br{\'e}zis and Browder in \cite{Brezis_Browder_1978_One_Property_of_Sobolev_Spaces} to obtain that
\begin{equation}\label{inequality_using_monotonicity_of_b}
\begin{split}
&\mathrel{\phantom{=}} a(u^N-\tilde u^N,u^N+u^L-\tilde u^N-\tilde u_{h_L}^L)\\
&=(b(x,\tilde u^N+u_{h_L}^L)-b(x,u^N+u^L),u^N+u^L-\tilde u^N-\tilde u_{h_L}^L)\leq 0
\end{split}
\end{equation}
where we have used the monotonicity of the nonlinearity: $(b(x,w)-b(x,z),w-z)\ge 0$, $\forall w,z\in H^1(\Omega)$. Using the boundedness of the bilinear form $a(.,.)$ we get 
\begin{align*}
&\mathrel{\phantom{=}} \vertiii{\nabla(u^N-\tilde u^N)}^2= a(u^N-\tilde u^N,u^N-\tilde u^N)\\
&= a(u^N-\tilde u^N,u^N+u^L-\tilde u^N-\tilde u_{h_L}^L)+a(u^N-\tilde u^N,\tilde u_{h_L}^L-u^L)\\
&\leq  0+\vertiii{\nabla(u^N-\tilde u^N)}\vertiii{\nabla(u^L-\tilde u_{h_L}^L)} .
\end{align*}
Thus,
\begin{align}\label{the_estimate_that_un_minus_tildeun_satisfies_EnergyNorm}
\vertiii{\nabla(u^N-\tilde u^N)}\leq \vertiii{\nabla(u^L-\tilde u_{h_L}^L)} .
\end{align}
Next, we estimate the term $\vertiii{\nabla(u^L-\tilde u_{h_L}^L)}$. By the triangle inequality, we have
\begin{align}\label{triangle_inequality_on_ul_minus_tildeu_hll_EnergyNorm}
\vertiii{\nabla(u^L-\tilde u_{h_L}^L)}\leq \vertiii{\nabla(u^L-\tilde u^L)}+\vertiii{\nabla(\tilde u^L-\tilde u_{h_L}^L)},
\end{align}
where to the second term we apply the a posteriori error estimate \eqref{A_posteriori_error_estimate_for_tildeuL_3term_regularization} or \eqref{A_posteriori_error_estimate_for_tildeuL_3term_regularization_squared_Majorant} for problem \eqref{tildeul_weak_formulation} and the first term is bounded as follows: subtract equation \eqref{tildeul_weak_formulation} from \eqref{definition_ul_weak_formulation_with_functional_gGamma_3term_regularization}, take $\phi=u^L-\tilde u^L$ and use Cauchy-Schwartz inequality to obtain
\begin{align*}
&\vertiii{\nabla(u^L-\tilde u^L)}^2\leq \sqrt{\epsilon_m}\|\nabla u^H-T(\nabla\tilde u^H)\|_{L^2(\Omega_m)}\vertiii{\nabla(u^L-\tilde u^L)}
\end{align*}
Thus, we get
\begin{align}\label{estimate_on_ul_minus_tildeul_EnergyNorm}
&\vertiii{\nabla(u^L-\tilde u^L)}\leq \sqrt{\epsilon_m}\|\nabla u^H-T(\nabla\tilde u^H)\|_{L^2(\Omega_m)}\leq \sqrt{\epsilon_m}M_{\oplus,H}\left(\tilde u^H,T(\nabla \tilde u^H)\right).
\end{align}
Finally, if we want to compute an approximation $\tilde u$ of $u$ with a prescribed error tolerance $\delta$, using \eqref{triangle_inequality_on_u_minus_tildeu_EnergyNorm}, \eqref{estimate_on_the_error_in_computing_un_EnergyNorm},  \eqref{the_estimate_that_un_minus_tildeun_satisfies_EnergyNorm}, \eqref{triangle_inequality_on_ul_minus_tildeu_hll_EnergyNorm}, and \eqref{estimate_on_ul_minus_tildeul_EnergyNorm} we obtain
\begin{align}\label{final_form_of_the_estimate_for_u_minus_tildeu_EnergyNorm}
\begin{aligned}
&\mathrel{\phantom{=}}\vertiii{\nabla(u-\tilde u_h)}\leq \vertiii{\nabla(u^L-\tilde u_{h_L}^L)}+\vertiii{\nabla(\tilde u^N-\tilde u_{h_N}^N)}+\vertiii{\nabla(u^L-\tilde u_{h_L}^L)}\\
&\leq 2\left(\vertiii{\nabla(u^L-\tilde u^L)}+\vertiii{\nabla(\tilde u^L-\tilde u_{h_L}^L)} \right)+\sqrt{2}M_{\oplus,N}(\tilde u_{h_N}^N,\tilde y_N^*)\\
&\leq 2\left(\sqrt{\epsilon_m}\|\nabla u^H-T(\nabla \tilde u^H)\|_{L^2(\Omega_m)}+M_{\oplus,L}(\tilde u^L_{h_L},y_0^*)\right)+\sqrt{2} M_{\oplus,N}(\tilde u_{h_N}^N,\tilde y_N^*)\\
&\leq 2\sqrt{\epsilon_m}M_{\oplus,H}\left(\tilde u^H,T(\nabla \tilde u^H)\right)+2M_{\oplus,L}(\tilde u^L_{h_L},y_0^*)+\sqrt{2} M_{\oplus,N}(\tilde u_{h_N}^N,\tilde y_N^*)\leq \delta .
\end{aligned}
\end{align}
For the $2$-term regularization, we have
\begin{align}\label{the_estimate_on_u_minus_tildeu_2term_regularization}
\vertiii{\nabla(u-\tilde u_h)}\leq 2M_{\oplus,L}(u_{h_L}^L,y_0^*)+\sqrt{2}M_{\oplus,N}(\tilde u_{h_N}^N,\tilde y_N^*)\leq \delta,
\end{align}
where $\tilde u_h=u_{h_L}^L+\tilde u_{h_N}^N$, $u_{h_L}^L$ denotes a conforming approximation of $u^L$ and $\tilde u_{h_N}^N$ a conforming approximation of $\tilde u^N$--the exact solution of \eqref{2_term_regularization_un_weak_formulation} containing $u_{h_L}^L$. Note that we can estimate the $H^1$-seminorm and the full $H^1$-norm of the
difference $u-\tilde u_h$ by introducing Friedrichs constant and the minimum and maximum values of the dielectric coefficient $\epsilon$.
\begin{remark}
We recall that the conforming (FEM) approximations of $u^L$ and $\tilde u^N$ in the 2-term splitting are denoted by $u_{h_L}^L$ and $\tilde u_{h_N}^N$, respectively,
and the conforming approximtions of $\tilde u^L$ and $\tilde u^N$ in the 3-term splitting by $\tilde u^L_{h_L}$ and $\tilde u^N_{h_N}$, respectively. To be more precise,
since the Dirichlet boundary condition (BC) on  $\tilde u^L$ is $\gamma_{\partial\Omega}(\tilde u^L)=g$ for the $3$-term regularization and
$\gamma_{\partial\Omega}(u^L)=g-G$ for the $2$-term regularization, it is clear that $g$ and $g-G$ are not exactly representable in most finite element spaces.
Thus, by a conforming FEM approximation $\tilde u^L_{h_L}$ we should understand $\tilde u^L_{h_L}=u^L_g+\tilde u^L_{0,h_L}$ for the $3$-term regularization
and $u^L_{h_L}=u^L_{g-G}+u^L_{0,h_L}$ for the $2$-term regularization, where $\tilde u^L_{0,h_L}$ solves a discrete FEM form of the homogenized equation
$a(\tilde u_0^L,v)=\langle g_\Gamma,\gamma_\Gamma(v)\rangle_\Gamma-a(u_g^L,v)$ for the $3$-term splitting and
$a(u_0^L,v)=\langle g_\Gamma,\gamma_\Gamma(v)\rangle_\Gamma-a(u_{g-G}^L,v)$ for the $2$-term splitting, respectively ($u_g^L$ and $u_{g-G}^L$ are defined
in Section~\ref{Section_Variational_form_of_the_problem}). 

We further recall that for the $2$-term splitting there is no approximation error in the interface condition and
so we use $u^L$ and $u^L_{h_L}$ instead of $\tilde u^L$ and $\tilde u^L_{h_L}$.
In the 2-term and 3-term splitting, the Dirichlet BC on $\tilde u^N$ is homogeneous so there is no problem with building a conforming FEM approximation $\tilde u^N_{h_N}$.
To build a conforming FEM approximation $\tilde u^H$ for $u^H$ we do the same as above.
\end{remark}
\subsubsection{Overall error in the regular component $u$ without the additional splitting in $u^L+u^N$}
Here, we estimate the overall error in the regular component $u$ for the 3-term splitting in the case where we do not perform the additional splitting in $u^L+u^N$.
We need to estimate the error $\vertiii{\nabla(u-\tilde u_h)}$, where $\tilde u_h$ is a conforming approximation of $\tilde u$, the solution of
problem~\eqref{u_weak_formulation}, where instead of $\nabla u^H$ we have $T(\nabla \tilde u^H)$. By the triangle inequality, we have
\begin{align}\label{triangle_inequality_u_minus_tildeu_h}
&\vertiii{\nabla(u-\tilde u_h)}\leq \vertiii{\nabla(u-\tilde u)}+\vertiii{\nabla(\tilde u-\tilde u_h)}.
\end{align}
The second term is estimated by \eqref{error_estimate_CEN_uRegular}. For the first term, after subtracting the equation 
\begin{align}\label{tildeu_weak_formulation}
\begin{aligned}
&a(\tilde u,v)+\int_{\Omega}{b(x,\tilde u)v dx}=\\
&-\int\limits_{\Omega_m}{\epsilon_m T \left(\nabla \tilde u^H\right)\cdot\nabla v dx}+\int\limits_{\Omega_s}{\epsilon_m\nabla G\cdot\nabla v dx}, \text{ for all } v\in H_0^1(\Omega) ,
\end{aligned}
\end{align}
(that comes from~\eqref{u_weak_formulation}), we obtain
\begin{align*}
&a(u-\tilde u,v)=\left(b(x,\tilde u)-b(x,u)),v\right)+\int\limits_{\Omega_m}{\epsilon_m \left( T\left(\nabla\tilde u^H\right)-\nabla u^H\right)\cdot \nabla v dx}, \text{ for all } v\in H_0^1(\Omega).
\end{align*}
Set here $v:=u-\tilde u\in H_0^1(\Omega)$. Using the monotonicity of $b(x,\cdot)$, we see that
\begin{align}
\vertiii{\nabla(u-\tilde u)}\leq \sqrt{\epsilon_m}\|T\left(\nabla \tilde u^H\right)-\nabla  u^H\|_{L^2(\Omega_m)}.
\end{align}
The overall error estimate for $u$ is as follows:
\begin{align}\label{overal_error_estiamte_for_uRegular_3_term_no_splitting_in_uLplusuN}
\vertiii{\nabla(u-\tilde u_h)}\leq \sqrt{\epsilon_m}M_{\oplus,H}\left(\tilde u^H,T(\nabla \tilde u^H)\right)+ \sqrt{2}M_\oplus(\tilde u_h,\tilde y^*)
\end{align}
\begin{remark}
Note that we can skip the condition that the test functions $v$ in \eqref{tildeu_weak_formulation} and \eqref{u_weak_formulation} are in $L^\infty(\Omega)$ because $u,\tilde u\in L^\infty(\Omega)$. We 
have already seen in Section \ref{Section_Variational_form_of_the_problem} that $u$ is in $L^\infty(\Omega)$. To see that $\tilde u$ is also in $L^\infty(\Omega)$ we can formally split $\tilde u$ in a component that solves a linear nonhomogeneous interface problem and a component that solves a nonlinear homogeneous problem. The first component will be in $L^\infty$ by virtue of Theorem B.2 in \cite{Kinderlehrer_Stampacchia} provided that $T\left(\nabla\tilde u^H\right) \in \left[L^s(\Omega)\right]^d$ for some $s>d$. The second component that solves the nonlinear homogeneous problem, is in $L^\infty(\Omega)$ by virtue of Propostion \ref{Proposition_1}.
\end{remark}

\section{Numerical results}\label{Section_Numerical_Results}

In this section we present three numerical examples based on the two term and three term regularizations. They show that nonlinear mathematical models in question can be studied by fully reliable computer simulation methods that provide results with guaranteed and explicitly known accuracy. In the first and third examples, we consider the system of two chromophores Alexa 488 and Alexa 594, which are used for protein labelling in biophysical experiments. The second experiment is conducted on an insulin protein (PDB ID: 1RWE). In all examples, we assume a solution consisting of NaCl with $k^2=10\, \AA^{-2}$, corresponding to ionic strength of $I_s\approx1.178$ molar. The ground state charges are obtained by CHARMM 32. In the first and third examples, we assume dielectric constnats $\epsilon_m=2,\,\epsilon_s=80$ and in the second example $\epsilon_m=20,\,\epsilon_s=80$. The numerical experiments are carried out in FreeFem++ developed and maintained by Frederich
Hecht~\cite{FreeFem}.
All Figures below are generated with the help of VisIt~\cite{VisIt}.
The computational domain $\Omega$ for all examples is a cube with a side lenght of $A=6 a_{\max}+24\, \AA$, where $a_\max$ is the maximum side length of the smallest bounding box for the molecule(s) with edges parallel to the coordinate axes. This amounts to $A=295.85\, \AA$ for the first and third examples and $A=268.21 \,\AA$ for the second one. The molecules are situated in the center of $\Omega$.  For the Friedrichs' constant $C_{F\Omega}$ on a cube we have that $C_{F\Omega}\leq \frac{A\sqrt{3}}{3\pi}$ (see \cite{Mikhlin_1986_book}). The Dirichlet boundary condition in \eqref{PBE_special_form} for all experiments is given by $\psi=0$ on $\partial\Omega$. The discretization used in the numerical tests to find conforming approximations $u_{h_L}^L$ and $\tilde u_{h_N}^N$ for $u^L$ and $\tilde u^N$, respectively, is based on standard linear ($P_1$) finite elements although the derived estimates apply to any conforming approximations, which could for example also be obtained from higher order finite element methods (hpFEM) or isogeometric analysis (IGA). The surface meshes are constructed with TMSmesh 2.1 \cite{Chen_Tu_Lu_2012, Chen_Liu_Lu_2018} which produces a Gaussian molecular surface. The surface mesh of the two chromophores is additionally optimized with the help of Mmgs \cite{mmg_Dapogny_Dobrzynski_Frey_2013} and the surface mesh of the insulin protein is optimised with MeshLab \cite{MeshLabGeneral}. %Mmgs \cite{mmg_Dapogny_Dobrzynski_Frey_2013}. 
The initial tetrahedral meshes are generated using TetGen \cite{TetGen} and then they are adapted with the help of mmg3d \cite{mmg3d_dobrzynski}. The shape of the molecules is not changed during adaptation. This is justified, since the molecule structure is only known with a certain precision from X-ray crystallography. It is also possible to use isoparametric elements to represent the molecular surface exactly. Then, in the mesh refining procedure new points will be inserted on the surface by splitting the curved elements on the interface~$\Gamma$.

\subsection{Example 1: Alexa 488 and Alexa 594}

The first system consists of two chromophores Alexa 594 and Alexa 488 with a total of 171 atoms in aqueous solution of NaCl. The parameters of the force fields of Alexa chromophores were created by an analogy approach from that of similar chemical groups in the CHARMM forse field (version v35b3). The coordinates of the molecules are taken from a time frame of molecular dynamic simulations. In the all-atom MD simulations the dyes were attached to a polyproline 11 and dissolved in water box with NaCl \cite{Sobakinskaya_Busch_Renger_2018}. 
The parameters of this example are $\epsilon_m=2$, $\epsilon_s=80$, $k_m^2=0\,\AA^{-2}$,  $k_s^2=10\, \AA^{-2}$ which corresponds to ionic strength $I_s=1.178$ molar.

\subsubsection{Finding $u^L$}\label{SubSubsection_2_term_Solving_for_uL_2DYES}

First, we solve adaptively \eqref{PBE_2term_regular_linear_part} to find an approximation $u_{h_L}^L$ of $u^L$. As an error indicator we use the second term in the error estimate computed over each element $K$ \eqref{A_posteriori_error_estimate_for_tildeuL_2term_regularization} 
\begin{align}\label{2_term_regularization_error_indicator_Linear_Part}
\eta_K^L=\vertiii{\epsilon\nabla u_{h_L}^L-y_{L}^*}_{*(K)}=\left(\int\limits_{K}{\frac{1}{\epsilon}\abs{\epsilon\nabla u_{h_L}^L-y_{L}^*}^2 dx}\right)^{\frac{1}{2}},
\end{align}
where $y_{L}^*=(\epsilon_m-\epsilon)\nabla G+y_0^*$.
To find $y_0^*\in H(\div;\Omega)$, we perform a minimization of the squared majorant $M_{\oplus,L}^2(u_{h_L}^L,y_0^*;\alpha)$ over $\alpha\in\mathbb R_{>0}$ and $y_0^*\in RT_0$ defined over the same mesh.
\begin{eqnarray}
\vertiii{\nabla(u^L-u_{h_L}^L)}^2 &\leq& M_{\oplus,L}^2(u_{h_L}^L,y_0^*) \nonumber\\
&\leq& (1+\alpha)\frac{C_{F\Omega}}{\epsilon_\min}\|\div y_0^*\|_{L^2(\Omega)}^2+\frac{(1+\alpha)}{\alpha}\vertiii{\epsilon\nabla u_{h_L}^L-(\epsilon_m-\epsilon)\nabla G-y_0^*}_*^2
\nonumber \\
&:=&M_{\oplus,L}^2(u_{h_L}^L,y_0^*;\alpha) \label{A_posteriori_error_estimate_for_uL_2term_regularization_squared_Majorant}
\end{eqnarray}
This procedure gives a very sharp bound from above for the error. Moreover, we have a simple and efficient lower bound for the energy norm of the error $\nabla (u^L-u_{h_L}^L)$. Indeed, let us denote by $J^L$ the quadratic functional whose unique minimizer over $H_{g-G}^1(\Omega)$ is the solution $u^L$ of \eqref{definition_ul_weak_formulation_with_functional_gGamma_2term_regularization} and which is defined by $J^L(v)=\int\limits_{\Omega}{\left(\frac{\epsilon}{2}\abs{\nabla u^L}-(\epsilon_m-\epsilon)\nabla G\cdot\nabla v\right)dx}$. Then, assuming that $u_{h_L}^L\in H_{g-G}^1(\Omega)$ (for example when $u_{h_L}^L=u_{g-G}^L+u_{0,h_L}^L$, where $u_{0,h_L}^L$ is the finite element solution of the homogenized version of \eqref{definition_ul_weak_formulation_with_functional_gGamma_2term_regularization} and $u_{g-G}^L$ is defined in Section \ref{Section_Variational_form_of_the_problem}),  from the equality 
\begin{align}\label{2_term_regularization_uL_equality_quadratic_functional_energy_norm_error}
\vertiii{\nabla (u_{h_L}^L-u^L)}^2=2\left(J^L(u_{h_L}^L)-J^L(u^L)\right)
\end{align}
it follows that for all $w\in H_{g-G}^1(\Omega)$
\begin{align}
\vertiii{\nabla (u_{h_L}^L-u^L)}^2\ge 2\left(J^L(u_{h_L}^L)-J^L(w)\right)=:M_{\ominus,L}^2(u_{h_L}^L,w).
\end{align}
For $w$ we always take the last available approximation $u_{h_L}^L$ from the adaptive procedure and compute the lower bound for the error on all previous levels. For convenience, we will denote the approximation $u_{h_L}^L$ on mesh level $i$ by $u_i^L$. Instead $y_L^*=(\epsilon_m-\epsilon)\nabla G+y_0^*$ and $y_0^*$ we write $y_{L,i}^*$ and $y_{0,i}^*$ where $i=0,1,..\bar p$, $i=0$ corresponds to the initial mesh, and $i=\bar p$ corresponds to the last mesh. The results after solving adaptively for $u^L$ are shown in the tables below where $\|v\|_0$ denotes the $L^2(\Omega)$ norm of the function $v$ and $\bar p=7$.
{\small
\begin{center}
\captionof{table}{Example 1} \label{Example1_Table1} 
\vspace{1ex}
%\begin{tabular}{ |p{1.9cm}|p{1.8cm}|p{2.cm}|p{2.1cm}|p{2cm}|p{2cm}|p{2cm}| }
\begin{tabular}{ |c|l|c|c|c|c|c| }
\hline
\multicolumn{7}{|c|}{Example 1: $k_m^2=0,\,k_s^2=10,\,\epsilon_m=2,\,\epsilon_s=80$} \\
\hline
 & & & & & &\\
level $i$ &\#elements &$\|u_i^L\|_0$ &$\vertiii{\nabla u_i^L}$ & $M_{\ominus,L}(u_i^L,u_{\bar p}^L)$ &
$M_{\oplus,L}(u_i^L,y_{0,i}^*)$& $J^L(u_i^L)$\\
 & & & & & &\\
\hline
0	&667 008 		    	&63584.27 	    &15679.01 		&3318.24			&3358.51		&-122925209.65 \\
1	&1 695 251 		&64014.26      &15977.11 	     &1267.71			&1369.07		&-127627031.79  \\
2	&3 803 582 	     &64064.47      &16006.50 		&819.032			&968.374		&-128095169.28 \\
3	&7 238 416    	     &64094.10      &16018.35 		&543.720			&749.503		&-128282760.11 \\
4	&10 268 886       &64109.01      &16022.61 		&401.463		     &653.047		&-128349989.97\\
5	&13 164 899       &64115.19 	   &16024.94 		&295.404			&593.411		&-128386944.42\\
6	&16 124 993       &64119.55      &16026.50 		&194.810			&549.973		&-128411600.81 \\
7	&19 531 518       &64122.38      &16027.69 		&0.00000			&514.037		&-128430576.31\\
\hline
\end{tabular}
\end{center}
}
We can also find guaranteed lower and upper bounds on the relative errors in energy and combined energy norm, as well as practical estimations for these quantities (see \cite{Kraus_Nakov_Repin_PBE1_2018}). The combined energy norm of the pair $(v,q)\in H_0^1(\Omega)\times \left[L^2(\Omega)\right]^d$ is defined by
\begin{align*}
\vertiii{(v,q)}_{\text{CEN}}:=\sqrt{\vertiii{\nabla v}^2+\vertiii{q}_*^2} .
\end{align*}
By ${\text{RE}}_{i,j,k,s}^{\text{L,Up}}$ we denote the guaranteed upper bound for the relative error in energy norm, by ${\text{RCEN}}_{i,j,s}^{\text{L,Up}}$ the guaranteed
upper bound on the relative error in combined energy norm, by ${\text{RE}}_{i,j,k,s}^{\text{L,Low}}$ the guaranteed lower bound on the relative energy norm error, and by
${\text{RCEN}}_{i,j,s}^{\text{L,Low}}$ the guaranteed lower bound for the relative error in combined energy norm where the indices $i,j,k,s$ correspond to the refinement
levels from which approximations for $u^L$ and $p_L^*$ are taken. For any $i,j,k,s\in\{0,1,2,...,\bar p\}$ we have
\begin{align*}
\vertiii{\nabla u_j^L}-M_{\oplus,L}(u_j^L,y_{0,s}^*)\leq\vertiii{\nabla u_j^L-\nabla (u_j^L-u^L)}\leq \vertiii{\nabla u_j^L}+M_{\oplus,L}(u_j^L,y_{0,s}^*).
\end{align*}
Therefore,
\begin{subequations}\label{2_term_uL_Sure_Bounds_For_Relative_Energy_Error}
\begin{eqnarray}
\frac{\vertiii{\nabla(u_i^L-u^L)}}{\vertiii{\nabla u^L}}&\leq&\frac{M_{\oplus,L}(u_i^L,y_{0,k}^*)}{\vertiii{\nabla u_j^L}-M_{\oplus,L}(u_j^L,y_{0,s}^*)}
=:{\text{RE}}_{i,j,k,s}^{\text{L,Up}}, \label{2_term_uL_Sure_Bounds_For_Relative_Energy_Error_1}\\
{\text{RE}}_{i,j,k,s}^{\text{L,Low}}:&= & \frac{M_{\ominus,L}(u_i^L,u_k^L)}{\vertiii{\nabla u_j^L}+M_{\oplus,L}(u_j^L,y_{0,s}^*)}\leq \frac{\vertiii{\nabla(u_i^L-u^L)}}{\vertiii{\nabla u^L}} , \label{2_term_uL_Sure_Bounds_For_Relative_Energy_Error_2}
\end{eqnarray}
\end{subequations}
where \eqref{2_term_uL_Sure_Bounds_For_Relative_Energy_Error_1} is valid if $\vertiii{\nabla u_j^L}-M_{\oplus,L}(u_j^L,y_{0,s}^*)>0$. For any level $i$, the above bounds are expected to be the sharpest when we take $j,k,s=\bar p$. In practice, on each level $i$, the best one can do is to take for ${\text{RE}}_{i,j,k,s}^{\text{L,Up}}$ $j=i,s=i,k=i$. Optionally, once the computations are done, i.e., we have reached level $\bar p$, one can return and recompute slightly sharper upper bounds for each $i=0,1,...,\bar p$ with $j=\bar p, s=\bar p$ and $k=i$. This results in only 2 arithmetic operations per level, provided that $\vertiii{\nabla u_{\bar p}^L}$ and $M_{\oplus,L}(u_i^L,y_{0,i}^*),\,i=0,1,...,\bar p$ are saved. On the other hand, ${\text{RE}}_{i,j,k,s}^{\text{L,Low}}$ is equal to zero if $k=i$ and it is expected that ${\text{RE}}_{i,j,k,s}^{\text{Low}}$ will be negative if $k<i$. Therefore, on each level $i$, to compute the best possible lower bounds for all previous levels $0,1,...,i-1$, we take $k=i, j=i,s=i$.
Further, we have the equality
\begin{align}\label{Prager_Synge_Equality}
\vertiii{\epsilon\nabla v-y^*}_*^2=\vertiii{\nabla(v-u)}^2 + \vertiii{y^*-p^*}_*^2 - 2\int\limits_{\Omega}{(y^*-p^*)\cdot\nabla(v-u) dx},
\end{align}
 which holds for any $v,u\in H^1(\Omega)$ and any $y^*,p^*\in \left[L^2(\Omega)\right]^d$. Also,
 \begin{align}\label{upper_bound_integral_in_Prager_Synge}
  &\mathrel{\phantom{=}}\abs{\int\limits_{\Omega}{(y_L^*-p_L^*)\cdot\nabla(u_{h_L}^L-u^L) dx}}\\
  &=\abs{\int\limits_{\Omega}{\div y_0^*(u_{h_L}^L-u^L)dx}}\leq \frac{C_{F\Omega}\|\div y_0^*\|_{L^2(\Omega)} }{\sqrt{\epsilon_{\min}}}M_{\oplus,L}(u_{h_L}^L,y_0^*),
 \end{align}
and, therefore, we obtain the estimate 
 \begin{align}\label{2_term_uL_Lower_Upper_Bounds_For_CEN}
 \begin{aligned}
 &\mathrel{\phantom{\leq}}\left(M_{\ominus,L}^{\text{CEN}}(u_{h_L}^L,y_0^*)\right)^2:=\vertiii{\epsilon\nabla u_{h_L}^L-y_L^*}_*^2-2\frac{C_{F\Omega}\|\div y_0^*\|_{L^2(\Omega)} }{\sqrt{\epsilon_{\min}}}M_{\oplus,L}(u_{h_L}^L,y_0^*)\\
 &\leq \vertiii{(u_{h_L}^L-u^L, y_L^*-p_L^*) }_{\text{CEN}}^2\\
 &\leq \vertiii{\epsilon\nabla u_{h_L}^L-y_L^*}_*^2+2\frac{C_{F\Omega}\|\div y_0^*\|_{L^2(\Omega)} }{\sqrt{\epsilon_{\min}}}M_{\oplus,L}(u_{h_L}^L,y_0^*)=:\left(M_{\oplus,L}^{\text{CEN}}(u_{h_L}^L,y_0^*)\right)^2.
 \end{aligned}
 \end{align}
Since $y_0^*$ is found by minimization of $M_{\oplus,L}^2(u_{h_L}^L,y_0^*;\alpha)$, in our experiments $\|\div y_0^*\|_{L^2(\Omega)}$ is usually of the order $10^{-5}$ to $10^{-4}$ and the above estimate turns out to be very sharp.
Then for any level $i=0,1,...,\bar p$, we can bound the relative error in the combined energy norm as follows
 \begin{align}\label{2_term_uL_Sure_Bounds_For_Relative_Combined_Energy_Norm_Error}
 \begin{aligned}
 &\frac{\vertiii{\left(u_i^L-u^L,y_{L,i}^*-p_L^*\right)}_{\text{CEN}}}{\vertiii{(u^L,p_L)}_{\text{CEN}}}\leq \frac{M_{\oplus,L}^{\text{CEN}}(u_i^L,y_{0,i}^*)}{\sqrt{2}\left(\vertiii{\nabla u_j^L}-M_{\oplus,L}(u_j^L,y_{0,s}^*)\right)} :={\text{RCEN}}_{i,j,s}^{\text{L,Up}}\\
  &{\text{RCEN}}_{i,j,s}^{\text{L,Low}}:=\frac{M_{\ominus,L}^{\text{CEN}}(u_i^L,y_{0,i}^*)}{\sqrt{2}\left(\vertiii{\nabla u_j^L}+M_{\oplus,L}(u_j^L,y_{0,s}^*)\right)}\leq \frac{\vertiii{\left(u_i^L-u^L,y_{L,i}^*-p_L^*\right)}_{\text{CEN}}}{\vertiii{(u^L,p_L)}_{\text{CEN}}}.
 \end{aligned}
 \end{align}
 For every level $i$, the sharpest estimates ${\text{RCEN}}_{i,j,s}^{\text{L,Up}}$ and ${\text{RCEN}}_{i,j,s}^{\text{L,Low}}$  are obtained when $j=i,s=i$. In the table below, we also present the practical estimation ${\text{P}}_{\text{rel},i}^{\text{L,CEN}}$ for the relative error in combined energy norm given by
 \begin{align}\label{Practical_estimation_for_CEN_Error}
 {\text{P}}_{\text{rel},i}^{\text{L,CEN}}:=\frac{\vertiii{\epsilon\nabla u_i^L-y_{L,i}^*}_*}{\sqrt{2}\vertiii{\nabla u_i^L}},\,\text{for all }i=0,1,...,\bar p.
 \end{align}
 
{\small
\begin{center}
\captionof{table}{Example 1} \label{Example1_Table2} 
\vspace{1ex}
%\begin{tabular}{ |p{1.9cm}|p{1.8cm}|p{2.cm}|p{2.1cm}|p{2cm}|p{2cm}|p{2cm}| }
\begin{tabular}{ |c|l|c|c|c|c|c| }
\hline
\multicolumn{7}{|c|}{Example 1: $k_m^2=0,\,k_s^2=10,\,\epsilon_m=2,\,\epsilon_s=80$} \\
\hline
 & & & & & &\\
level $i$ &\#elements & ${\text{RE}}_{i,\bar p,\bar p,\bar p}^{\text{L,Low}}[\%]$   & ${\text{RE}}_{i,\bar p,i,\bar p}^{\text{L,Up}}[\%]$  & ${\text{RCEN}}_{i,\bar p,\bar p}^{\text{L,Low}}[\%]$ & ${\text{P}}_{\text{rel},i}^{\text{L,CEN}}[\%]$ & ${\text{RCEN}}_{i,\bar p,\bar p}^{\text{L,Up}}[\%]$ \\
 & & & & & &\\
\hline
0	&667 008       	&20.0598 	    &21.6487 		&14.3546			&15.1455		&15.3079 \\
1	&1 695 251 	&7.66370      &8.82493 	     &5.85216			&6.05907		&6.24017  \\
2	&3 803 582 	&4.95130      &6.24207 		&4.13936			&4.27784		&4.41381 \\
3	&7 238 416    	&3.28696      &4.83125 		&3.20376			&3.30850		&3.41621 \\
4	&10 268 886   	&2.42697      &4.20949 		&2.79147		     &2.88196		&2.97656\\
5	&13 164 899   	&1.78581 	   &3.82509 		&2.53655			&2.61840		&2.70474\\
6	&16 124 993    &1.17768      &3.54509		&2.35088			&2.42650		&2.50675 \\
7	&19 531 518    &0.00000      &3.31345 		&2.19726			&2.26777		&2.34296\\
\hline
\end{tabular}
\end{center}
}
 
\subsubsection{Finding $\tilde u^N$}\label{SubSubsection_2_term_Solving_for_uN_2DYES}

Once we have obtained an approximation $u_{h_L}^L$ for $u^L$, we adaptively solve \eqref{2_term_regularization_un_weak_formulation} with $u_{h_L}^L$ instead of $u^L$ in it to find approximations $\tilde u_{h_N}^N$ of $\tilde u^N$. For $u_{h_L}^L$ we take the approximation $u_2^L$ from level 2. In this case, we have $M_{\oplus,L}(u_2^L,y_{0,2}^*)=968.374$, see Table \ref{Example1_Table1},
and ${\text{RE}}_{2,2,2,2}^{\text{L,Up}}=6.43946 \%$.
For $\tilde y_N^*\in H(\div;\Omega)$ with $\div(\tilde y_N^*)=0$ in $\Omega_m$, we use a patchwise equilibrated reconstruction of the numerical flux $\epsilon \nabla \tilde u_{h_N}^N$ based on \cite{Braess_Schoberl_2006}. More precisely, we find $\tilde y_N^*$ in the Raviart-Thomas space $RT_0$ over the same mesh, such that its divergence is equal to the $L^2$ orthogonal projection of $k^2\sinh(\tilde u_{h_N}^N+G+u_{h_L}^L)$ onto the space of piecewise constants. Since the computations on each patch are independent from the computations on the rest of the patches, this reconstruction is easy to implement in parallel. As an error indicator, we use the quantity $\eta_K^N$ \eqref{Explicit_form_of_Functional_error_equality}.
\begin{align}\label{2_term_regularization_error_indicator_NonLinear_Part}
\eta_K^N=\left(\vertiii{\epsilon\nabla \tilde u_{h_N}^N-\tilde y_N^*}_{*(K)}^2+2D_{F,K}(\tilde u_{h_N}^N,-\Lambda^* \tilde y_N^*)\right)^{\frac{1}{2}}.
\end{align}
where $D_{F,K}(\tilde u_{h_N}^N,-\Lambda^* \tilde y_N^*)$ is defined as in \eqref{Explicit_form_DF_v_minus_lambda_yStar} but with integration taking place only on elements $K\in \Omega_s$ and with $G+u_{h_L}^L$ instead of $\tilde u_{h_L}^L$. From \eqref{Explicit_form_of_Functional_error_equality} we have the following upper bounds for the error in energy and combined energy norm
\begin{align*}
&\vertiii{\nabla(\tilde u_{h_N}^N-\tilde u^N)}\leq \sqrt{2}M_{\oplus,N}(\tilde u_{h_N}^N,\tilde y_N^*)\\
&\vertiii{\left(\tilde u_{h_N}^N-\tilde u^N, \tilde y_N^*-\tilde p_N^*\right)}_{\text{CEN}}\leq \sqrt{2}M_{\oplus,N}(\tilde u_{h_N}^N,\tilde y_N^*).
\end{align*}
We will denote by $\tilde u_i^N$ the finite element approximations of $\tilde u^N$ and by $\tilde y_{N,i}^*$ the approximations of the flux $\tilde p_N^*$ at mesh refinement level $i,i=0,1,2,...,\bar p$, where $\bar p=6$. By $J_{h_L}^N: H_0^1(\Omega)\to\mathbb R\cup\{+\infty\}$ we denote the functional defined by
\begin{align}\label{definition_of_J_hL_N}
J_{h_L}^N(v):=\left\{
\begin{aligned}
&\int\limits_{\Omega}{\left[\frac{\epsilon(x)}{2}\abs{\nabla v}^2+k^2\cosh(v+G+u_{h_L}^L)\right]dx},\text{ if }  k^2\cosh(v+G+u_{h_L}^L)\in L^1(\Omega),\\
&+\infty, \text{ if } k^2\cosh(v+G+u_{h_L}^L) \notin L^1(\Omega).
\end{aligned}
\right.
\end{align}
The unique minimizer of $J_{h_L}^N$ over $H_0^1(\Omega)$ is the solution $\tilde u^N$ to the problem \eqref{2_term_regularization_un_weak_formulation} with $u_{h_L}^L$ instead of $u^L$ in it (see \cite{Kraus_Nakov_Repin_PBE1_2018}). The subindex $h_L$ in the notation for the functional $J_{h_L}^N$ corresponds to the mesh refinement level on which the approximation $u_{h_L}^L$ of $u^L$ is computed. Since in this case we take $u_2^L$ as an approximation of $u^L$, we are interested in the values of $J_2^N(\tilde u_i^N)$ on levels $i=0,1,...,6$ in the adaptive solution for $\tilde u^N$ (see Table \ref{Example1_Table3}).
%{\small
%\begin{center}
%\captionof{table}{Example 1} \label{Example1_Table3} 
%\vspace{1ex}
%%\begin{tabular}{ |p{1.9cm}|p{1.8cm}|p{2.cm}|p{2.1cm}|p{2cm}|p{2cm}|p{2cm}| }
%\begin{tabular}{ |c|l|c|c|c|c|c| }
%\hline
%\multicolumn{7}{|c|}{Example 1: $k_m^2=0,\,k_s^2=10,\,\epsilon_m=2,\,\epsilon_s=80$} \\
%\hline
% & & & & & &\\
%level $i$ &\#elts &$\|\tilde u_i^N\|_0$ &$\vertiii{\nabla \tilde u_i^N}$ & $\vertiii{\epsilon\nabla \tilde u_i^N-\tilde y_{N,i}^*}_*^2$ &
%$2M_{\oplus,N}^2(\tilde u_i^N,\tilde y_{N,i}^*)$& $J^N(\tilde u_i^N)$\\
% & & & & & &\\
%\hline
%0	&667 008		&927.667 	    &324.3302 		&24063.116			&37057.026		&259017030.567 \\
%1	&1 315 573 	&928.279      &320.4252 	     &4921.8825			&11792.234		&259013323.940  \\
%2	&5 800 985 	&928.384      &321.6557 		&1551.0627			&3828.2940		&259012091.454 \\
%3	&9 514 417    &928.394      &321.6757 		&1096.4029			&2806.0850		&259011967.304 \\
%4	&13 957 123  &928.399      &321.7411 		&869.45152		     &2236.0886		&259011894.508\\
%5	&18 286 791  &928.401 	    &321.7820 		&717.94712			&1889.1247		&259011849.222\\
%6	&22 883 680  &928.403 	    &321.8072 		&617.65922			&1645.7526		&259011818.931\\
%\hline
%\end{tabular}
%\end{center}
%}

{\small
\begin{center}
\captionof{table}{Example 1} \label{Example1_Table3} 
\vspace{1ex}
\begin{tabular}{ |p{0.98cm}|p{1.68cm}|p{1.1cm}|p{1.25cm}|p{2.7cm}|p{3.cm}|p{2.25cm}| }
%\begin{tabular}{ |c|l|c|c|c|c|c| }
\hline
\multicolumn{7}{|c|}{Example 1: $k_m^2=0,\,k_s^2=10,\,\epsilon_m=2,\,\epsilon_s=80$} \\
\hline
 & & & & & &\\
level $i$ &\#elements &$\|\tilde u_i^N\|_0$ &$\vertiii{\nabla \tilde u_i^N}$ & $\vertiii{\epsilon\nabla \tilde u_i^N-\tilde y_{N,i}^*}_*$ &
$\sqrt{2}M_{\oplus,N}(\tilde u_i^N,\tilde y_{N,i}^*)$& $J_2^N(\tilde u_i^N)$\\
 & & & & & &\\
\hline
0	&667 008	&927.667 	 &324.330 		& \hs155.122			& \hs192.502		&259017030.567 \\
1	&1 315 573 	&928.279      &320.425 	        & \hs 70.1561			&\hs 108.592		&259013323.940  \\
2	&5 800 985 	&928.384      &321.655 		& \hs 39.3835			& \hs 61.8732		&259012091.454 \\
3	&9 514 417    &928.394      &321.675 		& \hs 33.1119			& \hs 52.9724		&259011967.304 \\
4	&13 957 123  &928.399      &321.741 		& \hs 29.4864		        & \hs 47.2872		&259011894.508\\
5	&18 286 791  &928.401 	 &321.782 		& \hs 26.7945			& \hs 43.4640		&259011849.222\\
6	&22 883 680  &928.403 	 &321.807 		& \hs 24.8527			& \hs 40.5678		&259011818.931\\
\hline
\end{tabular}
\end{center}
}
As for the linear part $u^L$ we can define the following guaranteed lower and upper bounds on the relative errors.
 \begin{align}\label{2_term_uN_Sure_Bounds_For_Relative_Energy_Error}
\frac{\vertiii{\nabla(\tilde u_i^N-\tilde u^N)}}{\vertiii{\nabla \tilde u^N}}\leq\frac{\sqrt{2}M_{\oplus,N}(\tilde u_i^N,\tilde y_{N,k}^*)}{\vertiii{\nabla \tilde u_j^N}-\sqrt{2}M_{\oplus,N}(\tilde u_j^N,\tilde y_{N,s}^*)}
=:{\text{RE}}_{i,j,k,s}^{\text{N,Up}}
\end{align}
Using \eqref{strengthened_error_inequality_Th7_2_1} and the estimates
\begin{align*}
&\mathrel{\phantom{\leq}}\vertiii{(\tilde u_i^N,\tilde y_{N,s}^*)}_{\text{CEN}}-\sqrt{2}M_{\oplus,N}\left(\tilde u_i^N,\tilde y_{N,s}^*\right)\\
&\leq \vertiii{(\tilde u^N,\tilde p_N^*)}_{\text{CEN}}\leq \vertiii{(\tilde u_i^N,\tilde y_{N,s}^*)}_{\text{CEN}}+\sqrt{2}M_{\oplus,N}\left(\tilde u_i^N,\tilde y_{N,s}^*\right),\,\text{for all }i,s=0,1,...,\bar p
\end{align*}
we can define the following guaranteed lower and upper bounds for the relative error in combined energy norm:
\begin{align*}%\label{2_term_uN_Sure_Bounds_For_Relative_Combined_Energy_Norm_Error}
\begin{aligned}
&\frac{\vertiii{\left(\tilde u_i^N-\tilde u^N,\tilde y_{N,i}^*-\tilde p_N^*\right)}_{\text{CEN}}}{\vertiii{\left(\tilde u^N,\tilde p_N^*\right)}_{\text{CEN}}}\leq \frac{\sqrt{2}M_{\oplus,N}\left(\tilde u_i^N,\tilde y_{N,i}^*\right)}{\vertiii{\left(\tilde u_j^N,\tilde y_{N,s}^*\right)}_{\text{CEN}}-\sqrt{2}M_{\oplus,N}\left(\tilde u_j^N,\tilde y_{N,s}^*\right)} :={\text{RCEN}}_{i,j,s}^{\text{N,Up}}\\
&{\text{RCEN}}_{i,j,s}^{\text{N,Low}}:=\frac{\frac{1}{\sqrt{2}}\vertiii{\epsilon\nabla \tilde u_i^N-\tilde y_{N,i}^*}}{\vertiii{\left(\tilde u_j^N,\tilde y_{N,s}^*\right)}_{\text{CEN}}+\sqrt{2}M_{\oplus,N}\left(\tilde u_j^N,\tilde y_{N,s}^*\right)}\leq \frac{\vertiii{\left(\tilde u_i^N-\tilde u^N,\tilde y_{N,i}^*-\tilde p_N^*\right)}_{\text{CEN}}}{\vertiii{\left(\tilde u^N,\tilde p_N^*\right)}_{\text{CEN}}}.
\end{aligned}
\end{align*}
The sharpest values for ${\text{RE}}_{i,j,k,s}^{\text{N,Up}}$, ${\text{RCEN}}_{i,j,s}^{\text{N,Up}}$, and ${\text{RCEN}}_{i,j,s}^{\text{N,Low}}$ at each level $i$ are expected to be obtained when $j=\bar p, k=\bar p, s=\bar p$ (assuming that we do not have another better approximation $\tilde y_{N}^*$ for the flux $\tilde p_N^*$). 
%In this case, one has to reavaluate the majorant $M_{\oplus,N}(\tilde u_i^N,\tilde y_{N,\bar p}^*)$ at each level $i$. A much cheaper approach to obtain sharper bounds is to keep $k=i$ and take $j=\bar p,s=\bar p$ since all needed quantities are already computed. 
%In practice, at each level $i$, one uses $j=i, s=i$. Similarly, the sharpest values for ${\text{RCEN}}_{i,j,s}^{\text{N,Low}}$ (assuming we do not have a better approximation $\tilde y_N^*$ to $\tilde p_N^*$) at a given mesh refinement level $i$ are expected when $j=\bar p, s=\bar p$ and in practice one uses $j=i,s=i$.
 The practical estimation $P_{\text{rel},i}^{\text{N,CEN}}$ for the relative error in combined energy norm is given by
 \begin{align}\label{2term_uN_Practical_estimation_for_CEN_Error}
 {\text{P}}_{\text{rel},i}^{\text{N,CEN}}:=\frac{\vertiii{\epsilon\nabla \tilde u_i^N-\tilde y_{N,i}^*}_*}{\sqrt{2}\vertiii{\nabla \tilde u_i^N}},\,\text{for all }i=0,1,...,\bar p.
 \end{align}
We also introduce a practical upper bound
\begin{align}\label{2term_uN_Practical_estimation_of_Upper_Bound_for_Energy_Error}
{\text{PRE}}_{i,j}^{\text{N,Up}}:=\frac{\vertiii{\epsilon\nabla \tilde u_i^N-\tilde y_{N,j}^*}_*}{\vertiii{\nabla \tilde u_i^N}},\,\text{for all }i,j=0,1,...,\bar p ,
\end{align}
for the relative error in energy norm which is based on the relation
$$\vertiii{\nabla (\tilde u_i^N-\tilde u^N)}\leq \vertiii{\left(\tilde u_i^N-\tilde u^N,\tilde y_{N,i}^*-\tilde p_N^*\right)}\approx \vertiii{\epsilon\nabla \tilde u_i^N-\tilde y_{N,i}^*}_*$$
and is useful when it is suspected that the guaranteed upper bound for the relative error overestimates the real error.
 
The above introduced bounds on the relative errors are presented in Table~\ref{Example1_Table4} and Table~\ref{Example1_Table5}.
{\small
\begin{center}
\captionof{table}{Example 1} \label{Example1_Table4} 
\vspace{1ex}
\begin{tabular}{ |p{0.98cm}|p{1.68cm}|p{2cm}|p{1.74cm}|p{2.5cm}|p{1.76cm}|p{2.3cm}| }
%\begin{tabular}{ |c|l|c|c|c|c|c| }
\hline
\multicolumn{7}{|c|}{Example 1: $k_m^2=0,\,k_s^2=10,\,\epsilon_m=2,\,\epsilon_s=80$} \\
\hline
 & & & & & &\\
level $i$ &\#elements &${\text{PRE}}_{i,i}^{\text{N,Up}}[\%]$ & ${\text{RE}}_{i,i,i,i}^{\text{N,Up}}[\%]$ & ${\text{RCEN}}_{i,i,i}^{\text{N,Low}}[\%]$ & ${\text{P}}_{\text{rel},i}^{\text{N,CEN}}[\%]$& ${\text{RCEN}}_{i,i,i}^{\text{N,Up}}[\%]$\\
 & & & & & &\\
\hline
0	&667 008		&47.828 	   &146.02 	&16.227			&33.819		&66.164 \\
1	&1 315 573 	&21.894      &51.263	     &8.7551			&15.481		&31.076  \\
2	&5 800 985 	&12.243      &23.817 		&5.3764			&8.6578		&15.695 \\
3	&9 514 417    &10.293      &19.714		&4.6026			&7.2786		&13.152 \\
4	&13 957 123  &9.1646      &17.229 		&4.1455		     & 6.4803		&11.579\\
5	&18 286 791  &8.3269 	  &15.616 		&3.7964			&5.8880		&10.546\\
6	&22 883 680  &7.7228	  &14.424		&3.5422			&5.4608		&9.7758\\
\hline
\end{tabular}
\end{center}
}

{\small
\begin{center}
\captionof{table}{Example 1} \label{Example1_Table5} 
\vspace{1ex}
\begin{tabular}{ |p{0.98cm}|p{1.68cm}|p{1.2cm}|p{1.74cm}|p{1.66cm}|p{2.7cm}|p{3cm}| }
%\begin{tabular}{ |c|l|c|c|c|c|c| }
\hline
\multicolumn{7}{|c|}{Example 1: $k_m^2=0,\,k_s^2=10,\,\epsilon_m=2,\,\epsilon_s=80$} \\
\hline
 & & & & & &\\
level $i$ &\#elements & ${\text{RE}}_{i,\bar p,\bar p,\bar p}^{\text{N,Up}}$ & ${\text{RCEN}}_{i,\bar p,\bar p}^{\text{N,Low}}$ & ${\text{RCEN}}_{i,\bar p,\bar p}^{\text{N,Up}}$ & $\vertiii{\epsilon\nabla \tilde u_i^N-\tilde y_{N,\bar p}^*}_*$  & $\sqrt{2}M_{\oplus,N}(\tilde u_i^N,\tilde y_{N,\bar p}^*)$ \\
 & &$[\%]$ &$[\%]$ &$[\%]$ & &\\
\hline
0	&667 008		   &40.974 	&22.129			&46.434	  	&89.292    &115.23\\
1	&1 315 573 	   &24.311	     &10.008			&26.194		&51.983	&68.373\\
2	&5 800 985 	   &16.737 	&5.6183			&14.924		&32.512	&47.071\\
3	&9 514 417       &15.855		&4.7236			&12.777		&29.023	&44.591\\
4	&13 957 123     &15.269 		&4.2064		     &11.406		&27.033	&42.943\\
5	&18 286 791     &14.841 		&3.8224			&10.484		&25.734	&41.741\\
6	&22 883 680  	  &14.424		&3.5422			&9.7758		&24.852	&40.567\\
\hline
\end{tabular}
\end{center}
}
Finally, according to \eqref{the_estimate_on_u_minus_tildeu_2term_regularization} the overall error in the regular component $u$ will be 
\begin{align*}
&\vertiii{\nabla(u-\tilde u_h)}\leq 2M_{\oplus,L}(u_{h_L}^L,y_0^*)+\sqrt{2}M_{\oplus,N}(\tilde u_{h_N}^N,\tilde y_N^*)=2M_{\oplus,L}(u_2^L,y_{0,2}^*)+\sqrt{2}M_{\oplus,N}(\tilde u_6^N,\tilde y_{N,6}^*)\\
&=2\times 968.37+40.57 =1977.31
\end{align*}
For comparisson, the energy norm of the approximate regular component $\tilde u_h=u_2^L+\tilde u_6^N$ is $\vertiii{\nabla \tilde u_h}=16276.2$. This means that the relative error in energy norm is no more than approximately $1977.31/16276.2=12.15 \%$.

\subsection{Example 1 (Alexa 488 and Alexa 594) recomputed with $u_4^L$} 
Here, we recompute an approximation $\tilde u_h$ of the regular component $u$ from Example 1. This time we take $u_4^L$ as an approximation of $u^L$ and solve with it for $\tilde u^N$.
For  $u_4^L$, we have $M_{\oplus,L}(u_4^L,y_{0,4}^*)=653.047$ and ${\text{RE}}_{4,4,4,4}^{\text{L,Up}}=4.2489 \%$. The final level is $\bar p=3$.
{\small
\begin{center}
\captionof{table}{Example 1} \label{Example1_Table6} 
\vspace{1ex}
\begin{tabular}{ |p{0.98cm}|p{1.68cm}|p{1.1cm}|p{1.25cm}|p{2.7cm}|p{3.cm}|p{2.25cm}| }
%\begin{tabular}{ |c|l|c|c|c|c|c| }
\hline
\multicolumn{7}{|c|}{Example 1: $k_m^2=0,\,k_s^2=10,\,\epsilon_m=2,\,\epsilon_s=80$} \\
\hline
 & & & & & &\\
level $i$ &\#elements &$\|\tilde u_i^N\|_0$ &$\vertiii{\nabla \tilde u_i^N}$ & $\vertiii{\epsilon\nabla \tilde u_i^N-\tilde y_{N,i}^*}_*$ &
$\sqrt{2}M_{\oplus,N}(\tilde u_i^N,\tilde y_{N,i}^*)$& $J_4^N(\tilde u_i^N)$\\
 & & & & & &\\
\hline
0	&667 008	&880.446 	    &314.647 		&128.067			&142.404		&259007595.769 \\
1	&1 389 691 	&880.942      &309.468 	     &39.1633			&60.4019		&259004495.901  \\
2	&5 706 468 	&880.989      &309.166 		&24.0717			&40.2461		&259004177.333 \\
3	&8 606 657    &880.992      &309.138 		&20.5852			&35.7841		&259004134.610 \\
\hline
\end{tabular}
\end{center}
}

{\small
\begin{center}
\captionof{table}{Example 1} \label{Example1_Table7} 
\vspace{1ex}
\begin{tabular}{ |p{0.98cm}|p{1.68cm}|p{2cm}|p{1.74cm}|p{2.5cm}|p{1.76cm}|p{2.3cm}| }
%\begin{tabular}{ |c|l|c|c|c|c|c| }
\hline
\multicolumn{7}{|c|}{Example 1: $k_m^2=0,\,k_s^2=10,\,\epsilon_m=2,\,\epsilon_s=80$} \\
\hline
 & & & & & &\\
level $i$ &\#elements &${\text{PRE}}_{i,i}^{\text{N,Up}}[\%]$ & ${\text{RE}}_{i,i,i,i}^{\text{N,Up}}[\%]$ & ${\text{RCEN}}_{i,i,i}^{\text{N,Low}}[\%]$ & ${\text{P}}_{\text{rel},i}^{\text{N,CEN}}[\%]$& ${\text{RCEN}}_{i,i,i}^{\text{N,Up}}[\%]$\\
 & & & & & &\\
\hline
0	&667 008		&40.702      &82.676 		&14.972			&28.780		&44.499 \\
1	&1 389 691 	&12.655      &24.251       &5.5423			&8.9484		&15.943  \\
2	&5 706 468 	&7.7860      &14.965 		&3.5610			&5.5055		&10.125 \\
3	&8 606 657    &6.6589      &13.090		&3.0751			&4.7085		&8.9065 \\
\hline
\end{tabular}
\end{center}
}

{\small
\begin{center}
\captionof{table}{Example 1} \label{Example1_Table8} 
\vspace{1ex}
\begin{tabular}{ |p{0.98cm}|p{1.68cm}|p{1.2cm}|p{1.74cm}|p{1.66cm}|p{2.7cm}|p{3cm}| }
%\begin{tabular}{ |c|l|c|c|c|c|c| }
\hline
\multicolumn{7}{|c|}{Example 1: $k_m^2=0,\,k_s^2=10,\,\epsilon_m=2,\,\epsilon_s=80$} \\
\hline
 & & & & & &\\
level $i$ &\#elements & ${\text{RE}}_{i,\bar p,\bar p,\bar p}^{\text{N,Up}}$ & ${\text{RCEN}}_{i,\bar p,\bar p}^{\text{N,Low}}$ & ${\text{RCEN}}_{i,\bar p,\bar p}^{\text{N,Up}}$ & $\vertiii{\epsilon\nabla \tilde u_i^N-\tilde y_{N,\bar p}^*}_*$  & $\sqrt{2}M_{\oplus,N}(\tilde u_i^N,\tilde y_{N,\bar p}^*)$ \\
 & &$[\%]$ &$[\%]$ &$[\%]$ & &\\
\hline
0	&667 008		   &33.507 	&19.131			&35.443	  	&73.624    &91.593\\
1	&1 389 691 	   &16.777      &5.8504			&15.033		&31.708	&45.861\\
2	&5 706 468 	   &13.637 	&3.5959			&10.017		&22.183	&37.279\\
3	&8 606 657       &13.090		&3.0751			&8.9065		&20.585	&35.784\\
\hline
\end{tabular}
\end{center}
}
In Table \ref{Example1_Table8}, it can be seen that the convergence of $\tilde u_{h_N}^N$ to $\tilde u^N$ is faster compared to the case when we used a worse approximation for $u^L$.
Finally, according to \eqref{the_estimate_on_u_minus_tildeu_2term_regularization} the overall error in the regular component $u$ can be estimated by
\begin{eqnarray*}
\vertiii{\nabla(u-\tilde u_h)}&\leq& 2M_{\oplus,L}(u_{h_L}^L,y_0^*)+\sqrt{2}M_{\oplus,N}(\tilde u_{h_N}^N,\tilde y_N^*) \\
&=&2M_{\oplus,L}(u_4^L,y_{0,4}^*)+\sqrt{2}M_{\oplus,N}(\tilde u_3^N,\tilde y_{N,3}^*)\\
&=&2\times 653.047+35.784 =1341.878
\end{eqnarray*}
in this example.
For comparisson, the energy norm of the approximate regular component $\tilde u_h=u_4^L+\tilde u_3^N$ is $\vertiii{\nabla \tilde u_h}=16298.534$. This means that the relative error in energy norm is no more than approximately $1341.878/16298.534=8.23 \%$.

\subsection{Example 2: Insulin protein (PDB ID: 1RWE)}
For the second application, we consider an insulin molecule. The system was prepared by using the crystal structure of the insulin from Protein Data Bank (ID code 1RWE). The CHARMM-GUI web server was employed to add hydrogens to the system. The total number of atoms (with the added hydrogens) is 1590. The charges for calculations were taken from the psf file, created by the  CHARMM-GUI.
In this example, the parameters are $\epsilon_m=20$, $\epsilon_s=80$, $k_m^2=0\,\AA^{-2}$, $k_s^2=10\, \AA^{-2}$ which corresponds to ionic strength $I_s\approx 1.178$ molar. 

\subsubsection{Finding $u^L$}
Here, we use the same notation for the lower and upper bounds for the errors as in Section \ref{SubSubsection_2_term_Solving_for_uL_2DYES} and the last level of refinement on which we have computed an approximation of $u^L$ is $\bar p=7$.
{\small
\begin{center}
\captionof{table}{Example 2} \label{Example2_Table1} 
\vspace{1ex}
%\begin{tabular}{ |p{1.9cm}|p{1.8cm}|p{2.cm}|p{2.1cm}|p{2cm}|p{2cm}|p{2cm}| }
\begin{tabular}{ |c|l|c|c|c|c|c| }
\hline
\multicolumn{7}{|c|}{Example 2: $k_m^2=0,\,k_s^2=10,\,\epsilon_m=20,\,\epsilon_s=80$} \\
\hline
 & & & & & &\\
level $i$ &\#elements &$\|u_i^L\|_0$ &$\vertiii{\nabla u_i^L}$ & $M_{\ominus,L}(u_i^L,u_{\bar p}^L)$ &
$M_{\oplus,L}(u_i^L,y_{0,i}^*)$& $J^L(u_i^L)$\\
 & & & & & &\\
\hline
0	&477 557	&4418.09 	&2096.18 		&1459.88			&1538.92		&-2196016.775 \\
1	&1 097 597 	&4442.42      &2416.23 	       &828.922			&929.997		&-2918092.174  \\
2	&1 864 757 	&4443.71      &2497.76 		&535.229			&641.975		&-3118413.187\\
3	&2 759 585    &4445.16      &2530.36 		&350.049			&474.065		&-3200381.471 \\
4	&3 786 443    &4446.92      &2542.89 		&242.798		       &389.654		&-3232173.083\\
5	&4 980 436    &4449.12 	 &2548.63 		&172.528			&344.466		&-3246765.623\\
6	&6 351 963    &4450.91 	 &2552.02 		&111.874			&314.884		&-3255390.680\\
7	&7 940 767    &4452.58 	 &2554.47 		&0.00000			&291.514		&-3261648.671\\
\hline
\end{tabular}
\end{center}
}

 {\small
\begin{center}
\captionof{table}{Example 2} \label{Example2_Table2} 
\vspace{1ex}
%\begin{tabular}{ |p{1.9cm}|p{1.8cm}|p{2.cm}|p{2.1cm}|p{2cm}|p{2cm}|p{2cm}| }
\begin{tabular}{ |c|l|c|c|c|c|c| }
\hline
\multicolumn{7}{|c|}{Example 2: $k_m^2=0,\,k_s^2=10,\,\epsilon_m=20,\,\epsilon_s=80$} \\
\hline
 & & & & & &\\
level $i$ &\#elements & ${\text{RE}}_{i,\bar p,\bar p,\bar p}^{\text{L,Low}}[\%]$   & ${\text{RE}}_{i,\bar p,i,\bar p}^{\text{L,Up}}[\%]$  & ${\text{RCEN}}_{i,\bar p,\bar p}^{\text{L,Low}}[\%]$ & ${\text{P}}_{\text{rel},i}^{\text{L,CEN}}[\%]$ & ${\text{RCEN}}_{i,\bar p,\bar p}^{\text{L,Up}}[\%]$ \\
 & & & & & &\\
\hline
0	&477 557	&51.29 	    	&68.00		&38.23			&51.91			&48.08 \\
1	&1 097 597 	&29.12     	&41.09 	        &23.10			&27.21			&29.05 \\
2	&1 864 757  	&18.80      	&28.36 		&15.94			&18.17			&20.05\\
3	&2 759 585    &12.29     	&20.94 		&11.77			&13.24			&14.81 \\
4	&3 786 443    &8.531        	&17.21 		&9.681			&10.83			&12.17\\
5	&4 980 436    &6.062 	    	&15.22 		&8.558			&9.557			&10.76\\
6	&6 351 963    &3.930 		&13.91 		&7.823			&8.724			&9.839\\
7	&7 940 767    &0.000	 	&12.88 		&7.242			&8.069			&9.108\\
\hline
\end{tabular}
\end{center}
}
Note that in this example, mesh refinements after level $4$, decrease the error slower than up to level $4$.
This is because the error is already equilibrated on the computational domain $\Omega$ and we should use
a more aggressive refinement strategy.

\subsubsection{Finding $\tilde u^N$}

%We reuse the notation introduced in Section \ref{SubSubsection_2_term_Solving_for_uN_2DYES}.
As an approximation $u_{h_L}^L$ of $u^L$, the notation is the same as in Section~\ref{SubSubsection_2_term_Solving_for_uN_2DYES},
we take the finite element  function $u_7^L$ from level 7. In this case, we have $M_{\oplus,L}(u_7^L,y_{0,7}^*)=291.514$
(see Table \ref{Example2_Table1}), ${\text{RE}}_{7,7,7,7}^{\text{L,Up}}=12.88 \%$, and the last refinement level for $\tilde u^N$
is $\bar p=6$.

{\small
\begin{center}
\captionof{table}{Example 2} \label{Example2_Table3} 
\vspace{1ex}
\begin{tabular}{ |p{0.98cm}|p{1.68cm}|p{1cm}|p{1.35cm}|p{2.7cm}|p{3.cm}|p{2.25cm}| }
%\begin{tabular}{ |c|l|c|c|c|c|c| }
\hline
\multicolumn{7}{|c|}{Example 2: $k_m^2=0,\,k_s^2=10,\,\epsilon_m=20,\,\epsilon_s=80$} \\
\hline
 & & & & & &\\
level $i$ &\#elements &$\|\tilde u_i^N\|_0$ &$\vertiii{\nabla \tilde u_i^N}$ & $\vertiii{\epsilon\nabla \tilde u_i^N-\tilde y_{N,i}^*}_*$ &
$\sqrt{2}M_{\oplus,N}(\tilde u_i^N,\tilde y_{N,i}^*)$& $J_7^N(\tilde u_i^N)$\\
 & & & & & &\\
\hline
0	&477 557	&685.05		  &704.41	 	&9139.8			&9256.1		&193257098.873 \\
1	&737 401 	&687.24      	  &601.23 	       &720.81			&854.60		&193070043.803 \\
2	&1 217 490 	&688.25          &558.35 		&534.08			&602.82		&193038582.359 \\
3	&1 739 923   &688.44         &514.49 		&285.75			&315.90		&193012502.711\\
4	&2 410 640  	&688.49      	  &501.62 		&161.35		       &180.60		&193005499.138\\
5	&3 266 704  	&688.50 	  &496.92 		&88.515			&103.06		&193002971.766\\
6	&4 379 634  	&688.51 	  &495.89 		&69.203			&80.683		&193002341.283\\
\hline
\end{tabular}
\end{center}
}

{\small
\begin{center}
\captionof{table}{Example 2} \label{Example2_Table4} 
\vspace{1ex}
\begin{tabular}{ |p{0.98cm}|p{1.68cm}|p{2cm}|p{1.74cm}|p{2.5cm}|p{1.76cm}|p{2.3cm}| }
%\begin{tabular}{ |c|l|c|c|c|c|c| }
\hline
\multicolumn{7}{|c|}{Example 2: $k_m^2=0,\,k_s^2=10,\,\epsilon_m=20,\,\epsilon_s=80$} \\
\hline
 & & & & & &\\
level $i$ &\#elements &${\text{PRE}}_{i,i}^{\text{N,Up}}[\%]$ & ${\text{RE}}_{i,\bar p,i,\bar p}^{\text{N,Up}}[\%]$ & ${\text{RCEN}}_{i,i,i}^{\text{N,Low}}[\%]$ & ${\text{P}}_{\text{rel},i}^{\text{N,CEN}}[\%]$& ${\text{RCEN}}_{i,i,i}^{\text{N,Up}}[\%]$\\
 & & & & & &\\
\hline
0	&477 557	&1297	       &2229	 	&34.98		&917.4		&-\\
1	&737 401 	&119.8          &205.8 	       	&26.54		&84.77		&405.0\\
2	&1 217 490 	&95.65          &145.1 		&24.67		&67.63		&185.4 \\
3	&1 739 923   &55.54          &76.08 		&18.64		&39.27		&69.90\\
4	&2 410 640  	&32.16          &43.49 		&12.62	       &22.74		&33.29\\
5	&3 266 704  	&17.81 	  	&24.82 		&7.730		&12.59		&17.07\\
6	&4 379 634  	&13.95 	  	&19.43 		&6.236		&9.867		&12.94\\
\hline
\end{tabular}
\end{center}
}

%{\small
%\begin{center}
%\captionof{table}{Example 2} \label{Example2_Table5} 
%\vspace{1ex}
%\begin{tabular}{ |p{0.98cm}|p{1.68cm}|p{1.2cm}|p{1.74cm}|p{1.66cm}|p{2.7cm}|p{3cm}| }
%%\begin{tabular}{ |c|l|c|c|c|c|c| }
%\hline
%\multicolumn{7}{|c|}{Example 1: $k_m^2=0,\,k_s^2=10,\,\epsilon_m=20,\,\epsilon_s=80$} \\
%\hline
% & & & & & &\\
%level $i$ &\#elements & ${\text{RE}}_{i,\bar p,\bar p,\bar p}^{\text{N,Up}}$ & ${\text{RCEN}}_{i,\bar p,\bar p}^{\text{N,Low}}$ & ${\text{RCEN}}_{i,\bar p,\bar p}^{\text{N,Up}}$ & $\vertiii{\epsilon\nabla \tilde u_i^N-\tilde y_{N,\bar p}^*}_*$  & $\sqrt{2}M_{\oplus,N}(\tilde u_i^N,\tilde y_{N,\bar p}^*)$ \\
% & &$[\%]$ &$[\%]$ &$[\%]$ & &\\
%\hline
%0	&477 557	&1297	       &2229	 	&34.98		&917.4		&-\\
%1	&737 401 	&119.8          &205.8 	       	&26.54		&84.77		&405.0\\
%2	&1 217 490 	&95.65          &145.1 		&24.67		&67.63		&185.4 \\
%3	&1 739 923   &55.54          &76.08 		&18.64		&39.27		&69.90\\
%4	&2 410 640  	&32.16          &43.49 		&12.62	       &22.74		&33.29\\
%5	&3 266 704  	&17.81 	  	&24.82 		&7.730		&12.59		&17.07\\
%6	&4 379 634  	&13.95 	  	&19.43 		&6.236		&9.867		&12.94\\
%\hline
%\end{tabular}
%\end{center}
%}
%According to \ref{the_estimate_on_u_minus_tildeu_2term_regularization} the overall error in the regular component $u$ will be 
%\begin{align*}
%&\vertiii{\nabla(u-\tilde u)}\leq 2M_{\oplus,L}(u_{h_L}^L,y_0^*)+\sqrt{2}M_{\oplus,N}(\tilde u_{h_N}^N,\tilde y_N^*)=2M_{\oplus,L}(u_3^L,y_{0,3}^*)+\sqrt{2}M_{\oplus,N}(\tilde u_5^N,\tilde y_{N,5}^*)\\
%&=2\times 291.514+80.683 =663.711
%\end{align*}
{\small
\begin{center}
\captionof{table}{Example 2} \label{Example2_Table5} 
\vspace{1ex}
\begin{tabular}{ |p{0.98cm}|p{1.68cm}|p{1.2cm}|p{1.74cm}|p{1.66cm}|p{2.7cm}|p{3cm}| }
%\begin{tabular}{ |c|l|c|c|c|c|c| }
\hline
\multicolumn{7}{|c|}{Example 2: $k_m^2=0,\,k_s^2=10,\,\epsilon_m=20,\,\epsilon_s=80$} \\
\hline
& & & & & &\\
level $i$ &\#elements & ${\text{RE}}_{i,\bar p,\bar p,\bar p}^{\text{N,Up}}$ & ${\text{RCEN}}_{i,\bar p,\bar p}^{\text{N,Low}}$ & ${\text{RCEN}}_{i,\bar p,\bar p}^{\text{N,Up}}$ & $\vertiii{\epsilon\nabla \tilde u_i^N-\tilde y_{N,\bar p}^*}_*$  & $\sqrt{2}M_{\oplus,N}(\tilde u_i^N,\tilde y_{N,\bar p}^*)$ \\
& &$[\%]$ &$[\%]$ &$[\%]$ & &\\
\hline
0	&477 557	&184.9	       &823.6	 	&1485		&583.82		&768.12\\
1	&737 401 	&95.06          &64.95 	       	&137.1		&356.29		&394.72\\

2	&1 217 490 	&67.72          &48.12 		&96.71		&248.68		&281.18 \\
3	&1 739 923   &39.23          &25.75 		&50.68		&141.60		&162.90\\
4	&2 410 640  	&26.73          &14.54 		&28.97	       &97.204		&110.99\\
5	&3 266 704  	&20.92 	  	&7.976 		&16.53		&76.236		&86.878\\
6	&4 379 634  	&19.43 	  	&6.236 		&12.94		&69.203		&80.683\\
\hline
\end{tabular}
\end{center}
}
According to~\eqref{the_estimate_on_u_minus_tildeu_2term_regularization} the overall error in the regular component $u$ can be estimated as follows:
\begin{eqnarray*}
\vertiii{\nabla(u-\tilde u_h)} &\leq& 2M_{\oplus,L}(u_{h_L}^L,y_0^*)+\sqrt{2}M_{\oplus,N}(\tilde u_{h_N}^N,\tilde y_N^*)\\
&=&2M_{\oplus,L}(u_7^L,y_{0,7}^*)+\sqrt{2}M_{\oplus,N}(\tilde u_6^N,\tilde y_{N,6}^*) =2\times 291.514+80.683 =663.711,
\end{eqnarray*}
see also Table~\ref{Example2_Table1}.

For comparisson, the energy norm of the approximate regular component $\tilde u_h=u_7^L+\tilde u_6^N$ is $\vertiii{\nabla \tilde u_h}=2940.55$. This means that the relative error is no more than
approximately $663.711/2940.55=22.57 \%$. We should note that this estimate is rather conservative. Also, the initial surface mesh that is used in this experiment could be further improved
which will also positively influence the quality of the finite element approximations.
\begin{figure}[H]
    \centering
    \begin{minipage}{1\textwidth}
        \centering
        \captionsetup{width=0.8\linewidth}
        \includegraphics[width=0.8\linewidth]{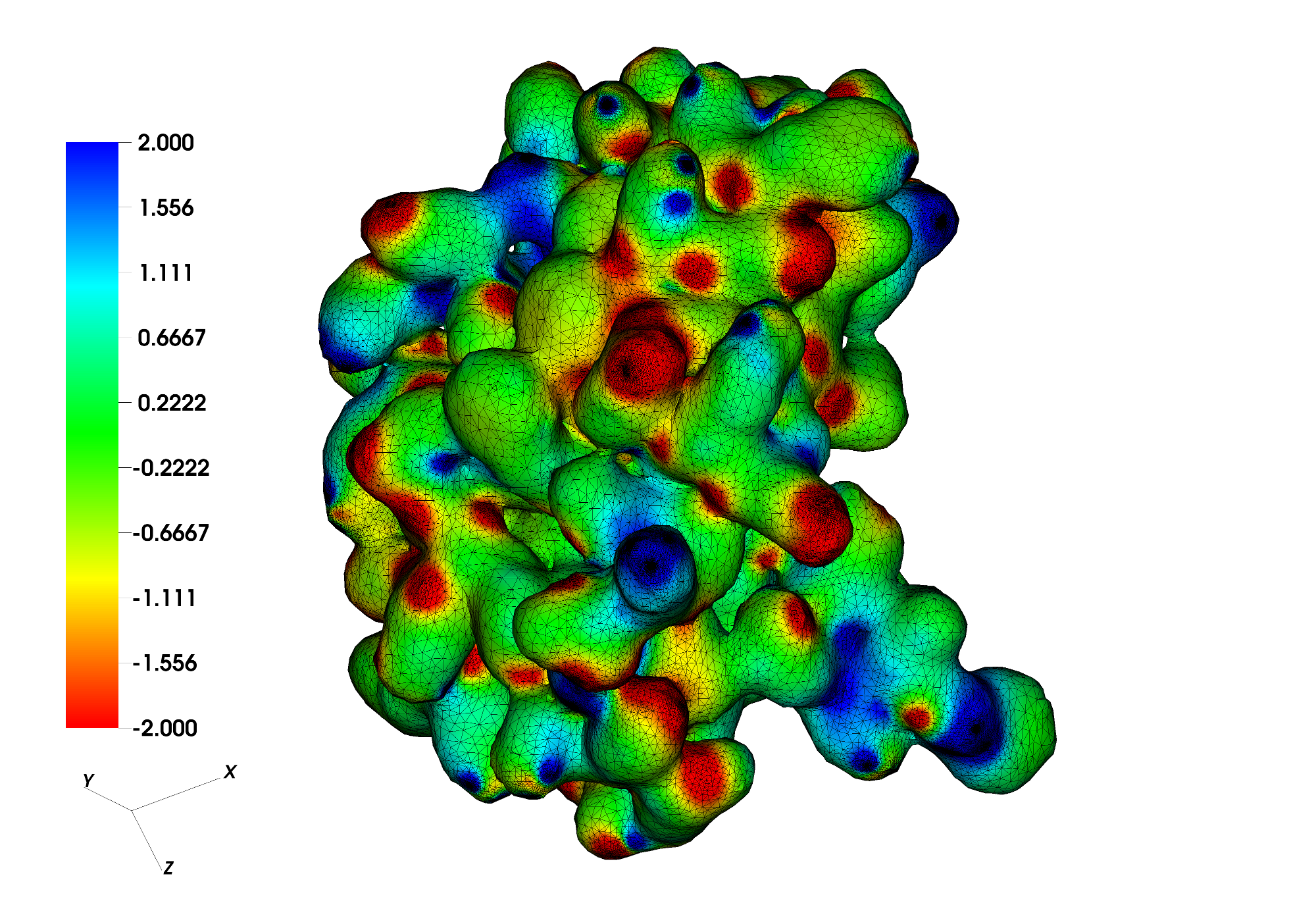}
        \caption{{\small Full potential surface map of the insulin protein (PDB ID: 1RWE) in units $k_BT/e_c$. Blue color indicates a positive potential (values$>2 k_BT/e_c$) and red color indicates negative potential (values $<-2 K_BT/e_c$).}}
        \label{Insulin_Potential_Surface_Map}
    \end{minipage}%
    \end{figure}

\begin{figure}[H]
    \centering
    \begin{minipage}{1\textwidth}
        \centering
        \captionsetup{width=0.8\linewidth}
      \includegraphics[width=0.8\linewidth]{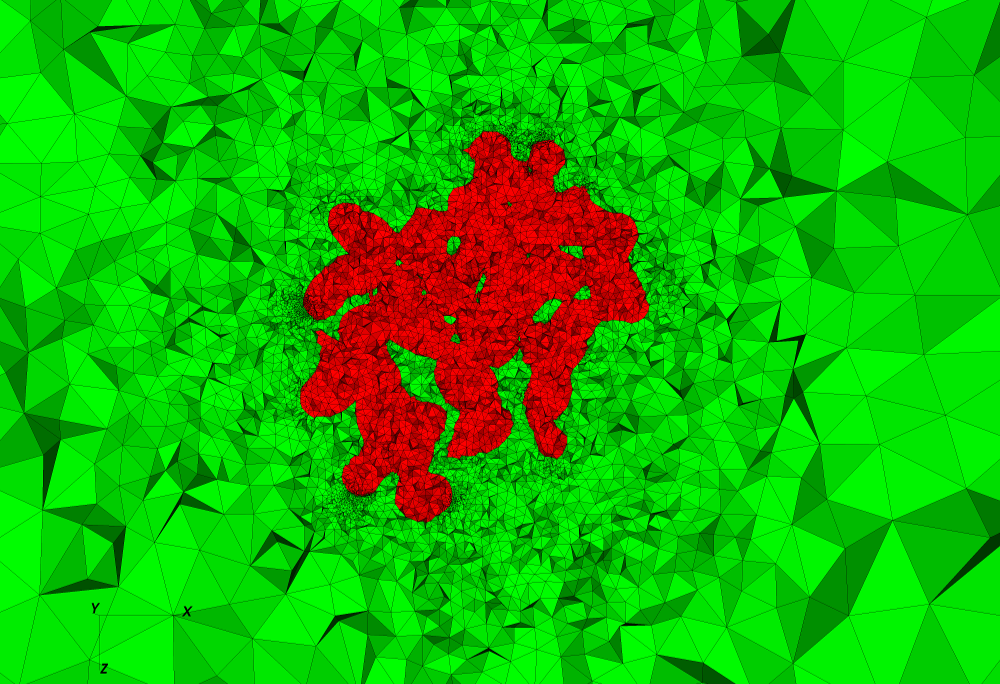}
        \caption{{\small Cross section of the mesh with the plane $y=15\, \AA$ at level 2 in the mesh refinement procedure for finding the component $\tilde u^N$ in Example 2. The molecule region $\Omega_m$ is marked red.}}
        \label{Insulin_Regions_Level_2}
    \end{minipage}%
    \end{figure}
    
    \begin{figure}[H]
    \centering
    \begin{minipage}{1\textwidth}
        \centering
        \captionsetup{width=0.8\linewidth}
      \includegraphics[width=0.8\linewidth]{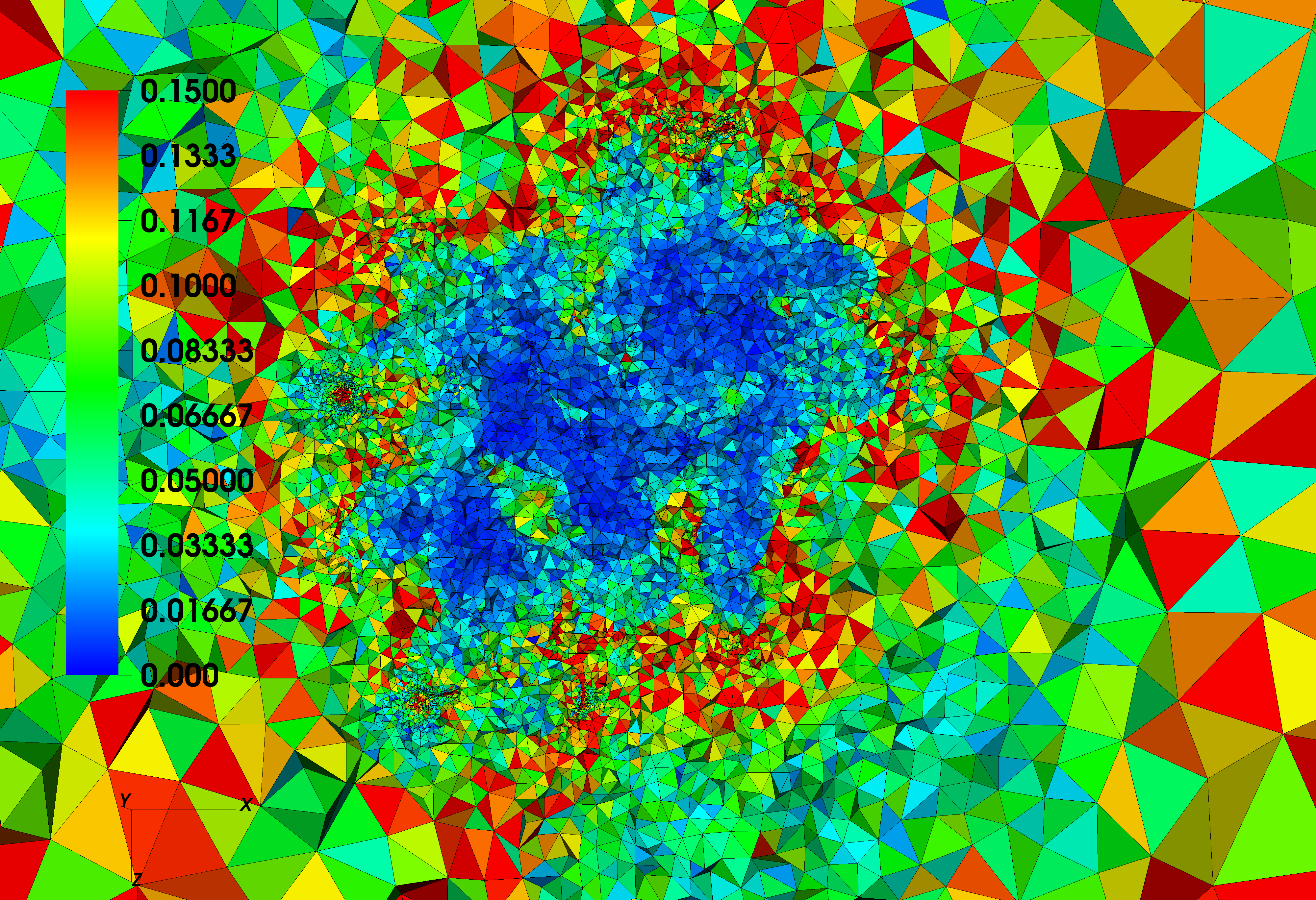}
        \caption{{\small Cross section of the mesh with the plane $y=15\, \AA$ at level 2 in the mesh refinement procedure for finding the component $\tilde u^N$ in Example 2. Error indicator as a piecewise constant function.}}
        \label{Insulin_Func_Error_Indicator_Level_2}
    \end{minipage}%
    \end{figure}
    
\subsection{Example 3: (Alexa 488 and Alexa 594) 3-term regularization without additional splitting in $u^L+u^N$}
In this example, we solve the PBE for the system consisting of the two chromophores Alexa 488 and Alexa 594 but this time we utilize the 3-term regularization scheme without further splitting the regular component $u$ in $u^L+u^N$. The parameters are the same as in the first example, i.e., $\epsilon_m=2$, $\epsilon_s=80$, $k_m^2=0\,\AA^{-2}$, $k_s^2=10\, \AA^{-2}$ which corresponds to ionic strength $I_s\approx 1.178$ molar. 

\subsubsection{Finding the harmonic component $u^H$}

According to the 3-term regularization, we first have to obtain a conforming approximation $\tilde u^H$ of $u^H$ by solving problem \eqref{PBE_special_form_harmonic}. We solve this problem on a sequence of adapted meshes using the error majornat \eqref{A_posteriori_estimate_primal_energy_norm_harmonic_component} and the derived from it error indicator. To reconstruct the numerical flux $\nabla \tilde u^H$ and obtain $T(\nabla \tilde u^H)$ we can either minimize the majorant \eqref{A_posteriori_estimate_primal_energy_norm_harmonic_component} over a subspace of $H(\div;\Omega_m)$, like $RT_0$, or apply some patchwise flux reconstruction technique. We notice that minimization of the majorant over $RT_0$ defined on the same mesh and applying the patchwise equilibrated flux reconstruction from \cite{Braess_Schoberl_2006} yields practically the same results. We have computed $\tilde u^H$ on a final mesh with $20\,615\,534$ tetrahedrons. The corresponding value for the majorant in \eqref{A_posteriori_estimate_dual_energy_norm_harmonic_component} is $M_{\oplus,H}\left(\tilde u^H,T(\nabla \tilde u^H)\right)=43.085$, where $T(\nabla \tilde u^H)$ is obtained by the flux reconstruction in \cite{Braess_Schoberl_2006} and thus $\div\left(T(\nabla \tilde u^H)\right)=0$. Since in this case $\|\nabla\left(\tilde u^H-u^H\right)\|_{L^2(\Omega_m)}\leq M_{\oplus,H}\left(\tilde u^H,T(\nabla \tilde u^H)\right)=43.085$, we obtain a guaranteed upper bound on the relative error in energy norm
\[
\frac{\|\nabla\left(\tilde u^H-u^H\right)\|_{L^2(\Omega_m)}}{\|\nabla u^H\|_{L^2(\Omega_m)}}\leq 3.16\, \%.
\]
%\begin{figure}[!ht]
%    \centering
%    \begin{minipage}{1\textwidth}
%        \centering
%        \captionsetup{width=0.8\linewidth}
%      \includegraphics[width=0.8\linewidth]{run01_2_Poisson_Harmonic_3_Term_2DYES_Ground_Charges0003.png}
%      %\includegraphics[width=1\linewidth]{Majorant_Errors_vs_DOFs_LineWidth_1p5_Cropped.pdf}
%        \caption{{\small The harmonic component $\tilde u^H$  on the surface of the two chromophores in units $k_BT/e$.}}
%        \label{run01_2_Poisson_Harmonic_3_Term_2DYES_Ground_Charges0003}
%    \end{minipage}%
%    \end{figure}

\subsubsection{Finding the regular component $\tilde{u}$}\label{subsec:finding_regular_comp}

Now, we find a conforming approximation $\tilde u_h\in V_h$ of $\tilde u$, the exact solution of problem \eqref{tildeu_weak_formulation}, by adaptively solving the Galerkin problem \eqref{Discrete_Weak_Formulation_uRegular_3_term}, where $V_h\subset H_0^1(\Omega)$ is the $P_1$ Lagrange finite element space . By $\tilde u_i$ and $\tilde y^*_i$ we denote the finite element approximations at mesh refinement level $i$ to $\tilde u$ and $\tilde p^*$, respectively. Here $\bar p=6$. To find a good approximation $\tilde y^*$  of the exact flux $\tilde p^*$ of the form \eqref{special_form_yStar_uRegular}, i.e., $\tilde y^*$ such that 
\begin{align*}
\tilde y^*=-\epsilon_m T(\nabla\tilde u^H) \mathbbm{1}_{\Omega_m}+\epsilon_m\nabla G\mathbbm{1}_{\Omega_s}+\tilde y_0^*, \text{ for } \tilde y_0^*\in H(\div;\Omega) \text{ with } \div \tilde y_0^*=0 \text{ in } \Omega_m,
\end{align*}
we can minimize the majorant $M^2_\oplus(\tilde u,\tilde y^*)$ or, equivalently, minimize the functional $-I^*(\tilde y^*)=G^*(\tilde y^*)+F^*(-\Lambda^*\tilde y^*)$ in $\tilde y_0^*$ over a subspace of $H(\div;\Omega)$, like $RT0, RT1$, by additionally enforcing the condition that $\div \tilde y_0^* = 0$ in $\Omega_m$. Another, computationally favourable, approach, which can be easily realized in parallel, is to apply a patchwise flux reconstruction that will also yield $\div \tilde y_0^*=0$ in $\Omega_m$. Such an approach is motivated by observing that the exact $\tilde p_0^*:=\tilde p^*+\epsilon_m T(\nabla\tilde u^H) \mathbbm{1}_{\Omega_m}-\epsilon_m\nabla G\mathbbm{1}_{\Omega_s}$ satisfies the integral identity
\begin{align}
&\int\limits_{\Omega}{\tilde p_0^*\cdot\nabla v dx}=-\int\limits_{\Omega}{k^2\sinh(\tilde u) v dx},\text{ for all } v\in H_0^1(\Omega).
\end{align}
If we define the function $q:=\epsilon\nabla\tilde u_h+\epsilon_m T\left(\nabla\tilde u^H\right) \mathbbm{1}_{\Omega_m}-\epsilon_m\nabla G\mathbbm{1}_{\Omega_s}$, then we can see that 
\begin{align*}
&\int\limits_{\Omega}{q\cdot\nabla v dx}=-\int\limits_{\Omega}{k^2\sinh(\tilde u_h) v dx},\text{ for all } v\in V_h.
\end{align*}
Now, we define $\Pi_{L_h}(q)$ to be its $L^2$ projection over the space $L_h$ of piecewise constant functions over the same mesh on which $V_h$ is defined. 
%Then, since 
%\begin{align*}
%&-\int\limits_{\Omega_m}{\epsilon_m T\left(\nabla\tilde u^H\right)\cdot\nabla v dx}+\int\limits_{\Omega_s}{\epsilon_m\nabla G\cdot \nabla v dx} =-\int\limits_{\Omega_m}{\Pi_{L_h}\left(\epsilon_m T\left(\nabla\tilde u^H\right)\right)\cdot\nabla v dx}+\int\limits_{\Omega_s}{\Pi_{L_h}\left(\epsilon_m\nabla G\right)\cdot\nabla v dx} ,
%\end{align*}
Since $\Pi_{L_h}\left(q\right)$ satisfies the problem
\begin{align}
&\int\limits_{\Omega}{\Pi_{L_h}(q) \cdot\nabla v dx}=-\int\limits_{\Omega}{k^2\sinh(\tilde u_h) v dx},\text{ for all } v\in V_h,
\end{align}
we define $\tilde y_0^*\in RT_0$ by applying the patchwise equilibrated flux reconstruction in \cite{Braess_Schoberl_2006} to the numerical flux $\Pi_{L_h}(q)\in L_h$ . Notice that since $k=0$ in $\Omega_m$, the obtained $\tilde y_0^*$ satisfies the realtion $\div\tilde y_0^*=0$ in $\Omega_m$.
We define in a similar fashion, as for the component $\tilde{u}^N$, the quantities ${\text{PRE}}_{i,i}^{\text{Up}}$, ${\text{RE}}_{i,j,k,s}^{\text{Up}}$, ${\text{RCEN}}_{i,j,s}^{\text{Low}}$, ${\text{P}}_{\text{rel},i}^{\text{CEN}}$, ${\text{RCEN}}_{i,j,s}^{\text{Up}}$.
{\small
\begin{center}
\captionof{table}{Example 3} \label{Example3_Table1} 
\vspace{1ex}
%\begin{tabular}{ |p{1.9cm}|p{1.8cm}|p{2.cm}|p{2.1cm}|p{2cm}|p{2cm}|p{2cm}| }
\begin{tabular}{ |c|l|c|c|c|c|c| }
\hline
\multicolumn{7}{|c|}{Example 3: $k_m^2=0,\,k_s^2=10,\,\epsilon_m=2,\,\epsilon_s=80$} \\
\hline
 & & & & & &\\
level $i$ &\#elements &$\|\tilde u_i\|_0$ &$\vertiii{\nabla \tilde u_i}$ & $\vertiii{\epsilon\nabla \tilde u_i-\tilde y_{i}^*}_*$ &
$\sqrt{2}M_{\oplus}(\tilde u_i,\tilde y_{i}^*)$& $\tilde J(\tilde u_i)$\\
 & & & & & &\\
\hline
0	&667 008	 &67.614      &361.00 		&200.91		&201.89		&258854115.518 \\
1	&1 384 294 	 &68.342      &363.74 	        &127.75		&128.44		&258852218.914  \\
2	&2 162 155 	 &68.445      &364.54 		&108.08		&108.72		&258851876.473 \\
3	&2 889 668    &68.512      &365.00 		&97.878		&98.456		&258851659.119\\
4	&3 591 411  	 &68.554      &365.26 		&91.161		&91.700		&258851527.684\\
5	&4 322 074  	 &68.583      &365.44 		&86.080		&86.587		&258851434.952\\
6	&5 142 768   &68.607      &365.58	 	&81.773		&82.252	        &258851360.368\\\
7	&6 108 164   &68.626      &365.70	 	&77.910		&78.362	        &258851296.239\\
8	&7 276 168   &68.643      &365.81	 	&74.263		&74.690        &258851239.273\\
9	&8 690 783   &68.658      &365.91	 	&71.903		&72.299	        &258851187.598\\
10	&10 390 012  &68.671      &366.01	 	&67.538		&67.917	        &258851140.495\\
11	&12 431 042  &68.683      &366.09	 	&64.457		&64.814	        &258851097.844\\
12	&14 871 214  &68.695      &366.17	 	&61.558		&61.893        &258851059.718\\
\hline
\end{tabular}
\end{center}
}

{\small
\begin{center}
\captionof{table}{Example 3} \label{Example3_Table2} 
\vspace{1ex}
%\begin{tabular}{ |p{1.9cm}|p{1.8cm}|p{2.cm}|p{2.1cm}|p{2cm}|p{2cm}|p{2cm}| }
\begin{tabular}{ |c|l|c|c|c|c|c| }
\hline
\multicolumn{7}{|c|}{Example 3: $k_m^2=0,\,k_s^2=10,\,\epsilon_m=2,\,\epsilon_s=80$} \\
\hline
 & & & & & &\\
level $i$ &\#elements &${\text{PRE}}_{i,i}^{\text{Up}}[\%]$ & ${\text{RE}}_{i,i,i,i}^{\text{Up}}[\%]$ & ${\text{RCEN}}_{i,i,i}^{\text{Low}}[\%]$ & ${\text{P}}_{\text{rel},i}^{\text{CEN}}[\%]$& ${\text{RCEN}}_{i,i,i}^{\text{Up}}[\%]$\\
 & & & & & &\\
\hline
0	&667 008	 &55.65      &126.8		&18.91		&39.35		&58.10 \\
1	&1 384 294 	 &35.12      &54.59	        &13.71		&24.83		&31.95 \\
2	&2 162 155 	 &29.65      &42.50 		&12.02		&20.96		&25.99 \\
3	&2 889 668    &26.81      &36.93		&11.09		&18.96		&23.05 \\
4	&3 591 411    &24.95      &33.52		&10.45		&17.64		&21.17\\
5	&4 322 074    &23.55      &31.05 		&9.967		&16.65		&19.79\\
6	&5 142 768    &22.36      &29.03 		&9.544		&15.81	        &18.63\\
7	&6 108 164    &21.30      &27.27 		&9.158		&15.06	        &17.61\\
8	&7 276 168  	 &20.30      &25.65 		&8.789		&14.35	        &16.66\\
9	&8 690 783  	 &19.65      &24.62 		&8.547		&13.89	        &16.05\\
10	&10 390 012   &18.45      &22.78 		&8.094		&13.04	        &14.95\\
11	&12 431 042   &17.60      &21.51 		&7.769		&12.44	        &14.18\\
12	&14 871 214   &16.81      &20.34 		&7.460		&11.88	        &13.46\\
\hline
\end{tabular}
\end{center}
}

{\small
\begin{center}
\captionof{table}{Example 3} \label{Example3_Table3} 
\vspace{1ex}
%\begin{tabular}{ |p{1.9cm}|p{1.8cm}|p{2.cm}|p{2.1cm}|p{2cm}|p{2cm}|p{2cm}| }
\begin{tabular}{ |p{0.97cm}|p{1.68cm}|p{1.8cm}|p{2.1cm}|p{2.1cm}|p{2.2cm}|p{2.15cm}| }
%\begin{tabular}{ |c|l|c|c|c|c|c| }
\hline
\multicolumn{7}{|c|}{Example 3: $k_m^2=0,\,k_s^2=10,\,\epsilon_m=2,\,\epsilon_s=80$} \\
\hline
 & & & & & &\\
level $i$ &\#elements & ${\text{RE}}_{i,\bar p,\bar p,\bar p}^{\text{Up}}[\%]$ & ${\text{RCEN}}_{i,\bar p,\bar p}^{\text{Low}}[\%]$
& ${\text{RCEN}}_{i,\bar p,\bar p}^{\text{Up}}[\%]$ & $\vertiii{\epsilon\nabla \tilde u_i-\tilde y_{\bar p}^*}_*$  & $\sqrt{2}M_{\oplus}(\tilde u_i,\tilde y_{\bar p}^*)$ \\
% & &$[\%]$ &$[\%]$ &$[\%]$ & &\\
\hline
0	&667 008	 &32.75      &24.34		&43.92		&99.184		&99.665 \\
1	&1 384 294 	 &25.75      &15.48	        &27.94		&78.070		&78.374 \\
2	&2 162 155 	 &24.27      &13.09 		&23.65		&73.546		&73.857 \\
3	&2 889 668    &23.28      &11.86		&21.41		&70.520		&70.840 \\
4	&3 591 411    &22.66     &11.04		&19.94		&68.630		&68.957\\
5	&4 322 074    &22.21      &10.43 		&18.83		&67.274		&67.607\\
6	&5 142 768    &21.85      &9.910 		&17.89		&66.170	        &66.507\\
7	&6 108 164    &21.54      &9.442 		&17.04		&65.209	        &65.550\\
8	&7 276 168  	 &21.26      &9.000 		&16.24		&64.348	        &64.692\\
9	&8 690 783  	 &21.00      &8.714 		&15.72		&63.564	        &63.910\\
10	&10 390 012   &20.76      &8.185 		&14.77		&62.843	        &63.191\\
11	&12 431 042   &20.54      &7.811		&14.10		&62.179	        &62.525\\
12	&14 871 214   &20.34      &7.460 		&13.46		&61.558	        &61.893\\
\hline
\end{tabular}
\end{center}
}
We should note that the bounds on the error in energy norm obtained by the majorant $\sqrt{2}M_\oplus(\tilde u_i,\tilde y_i^*)$ in Table \ref{Example3_Table1} are rather conservative and they could be improved by applying a flux reconstruction involving a higher order Raviart-Thomas spaces, like $RT_1$. To obtain an idea of how much the error is overestimated, we can compare the values $\sqrt{2}M_\oplus(\tilde u_0,\tilde y^*_{\bar p})$ and $\sqrt{2}M_\oplus(\tilde u_1,\tilde y^*_{\bar p})$ of the majorant, evaluated with the last available approximation~$\tilde y^*_{\bar p}$ of the exact flux~$\tilde p^*$ to the values $\sqrt{2}M_\oplus(\tilde u_0,\tilde y^*_0)$ and $\sqrt{2}M_\oplus(\tilde u_1,\tilde y^*_1)$ evaluated with the current $\tilde y^*_0$ and $\tilde y^*_1$ at the first two meshes. We see that the overestimation is between around $201.89/99.665\approx 2.025$ and $128.44/78.374\approx 1.64$ times. This means that it is safe to assume that the real error $\vertiii{\nabla(\tilde u-\tilde u_{12})}$ at the last level $\bar p=12$ is no more than approximately $\sqrt{2}M_{\oplus}(\tilde u_{\bar p},\tilde y^*_{\bar p})/2\approx 31$. With this in mind, we can obtain an overall guaranteed bound on the error in energy norm for the regular component $u$ by using \eqref{overal_error_estiamte_for_uRegular_3_term_no_splitting_in_uLplusuN}:
\begin{align*}
&\mathrel{\phantom{=}}\vertiii{\nabla(u-\tilde u_h)}\leq \sqrt{\epsilon_m}M_{\oplus,H}\left(\tilde u^H,T(\nabla \tilde u^H)\right)+ \sqrt{2}M_\oplus(\tilde u_h,\tilde y^*)\\
&=\sqrt{2}\times 43.085 + 61.893 = 122.824
\end{align*}
%$40.52\,\%$ 
\begin{figure}[!htb]
    \centering
    \begin{minipage}{1\textwidth}
        \centering
        \captionsetup{width=0.8\linewidth}
      \includegraphics[width=0.8\linewidth]{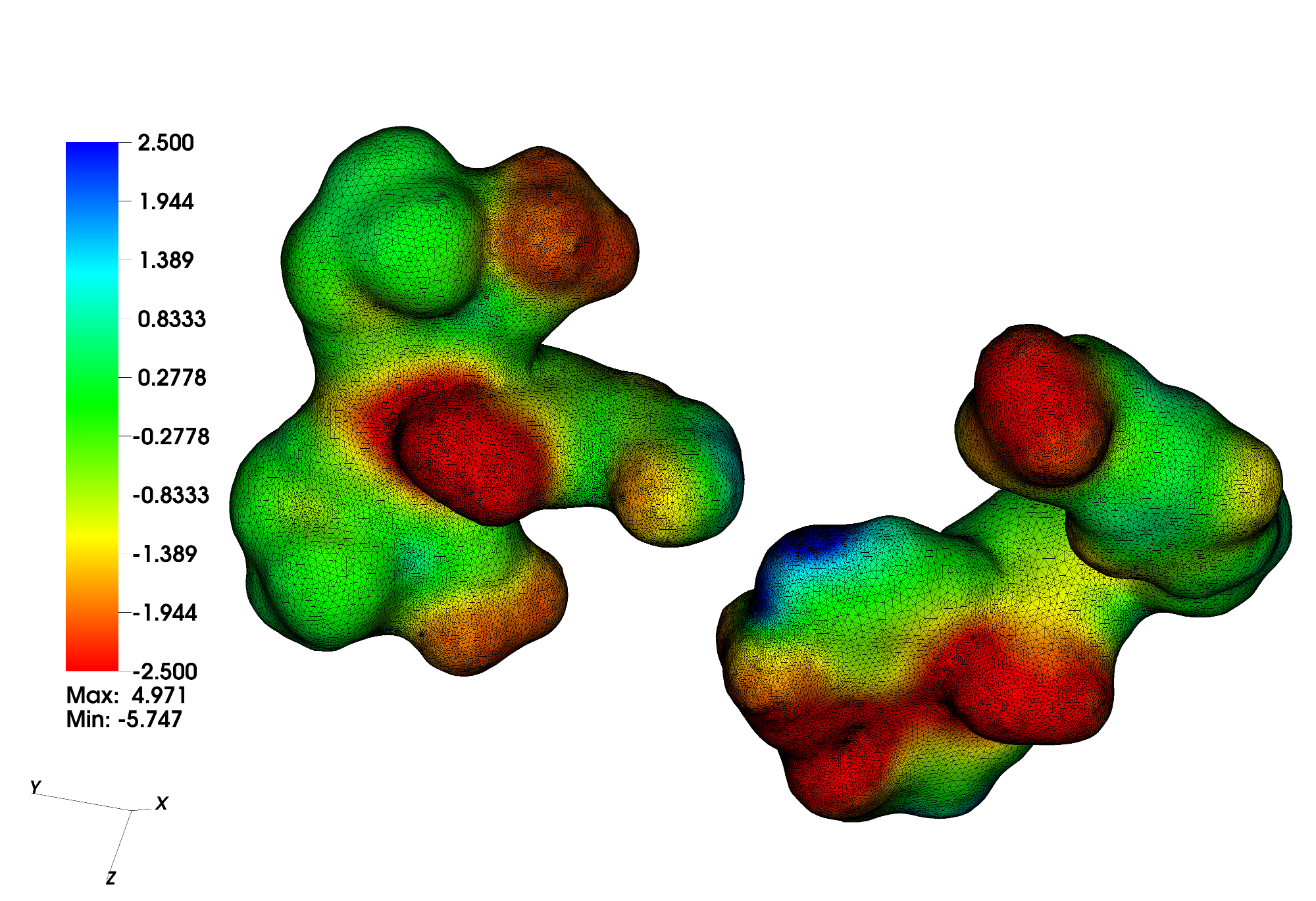}
        \caption{{\small Full potential surface map with the 3-term regularization (without additional splitting into $u^L+u^N$) for the system Alexa 488 and Alexa 594  in units $k_BT/e_c$. Blue color indicates a positive potential (values$>2.5 k_BT/e_c$) and red color indicates negative potential (values $<-2.5 K_BT/e_c$).}}
        \label{run01_6_3_Term_Regular_2DYES_epsm2_epss80_ks_squared_100000}
    \end{minipage}%
    \end{figure}
    
    \begin{figure}[!htb]
    \centering
    \begin{minipage}{1\textwidth}
        \centering
        \captionsetup{width=0.8\linewidth}
      \includegraphics[width=0.8\linewidth]{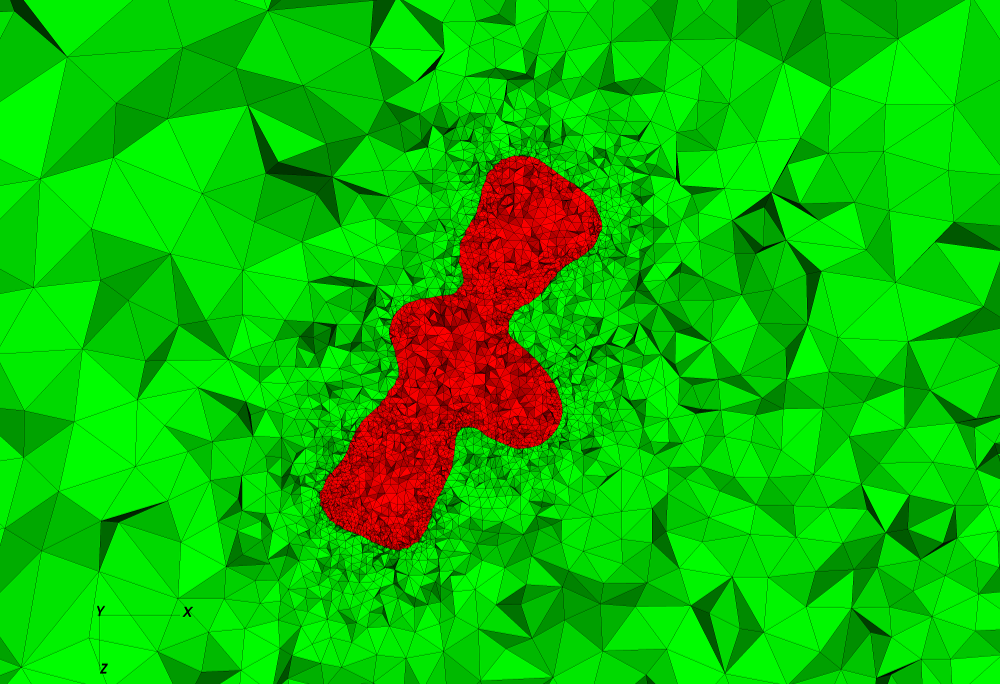}
        \caption{{\small Cross section of the mesh with the plane $y=3\, \AA$ at level 1 in the mesh refinement procedure for finding the component $\tilde u$ in Example 3. The molecule region $\Omega_m$ is marked red (Alexa 594).}}
        \label{run07_2DYES_3_Term_REGULAR_epsm2_epss80_ks_squared_Regions_Level_1}
    \end{minipage}%
    \end{figure}
    
    \begin{figure}[!htb]
    \centering
    \begin{minipage}{1\textwidth}
        \centering
        \captionsetup{width=0.8\linewidth}
      \includegraphics[width=0.8\linewidth]{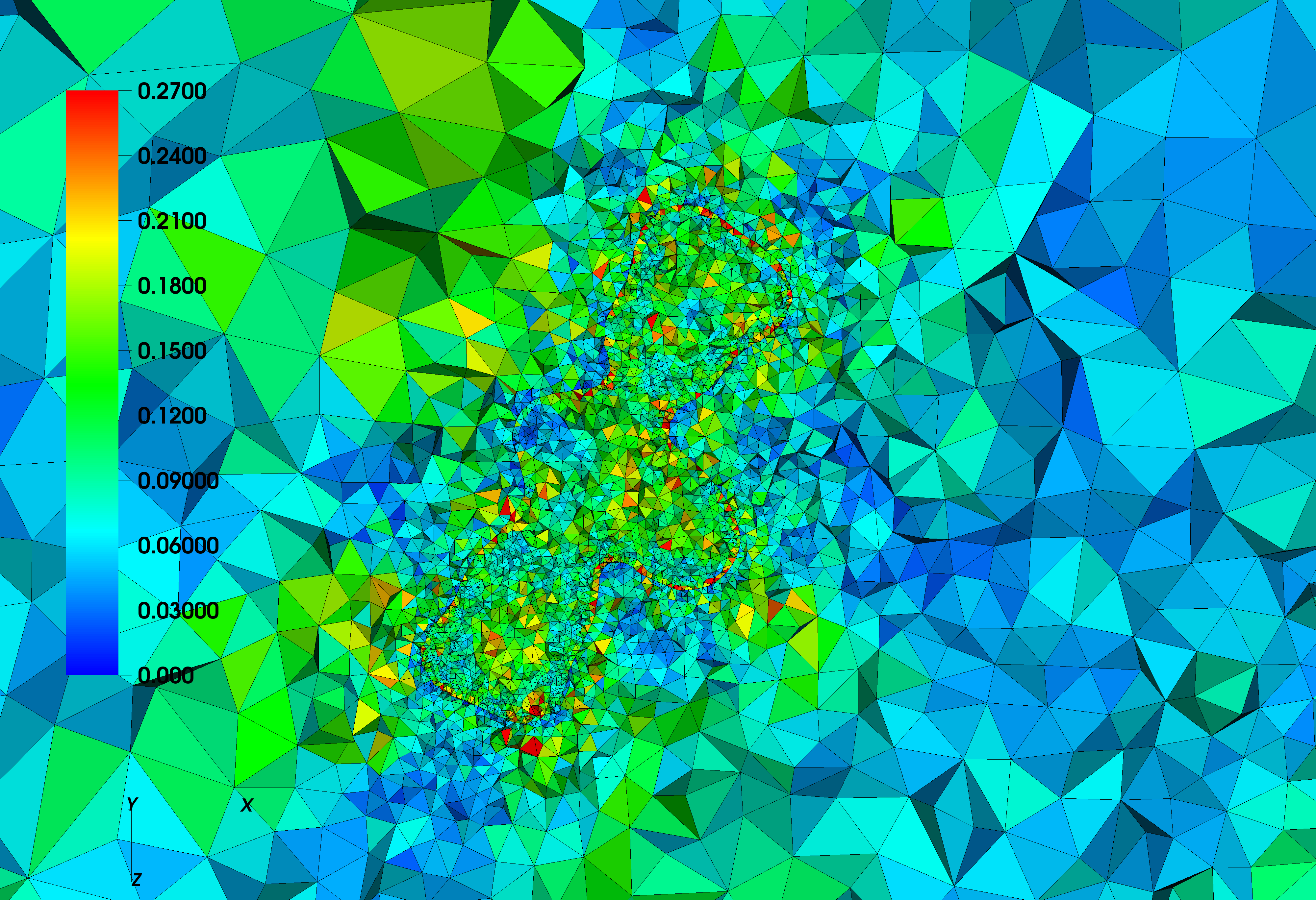}
        \caption{{\small Cross section of the mesh with the plane $y=3\, \AA$ at level 1 in the mesh refinement procedure for finding the component $\tilde u$ in Example 3. Error indicator as a piecewise constant function.}}
        \label{run07_2DYES_3_Term_REGULAR_epsm2_epss80_ks_squared_Func_Err_Indicator_Level_1}
    \end{minipage}%
    \end{figure}
    
\section{Conclusions}

We have analyzed the well posedness of two and three-term regularization schemes for the Poisson-Boltzmann equation and derived guaranteed and fully computable bounds on the error in energy
norm. For the 2-term regularization the dimensionless potential  $\tilde \phi$ is decomposed into $G+u$ where $G$ is analytically known and $u\in H^1(\Omega)$ is the regular component which has to
be approximated, that is, computed numerically. The regular component $u$ can be split additionally into $u^L+u^N$ where $u^L$ solves a linear nonhomogeneous interface problem and $u^N$ solves
a nonlinear homogeneous problem which depends on~$u^L$. For each of these two problems we have derived guaranted bounds on the error in energy norm. Moreover, for the nonlinear problem, we
have proved a continous dependence of the solution $u^N$ on perturbations in $u^L$. This property has been exploited to estimate the overall error in the regular component $u$. The derived error
estimate for $u$ is a linear combination of the majorants for the error in $u^L$ and the error in $u^N$ with perturbed $u^N$. Similarly, in the 3-term regularization scheme, the dimensionless potential
$\tilde \phi$ is decomposed into $G+u^H+u$. Here $u^H\in H^1(\Omega)$ is a harmonic function in the molecular domain $\Omega_m$, which has to be approximated numerically, and is equal to $-G$
in the solution domain~$\Omega_s$. Now, the regular component $u$ satisfies a nonlinear nonhomogeneous interface problem which depends on $u^H$. To solve this problem, we have analyzed two approaches. The first is to additionally make the splitting $u=u^L+u^N$ as we did in case of the 2-term regularization. In this case, we have analyzed how $u^L$ depends on perturbations in $u^H$,
and further, how $u^N$ depends on perturbations in $u^L$. Finally, we have derived an estimate for the overall error in the regular component $u$, which is a linear combination of the majorants for the
error in each of the components $u^H,u^L$, and $u^N$.
The second approach to derive estimates for the regular component $u$ of the solution of the nonlinear interface problem is to directly derive an error estimate for it and analyze its continuous
dependence on perturbations in $u^H$. In this case, the overall estimate for the error in energy norm is a linear combination of the majorants for the error in $u^H$ and the error in $u$ with
perturbed~$u^H$.

The a posteriori error analysis presented in this paper is based on the functional approach developed in~\cite{Repin_2000}. We have also utilized this approach to obtain a near best approximation
result for the two regularization schemes which is the basis for the analysis of qualified and unqualified convergence of finite element approximations. In other words, the generality of this method
allows for the derivation of both a posteriori and a priori error estimates for the considered class of problems.

We have presented three numerical tests performed on two realistic physical systems illustrating and validating our theoretical findings. The first system consists of the two chromophores, Alexa 488 and Alexa 594, and the second system of an insulin protein with a PDB ID 1RWE. The guaranteed error bounds which we have derived do not overestimate the error as drastically as residual based error estimates often do. To obtain a conforming approximation of the dual variable, we have utilized a patchwise flux reconstruction technique, cf.~\cite{Braess_Schoberl_2006}, which can be easily implemented on parallel machines by a scalable algorithm with linear complexity. This means that the proposed error estimation can be realized in a very efficient manner.

\bibliographystyle{plain}
\bibliography{references}
\end{document}